\documentclass[11pt]{amsart}
\usepackage{amsmath} \usepackage{amsxtra}
\usepackage{amscd} \usepackage{amsthm}  
\usepackage{amsfonts}\usepackage{amssymb} 
\usepackage{eucal}\usepackage{bm}
\usepackage{graphicx} 

\newtheorem{Cor}[subsubsection]{Corollary}
\newtheorem{Lm}[subsubsection]{Lemma}
\newtheorem{Pp}[subsubsection]{Proposition}
\newtheorem{Con}[subsubsection]{Conjecture}
\newtheorem{Thm}[subsubsection]{Theorem}
\newtheorem{Def}[subsubsection]{Definition}
\newtheorem{Rem}[subsubsection]{Remark}

\newcommand{\cA}{{\mathcal A}}
\newcommand{\cB}{{\mathcal B}}
\newcommand{\cC}{{\mathcal C}}
\newcommand{\cD}{{\mathcal D}}
\newcommand{\cH}{{\mathcal H}}
\newcommand{\cE}{{\mathcal E}}
\newcommand{\cG}{{\mathcal G}}

\newcommand{\cJ}{{\mathcal J}}
\newcommand{\cO}{{\mathcal O}}
\newcommand{\cL}{{\mathcal L}}
\newcommand{\cM}{{\mathcal M}}

\newcommand{\cF}{{\mathcal F}}
\newcommand{\cK}{{\mathcal K}}

\newcommand{\cQ}{{\mathcal Q}}

\newcommand{\cS}{{\mathcal S}}

\newcommand{\cU}{{\mathcal U}}
\newcommand{\cV}{{\mathcal V}}
\newcommand{\cW}{{\mathcal W}}
\newcommand{\cX}{{\mathcal X}}
\newcommand{\cY}{{\mathcal Y}}
\newcommand{\cZ}{{\mathcal Z}}

\renewcommand{\AA}{{\mathbb A}}

\newcommand{\ZZ}{{\mathbb Z}}

\newcommand{\PP}{{\mathbb P}}

\newcommand{\gp}{\mathfrak{p}}
\newcommand{\gq}{\mathfrak{q}}

\newcommand{\gr}{\mathfrak{r}}
\newcommand{\gs}{\mathfrak{s}}

\newcommand{\on}{\operatorname}
\newcommand{\Sh}{\on{Sh}}

\newcommand{\RCov}{\on{RCov}}
\newcommand{\mult}{\on{mult}}

\newcommand{\Rep}{{\on{Rep}}}

\newcommand{\Qlb}{\mathbb{\bar Q}_\ell}
\newcommand{\Gm}{\mathbb{G}_m}

\newcommand{\A}{\mathbb{A}}

\newcommand{\toup}[1]{\stackrel{#1}{\to}}
\newcommand{\hook}[1]{\stackrel{#1}{\hookrightarrow}}

\newcommand{\getsup}[1]{\stackrel{#1}{\gets}}

\newcommand{\Sp}{\on{\mathbb{S}p}}
\newcommand{\Spin}{\on{\mathbb{S}pin}}

\newcommand{\GSp}{\on{G\mathbb{S}p}}
\newcommand{\IC}{\on{IC}}

\newcommand{\Hom}{\on{Hom}}
\newcommand{\Mod}{\on{Mod}}

\newcommand{\End}{\on{End}}

\newcommand{\Sym}{\on{Sym}}
\newcommand{\SO}{\on{S\mathbb{O}}}
\newcommand{\GO}{\on{G\mathbb{O}}}
\newcommand{\Ker}{\on{Ker}}

\newcommand{\Aut}{\on{Aut}}

\newcommand{\RG}{\on{R\Gamma}}

\newcommand{\sym}{\on{sym}}
\newcommand{\triv}{\on{triv}}

\newcommand{\Pic}{\on{Pic}}
\newcommand{\uPic}{\on{\underline{Pic}}}
\newcommand{\Bun}{\on{Bun}}

\newcommand{\Bunb}{\on{\overline{Bun}} }
\newcommand{\Bunt}{\on{\widetilde\Bun}}

\newcommand{\Spec}{\on{Spec}}

\newcommand{\supp}{\on{supp}}

\newcommand{\END}{{{\mathcal E}nd}}

\newcommand{\GL}{\on{GL}}

\newcommand{\PGL}{\on{PGL}}

\newcommand{\Eis}{{\on{Eis}}}
\newcommand{\pr}{\on{pr}}
\newcommand{\id}{\on{id}}
\newcommand{\tr}{\on{tr}}
\newcommand{\QED}{$\square$} 
\newcommand{\Fq}{\mathbb{F}_q}  
\newcommand{\Fp}{\mathbb{F}_p}  

\newcommand{\iso}{{\widetilde\to}}

\newcommand{\comp}{\circ}
\newcommand{\Four}{\on{Four}}
\renewcommand{\H}{{\on{H}}}   

\newcommand{\DD}{\mathbb{D}}  
\newcommand{\D}{\on{D}}       
\newcommand{\wt}{\widetilde}

\newcommand{\select}[1]{{\it{#1}}}

\newcommand{\und}[1]{\underline{#1}}

\renewcommand{\div}{\on{div}}

\newcommand{\Av}{\on{Av}}
\newcommand{\ev}{\mathit{ev}}

\newcommand{\act}{\on{act}}
\newcommand{\dimrel}{\on{dim.rel}}

\newcommand{\SL}{\on{SL}}

\newcommand{\WP}{\on{WP}} 

\newcommand{\Modt}{\wt\Mod}
\newcommand{\ItMod}{\on{ItMod}}
\newcommand{\Int}{\on{IntMod}}
\newcommand{\La}{{\on{La}}}
\newcommand{\diag}{\on{diag}}

\newenvironment{Prf}{\par\noindent {\it Proof }}{\QED}
\newcommand{\Step}[1]{\par\noindent{\bf Step {#1}}.}

\newcommand{\nc}{\newcommand}
\nc{\ssec}{\subsection}
\nc{\sssec}{\subsubsection}

\begin{document}

\author{Sergey Lysenko}
\title{Geometric Waldspurger periods II}
\address{Institut Elie Cartan Lorraine, Universit\'e de Lorraine, 
 B.P. 239, F-54506 Vandoeuvre-l\`es-Nancy Cedex, France}
\email{Sergey.Lysenko@univ-lorraine.fr}
\begin{abstract} In this paper we extend the calculation of the geometric Waldspurger periods from our paper \cite{Ly1} to the case of ramified coverings. We give some applications to the study of Whittaker coeffcients of the theta-lifting of automorphic sheaves from $\PGL_2$ to the metaplectic group $\wt\SL_2$, they agree with our conjectures from \cite{Ly6}. In the process of the proof, we construct some new automorphic sheaves for $\GL_2$ in the ramified setting. We also formulate stronger conjectures about Waldspurger periods and geometric theta-lifting for the dual pair $(\wt\SL_2, \PGL_2)$.
\end{abstract} 
\thanks{The author would like to thank V. Lafforgue for constant and stimulating discussions. The author was supported by the ANR program ANR-13-BS01-0001-01.} 
\maketitle

\section{Introduction and main results}

\sssec{} In the present paper we prove some conjectures from the first part of this paper \cite{Ly1}. Let $X$ and $Y$ be smooth projective connected curves and $\phi: Y\to X$ a degree 2 ramified covering. Let $E$ be an irreducible rank two $\Qlb$-local system on $X$. Denote by $\Aut_E$ be the corresponding automorphic sheaf on the stack $\Bun_2$ of rank 2 vector bundles on $X$. We calculate the tensor (as well as symmetric and exterior) square of the Waldspurger periods of $\Aut_E$ for this covering, and so establish \cite[Conjecture~3]{Ly1} for irreducible rank 2 local systems on $X$. This generalizes \cite[Theorem~5]{Ly1} to the case of ramified coverings. 

 In the process of the proof, we get some new automorphic sheaves on $\GL_2$ in the ramified setting (not available in the literature to the best of our knowledge).

 We apply these results to the study of the geometric theta-lifting for the dual pair $(\wt\SL_2, \PGL_2)$. This application is our main motivation. So, this paper is mostly about the quantum geometric Langlands program for $\wt\SL_2$. Namely, consider the theta-lifting functor $F_G: \D^-(\Bun_{\PGL_2})_!\to \D^{\prec}(\Bunt_G)$ (cf. Section~\ref{Sect_1.2.2}). Here $\Bunt_G$ plays the role of the stack of $\wt\SL_2$-torsors on $X$. For an irreducible $\SL_2$-local system $E$ on $X$, let $\Aut_E$ be the automorphic sheaf on $\Bun_{\PGL_2}$ normalized as in \cite{FGV}. Our main result is Corollary~\ref{Cor_main_for_SL_2_SO_3}. It gives some description of (the tensor square of) the most non-degenerate Whittaker coefficients of $F_G(\Aut_E)$. It turns out that this Whittaker coefficient is placed in one cohomological degree, so is a vector space. We describe the tensor square of this vector space directly in terms of $E$. 
 
 This answer is to be compared with the conjectural extension of the Casselman-Shalika formula for the metaplectic Langlands program suggested in \cite[Section~11]{Ly7}, see also \cite{GL1, GL2}. 
  
\sssec{} We first formulate our results about the Waldspurger periods. Then we present our motivations, results and conjectures about $\wt\SL_2$, and finally we explain our strategy of the proof and some intermediate results that may be of independent interest.

\sssec{} Assume for a moment the ground field $k=\Fq$ finite of $q$ elements with $q$ odd. Let $G=\GL_2$. Let $E$ be an irreducible rank 2 local system on $X$, $f_E: \Bun_G(k)\to\Qlb$ be the function `trace of Frobenius' of $\Aut_E$. Let $\phi: Y\to X$ be a degree 2 possibly ramified covering, where $Y$ is smooth and projective. Write $\Pic Y$ for the Picard stack of $Y$. Let $\cJ$ be a rank one local system on $Y$ equipped with an isomorphism $N(\cJ)\,\iso\, \det E$, here $N(\cJ)$ is the norm of $\cJ$. Write $f_{\cJ}: (\Pic Y)(k)\to\Qlb$ for the corresponding character, the trace of Frobenius of the automorphic local system $A\cJ$ on $\Pic Y$ corresponding to $\cJ$. The Waldspurger period of $f_E$ is
$$
\int_{\cB\in (\Pic Y)(k)/(\Pic X)(k)} f_E(\phi_*\cB)f_{\cJ}^{-1}(\cB)d\cB\, .
$$
The function that we integrate does not change when $\cB$ is tensored by $\phi^*L$, $L\in\Pic X$, here $d\cB$ is a Haar measure. The theorem of Waldspurger says that the square of this period is equal (up to some explicit coefficient) to the central value of the $L$-function $L(\phi^*E\otimes\cJ^{-1}, \frac{1}{2})$, see \cite{W}. This result played a central role in the classification of the automorphic discrete spectrum of $\wt\SL_2$ by Waldspurger in \cite{W1, W2} (see also \cite{WTG} for a recent survey). We prove a geometric version of this result (cf. Theorem~\ref{Th_3}).

\sssec{General notation} Let $k$ denote an algebraically closed field of characteristic $p>2$, all the schemes (or stacks) we consider are defined over $k$. Fix a prime $\ell\ne p$. 

 For a stack $S$ locally of finite type write $\D(S)$ for the category introduced in \cite[Remark~3.21]{LO} and denoted $\D_c(S,\Qlb)$ in \select{loc.cit}. It should be thought of as the unbounded derived category of constructible $\Qlb$-sheaves on $S$. For $\ast=+,-, b$ we have the full triangulated subcategory $\D^{\ast}(S)\subset \D(S)$ denoted $\D_c^{\ast}(S,\Qlb)$ in \select{loc.cit.} Write $\D^{\ast}(S)_!\subset \D^{\ast}(S)$ for the full subcategory of objects which are extensions by zero from some open substack of finite type. Write $\D^{\prec}(S)\subset \D(S)$ for the full subcategory of complexes $K\in \D(S)$ such that for any open substack $U\subset S$ of finite type we have $K\mid_U\in \D^-(U)$. Write $\DD: \D^b(S)\to \D^b(S)$ for the Verdier duality functor.  

 Fix a nontrivial character $\psi: \Fp\to\Qlb^*$ and denote by
$\cL_{\psi}$ the corresponding Artin-Shreier sheaf on $\A^1$. 
Since we are working over an algebraically closed field, we systematically ignore the Tate twists. If $V\to S\gets V^*$ are dual rank $r$ vector bundles over a stack $S$, we normalize 
the Fourier transform $\Four_{\psi}: \D^{\prec}(V)\to\D^{\prec}(V^*)$ by 
$\Four_{\psi}(K)=(p_{V^*})_!(\xi^*\cL_{\psi}\otimes p_V^*K)[r]$,  
where $p_V, p_{V^*}$ are the projections, and $\xi: V\times_S V^*\to \A^1$ is the pairing.

 Let $X$ be a smooth projective connected curve of genus $g$. Write $\Omega$ for the canonical line bundle on $X$. For an algebraic group $G$ we denote by $\cF_G^0$ the trivial $G$-torsor on a given base. For a morphism of stacks $f:Y\to Z$ we denote by $\dimrel(f)$ the function of a connected component $C$ of $Y$ given by $\dim C-\dim C'$, where $C'$ is the connected component of $Z$ containing $f(C)$.
 
 Write $\Bun_k$ for the stack of rank $k$ vector bundles on $X$. For $k=1$ we also write $\Pic X$ for the Picard stack $\Bun_1$ of $X$. We have a line bundle $\cA_k$ on $\Bun_k$ with fibre $\det\RG(X,V)$ at $V\in\Bun_k$. View it as a $\ZZ/2\ZZ$-graded placed in  degree $\chi(V)\!\mod 2$.

\ssec{Geometric Waldspurger periods}  
\label{Sect_Geometric_Waldspurger_intro}

\sssec{} For $d\ge 0$ write $X^{(d)}$ for the $d$-th symmetric power of $X$, we view it as the scheme of degree $d$ effective divisors on $X$. Let $^{rss}X^{(d)}\subset X^{(d)}$ be the open subscheme of multiplicity free divisors. Pick $d\ge 0$, $D\in {^{rss}X^{(d)}}$. Let $\cE_{\phi}$ be a line bundle on $X$ equipped with $\cE_{\phi}^2\,\iso\, \cO(-D)$. Note that $\cO\oplus \cE_{\phi}$ is a $\cO_X$-algebra and set $Y=\Spec(\cO\oplus \cE_{\phi})$. Then $Y$ is a smooth projective curve, the projection $\phi: Y\to X$ is a degree two covering ramified exactly over $D$. Write $D_Y$ for the ramification divisor on $Y$, the restriction of $\phi$ gives an isomorphisms $D_Y\,\iso\, D$ of these finite $k$-schemes. Write $\sigma$ for the nontrivial automorphism of $Y$ over $X$. Write $\Omega_Y$ for the canonical line bundle on $Y$. 

Let $\Pic Y$ be the Picard stack of $Y$, $\phi_1: \Pic Y\to \Bun_2$ the map $\cB\mapsto \phi_*\cB$. We denote by $N: \Pic Y\to \Pic X$ the norm map given by $N(\cB)=\cE_{\phi}^{-1}\otimes \det(\phi_*\cB)$ for $\cB\in\Pic Y$. This is a homomorphism of group stacks.

\sssec{} Let $U_{\phi}$ be the group scheme on $X$ given by the cokernel of the natural map $\Gm\to \phi_*\Gm$. Let $\Bun_{U_{\phi}}$ be the stack of $U_{\phi}$-torsors on $X$. We denote by $e_{\phi}: \Pic Y\to \Bun_{U_{\phi}}$ the extension of scalars map. Recall from \cite{Ly1} that $\Bun_{U_{\phi}}$ can be seen as the stack classifying $\cB\in\Pic Y$ equipped with $N(\cB)\,\iso\, \cO_X$ and a compatible isomorphism $\gamma: \cB\mid_{D_Y}\,\iso\, \cO_{D_Y}$. This means that the square of $\gamma$ is the restriction $N(\cB)\mid_D\,\iso\, \cO_D$ of the above trivialization. Our convention is that under this identification $e_{\phi}$ sends $\cB$ to $\sigma^*\cB\otimes \cB^{-1}$ with the evident trivializations 
$$
N(\sigma^*\cB\otimes \cB^{-1})\,\iso\, \cO, \;\;\; \gamma: \sigma^*\cB\otimes \cB^{-1}\mid_{D_Y}\,\iso\, \cO_{D_Y}
$$
Note that $\Bun_{U_{\phi}}$ is naturally a group stack, we denote by $\mult: \Bun_{U_{\phi}}\times \Bun_{U_{\phi}}\to \Bun_{U_{\phi}}$ the multiplication map. By \cite[Appendix~A]{Ly1}, the connected components of $\Bun_{U_{\phi}}$ are indexed by $a\in\ZZ/2\ZZ$ and denoted $\Bun^a_{U_{\phi}}$. Here $\Bun^0_{U_{\phi}}$ is the connected component of unity. 

 If $d>0$ then $\Bun_{U_{\phi}}$ is a scheme and $\mult$ is proper. 
More precisely, write $\Pic^r Y$ for the component of $\Pic Y$ classifying line bundles of degree $r$, similarly for $X$. Let $\und{\Pic}^r Y$, $\und{\Pic}^r X$ be the corresponding coarse moduli spaces. For $d>0$ the natural map $\und{\Pic}^r X\to \und{\Pic}^{2r} Y$ is a closed immersion, and for any $r\in\ZZ$ one has $\und{\Pic}^{2r} Y/\und{\Pic}^r X\,\iso\, \Bun_{U_{\phi}}^0$. In this case if $\und{P}$ is the corresponding Prym variety then there is a Galois covering $\Bun^0_{U_{\phi}}\to \und{P}$ with Galois group isomorphic to $(\ZZ/2\ZZ)^{d-2}$, cf. \cite[Appendix~A]{Ly1} and \cite{Mu}.
 
\sssec{} 
\label{Sect_1.1.3_intro}
For a rank one local system $\cA$ on $X$ we denote by $A\cA$ the automorphic (character) local system on $\Pic X$ associated to $\cA$. For $r\ge 0$ we have the Abel-Jacobi map $X^{(r)}\to\Pic X$, $D_1\mapsto \cO(D_1)$, and the $*$-restriction of $A\cA$ under this map is canonically identified with the symmetric power $\cA^{(r)}$ of $\cA$. 
  
 For a notion of central character of $K\in \D(\Bun_2)$ we refer to \cite[Definition~3]{Ly1}. If $K\in \D(\Bun_2)$ has central character $\cA$, where $\cA$ is a rank one local system on $X$, then $\act^*K\,\iso\, A\cA\boxtimes K$ for the map $\act: \Bun_1\times\Bun_2\to \Bun_2$, $(L_1, L)\mapsto L\otimes L_1$. Recall the definition of the Waldspurger periods given in \cite[Definition~10]{Ly1}.
 
 If $\cJ$ is a rank one local system on $Y$ then $\cJ\otimes\sigma^*\cJ$ is equipped with natural descent data for $\phi: Y\to X$, and there is a rank one local system $N(\cJ)$ on $X$ equipped with $\phi^*N(\cJ)\,\iso\, \cJ\otimes \sigma^*\cJ$, the norm of $\cJ$. 
 
\begin{Def} 
\label{Def_WP}
Let $\cJ$ be a rank one local system on $Y$. Let $K\in\D(\Bun_2)$ be a complex with central character $N(\cJ)$. Then $A\cJ^{-1}\otimes\phi_1^*K$ is equipped with natural descent data for $e_{\phi}: \Pic Y\to \Bun_{U_{\phi}}$. Assume the following holds:
\begin{itemize}
\item[$(C_W)$] We are given $\cK_K\in \D(\Bun_{U_{\phi}})$ equipped with
$$
e_{\phi}^*\cK_K[\dimrel(e_{\phi})]\,\iso\, A\cJ^{-1}\otimes \phi_1^*K[\dimrel(\phi_1)]
$$
\end{itemize}
For $a\in\ZZ/2\ZZ$ the Waldspurger period of $K$ is
$$
WP^a(K,\cJ)=\RG_c(\Bun^a_{U_{\phi}}, \cK_K)
$$
\end{Def}

\sssec{}
\label{Sect_1.1.3_for_defs}
 For $r\ge 0$ let $m_{\phi, r}: Y^{(r)}\to \Bun_{U_{\phi}}$ be the map sending $D_1$ to $\cO(D_1-\sigma^*D_1)$ with the natural trivializations 
$$
N(\cO(D_1-\sigma^*D_1))\,\iso\, \cO_X,\;\;\;\; \cO(D_1-\sigma^*D_1)\mid_{D_Y}\,\iso\, \cO_{D_Y}
$$  
Let $E$ be a rank 2 irreducible local system on $X$. Write $\Aut_E$ for the corresponding automorphic sheaf on $\Bun_2$ normalized as in \cite{FGV}.
\begin{Thm} 
\label{Th_3}
Let $\cJ$ be a rank one local system on $Y$ equipped with $\det E\,\iso\, N(\cJ)$. The condition $(C_W)$ of Definition~\ref{Def_WP} is satisfied for $\Aut_E$, giving rise to $\cK:=\cK_{\Aut_E}\in \D(\Bun_{U_{\phi}})$. The complex $\cK$ is a direct sum of (possibly shifted) irreducible perverse sheaves on $\Bun_{U_{\phi}}$. One has
$$
\mult_!(\cK\boxtimes\cK)\,\iso\, \mathop{\oplus}\limits_{r\ge 0} \; (A\cJ^{-1})_{\Omega_Y}\otimes(m_{\phi, r})_!(\cJ\otimes\phi^*E^*)^{(r)}[r]
$$
Here $(A\cJ^{-1})_{\Omega_Y}$ denotes the $*$-fibre of $A\cJ^{-1}$ at $\Omega_Y\in\Pic Y$. In particular, for $a\in\ZZ/2\ZZ$ there is an isomorphism
$$
\mathop{\oplus}\limits_{
\substack{a_1+a_2=a, \\
a_i\in\ZZ/2\ZZ}}
 \!\WP^{a_1}(\Aut_E,\cJ)\otimes\WP^{a_2}(\Aut_E,\cJ)\,
\iso\, \mathop{\oplus}\limits_{\substack{r\ge 0, \\ a= r\!\!\!\!\!\mod \! 2}} \; (A\cJ^{-1})_{\Omega_Y}\otimes \RG(Y^{(r)}, (\cJ\otimes\phi^*E^*)^{(r)})[r]
$$ 
If $\phi^*E$ is irreducible then the latter complex is a vector space placed in cohomological degree zero.
\end{Thm} 

 This establishes \cite[Conjecture~3]{Ly1} for all irreducible rank 2 local systems on $X$. We also describe the symmetric and exteriour squares of Waldspurger periods in Remark~\ref{Rem_great}. Recall that for a local system $\cE$ on $X$ the geometric analog of the central $L$-value $L(\cE, \frac{1}{2})$ is the complex $\oplus_{r\ge 0} \RG(X^{(r)}, \cE^{(r)})[r](\frac{r}{2})$, so Theorem~\ref{Th_3} implies the Waldspurger formula by passing to the trace of Frobenius (up to the Tate twists that we ignored).

\begin{Rem} i) Let $\cJ$ be any rank one local system on $Y$. By \cite[Lemma~15]{Ly6}, one has canonically $\det\RG(Y, \cJ)\,\iso\, (A\cJ^{-1})_{\Omega_Y}$.

\smallskip\noindent
ii) In the situation of Theorem~\ref{Th_3} assume in addition that $
\H^0(Y, \cJ\otimes \phi^*E^*)=0$. Since $\sigma^*\cJ^{-1}\otimes \phi^*E\,\iso\, \cJ\otimes \phi^* E^*$, this yields $\H^2(Y, \cJ\otimes \phi^*E^*)=0$. Set $V=\H^1(Y, \cJ\otimes \phi^*E^*)$. The symplectic form $E\otimes E\to\det E$ induces a map 
\begin{equation}
\label{map_intro_on_H^2}
\H^2(Y, \cJ\otimes \sigma^*\cJ\otimes \phi^*E^*\otimes \phi^*E^*)\to \H^2(Y,\Qlb)\,\iso\, \Qlb
\end{equation}
Since the cup-product 
$$
V\otimes V\,\iso\, \H^1(Y, \cJ\otimes \phi^*E^*)\otimes \H^1(Y, \sigma^*\cJ\otimes \phi^*E^*)\to \H^2(Y, \cJ\otimes \sigma^*\cJ\otimes \phi^*E^*\otimes \phi^*E^*)
$$ 
is skew-symmetric, composing it with (\ref{map_intro_on_H^2}) one obtains a non-degenerate symmetric form $\Sym^2 V\to \Qlb$ on $V$. Let $\SO(V)$ be the special orthogonal subgroup of $\GL(V)$ preserving this form. We may view $V$ as a $\SO(V)$-torsor on $\Spec \Qlb$. One has $\dim V=4g_Y-4$, where $g_Y$ is the genus of $Y$. Let $\Spin(V)$ denote the simply-connected cover of $\SO(V)$. Let $\Gamma_{\alpha}, \Gamma_{\beta}$ be the half-spin representations of $\Spin(V)$ over $\Qlb$, here $\alpha, \beta$ are the correspnding fundamental weights of $\Spin(V)$ \cite[Section19.2, p. 291]{FH91}. Then 
$$
\Gamma_{\alpha}\otimes\Gamma_{\alpha}\oplus \Gamma_{\beta}\otimes\Gamma_{\beta}\,\iso\, \wedge^0 V\oplus \wedge^2 V\oplus \wedge^4 V\oplus\ldots
$$
and
$$
\Gamma_{\alpha}\otimes\Gamma_{\beta}\oplus \Gamma_{\beta}\otimes\Gamma_{\alpha}\,\iso\, \wedge^1 V\oplus \wedge^3 V\oplus \wedge^5 V\oplus\ldots
$$
Pick a trivialization of $(A\cJ)_{\Omega_Y}$. Then we see that there should be a numbering $\alpha_a$, $a\in\ZZ/2\ZZ$ of the half-spin fundamental weights of $\Spin V$ and isomorphisms 
$$
\WP^a(\Aut_E,\cJ)\,\iso\, \Gamma_{\alpha_a},
$$
where $\Gamma_{\alpha_a}$ is the irreducible half-spin representation of $\Spin(V)$ with highest weight $\alpha_a$. Viewing accordingly $\Gamma_{\alpha}\oplus\Gamma_{\beta}$ as the a $\ZZ/2\ZZ$-graded representation of $\Spin(V)$, the above is summarized as a $\ZZ/2\ZZ$-graded isomorphism
\begin{equation}
\label{iso_repr_theory_Spin}
(\Gamma_{\alpha}\oplus\Gamma_{\beta})^{\otimes 2}\,\iso\, \oplus_{r\ge 0} \; \wedge^r V,
\end{equation}
where the grading on the RHS is given by the parity of $r$. Note that $\dim V$ is divisible by 4. In this case the additional grading by $r\!\mod 4$ on the RHS of (\ref{iso_repr_theory_Spin}) corresponds to the additional grading of the LHS of (\ref{iso_repr_theory_Spin}) by the action of the involution permuting the two factors in the tensor product. We show in Remark~\ref{Rem_great} that the isomorphism of Theorem~\ref{Th_3} respects this additional grading. 

 A trivialization of $\det E$ yields a trivialization of $(A\cJ^{-1})_{\Omega_Y}$. Indeed, we get $\sigma^*\cJ\,\iso\,\cJ^{-1}$. The cup-product $\H^1(Y, \cJ)\otimes \H^1(X, \sigma^*\cJ)\to \H^2(Y, \Qlb)\,\iso\, \Qlb$ is a symplectic form on $\H^1(Y, \cJ)$, our claim easily follows.
 
\smallskip\noindent
iii) The $\GL_2$-variety $\GL_2/\phi_*\Gm$ is spherical (pointwise over $X$), so our setting and results could possibly be thought of in the perspective suggested by Sakellaridis-Venkatesh and Gaitsgory-Nadler in \cite{SV, GN}. 
This is also related to the study of the Waldspurger category in \cite{Ly2}. 
\end{Rem}

\ssec{Quantum geometric Langlands for $\wt\SL_2$}
\label{Sect_Quantum_intro}

\sssec{} Consider the group scheme $G=\Sp(\cO\oplus\Omega)$ on $X$. So, $\Bun_G$ is the stack classifying $M\in\Bun_2$ with an isomorphism $\det M\,\iso\,\Omega$. Let $\cA_G$ be the line bundle on $\Bun_G$ with fibre $\det\RG(X,M)$ at $M$. Write $\Bunt_G$ for the gerbe of square roots of $\cA_G$ over $\Bun_G$. It classifies $M\in\Bun_G$ and a line $\cB$ equipped with $\cB^2\,\iso\, \det\RG(X, M)$. 

 Let $\epsilon$ be the 2-automorphism of $\Bunt_G$ acting as $-1$ on $\cB$ and trivially on $M$. Write $\D_-(\Bunt_G)\subset \D(\Bunt_G)$ for the full subcategory of objects on which $\epsilon$ acts as $-1$. Recall that $\Rep(\SL_2)$ acts on $\D_-(\Bunt_G)$ naturally by Hecke functors (\cite{Ly4, Ly7}). The fundamental problem we are interested in for $\wt\SL_2$ is to find a spectral decomposition of $\D_-(\Bunt_G)$ under this action (see \cite{G1, A+} for many other related ideas).

 Write $\bar\epsilon$ for the 2-automorphism of $\Bunt_G$ that acts trivially on $\cB$ and by $-1$ on $M$. This makes sense as $-1$ acts trivially on $\det\RG(X, M)$. The action of $\bar\epsilon$ gives a $\ZZ/2\ZZ$-grading on $\D_-(\Bunt_G)$. The category $\Rep(\SL_2)$ is also $\ZZ/2\ZZ$-graded by the action of the center of $\SL_2$. By \cite[Lemma~2]{Ly6}, these gradings are compatible.

\sssec{} 
\label{Sect_1.2.2}
Let $H=\SO_3$ split. Write $F_G: \D^-(\Bun_H)_!\to\D^{\prec}(\Bunt_G)$ for the theta-lifting functor from \cite[Section~0.3.2]{Ly6}. We believe that $F_G$ commutes with the actions of $\Rep(\SL_2)$ on both sides, and this could possibly be verified as in \cite{Ly5}, but we did not check this. Unfortunately, the functor $F_G$ is not expected to be an equivalence. This is seen already in Conjecture~\ref{Con_theta_lift_to_Mp2} below. This could be a phenomenon similar to the one we observed for the quantum geometric Langlands for a torus \cite[Question 2.1, Section~5.2.8]{Ly8}, where the twisted derived category of $\Qlb$-sheaves is not equivalent to the corresponding untwisted one (predicted by the quantum Langlands correspondence \cite{G1}). 
 
\sssec{} For a local system $\cE$ on $X$ the geometric central $L$-value is defined as
$$
CL_{\cE}(X)=\oplus_{\theta\ge 0} \RG(X^{(\theta)}, \cE^{(\theta)})[\theta]
$$ 
We view it as $\ZZ/2\ZZ$-graded by the parity of $\theta$. Note that for local systems $\cE_i$ on $X$ one has canonically
$$
CL_{\cE_1\oplus \cE_2}(X)\,\iso\, CL_{\cE_1}(X)\otimes CL_{\cE_2}(X)
$$
 
  For an irreducible $\SL_2$-local system $E$ on $X$ denote by $\Aut_E$ the automorphic perverse sheaf on $\Bun_H$ normalized as in \cite[Section~0.3.3]{Ly6} by the property that its first Whittaker coefficient is `one'. We believe that the corresponding category of $E$-Hecke eigen-sheaves in $\D_-(\Bunt_G)$ contains a unique object, whose $\ZZ/2\ZZ$-graded pieces are irreducible perverse sheaves.
  
\begin{Con} 
\label{Con_theta_lift_to_Mp2}
Let $E$ be an irreducible $\SL_2$-local system on $X$. There is a $\ZZ/2\ZZ$-graded complex $SQ_E\in\D(\Spec k)$ equipped with a $\ZZ/2\ZZ$-graded isomorphism 
\begin{equation}
\label{iso_SQ_E_defining}
(SQ_E)^{\otimes 2}\,\iso\, CL_E(X)\, .
\end{equation} 
Then 
$$
F_G(\Aut_E)\,\iso\, SQ_E\otimes(\cF^+\oplus \cF^-),
$$
where $\cF^{\pm}$ are irreducible perverse sheaves on $\Bunt_G$.
The 2-automorphism $\bar\epsilon$ of $\Bunt_G$ acts on $\cF^{\pm}$ as $\pm$. Moreover, $\cF^+\oplus\cF^-\in \D_-(\Bunt_G)$ is a $E$-Hecke eigen-sheaf. 
\end{Con} 
\begin{Rem} If $E$ is an irreducible $\SL_2$-local system on $X$ then $CL_E(X)$ is placed in cohomological degree zero and identifies with the Clifford algebra of $\H^1(X, E)$. This Clifford algebra splits noncanonically, and the complex $SQ_E$ equipped with the isomorphism (\ref{iso_SQ_E_defining}) exists by the representation theory of the corresponding spinorial group, cf. Remark~\ref{Rem_1.2.13} below.
\end{Rem}

\sssec{}  Let $B\subset G$ be the parabolic group subscheme over $X$ preserving $\Omega$. The stack $\Bun_B$ classifies $\cE\in\Bun_1$ and an exact sequence on $X$
\begin{equation}
\label{sequence_ext_for_intro}
0\to\cE\otimes\Omega\to M\to \cE^{-1}\to 0
\end{equation}
Write $\tilde\nu_B: \Bun_B\to\Bunt_G$ for the map sending the above point to $(M,\cB)$, where $\cB=\det\RG(X, \cE\otimes\Omega)$ with the induced isomorphism $\cB^2\,\iso\, \det\RG(X,M)$. 

  Let $\cS_B$ be the stack classifying $\cE\in\Bun_1$ and $s_2: \cE^2\to\cO_X$. Then $\cS_B$ and $\Bun_B$ are dual generalized vector bundles over $\Bun_1$. Denote by $\Four_{\psi}: \D(\Bun_B)\,\iso\, \D(\cS_B)$ the corresponding Fourier transform. Let $\RCov^d\hook{} \cS_B$ be the open substack classifying $\cE\in\Bun_1, D\in {^{rss}X^{(d)}}$ and $s_2: \cE^2\,\iso\, \cO(-D)$. Let $\varepsilon$ be the 2-automorphism of $\RCov^d$ acting as $-1$ on $\cE$. Denote also by $\varepsilon$ the 2-automorphism of $\Bun_B$ acting as $-1$ on $\cE, M$. The above Fourier transform respects the actions of $\varepsilon$. Over a connected component of $\Bun_B$ containing a point (\ref{sequence_ext_for_intro}) if $\bar\epsilon$ acts on $K\in \D_-(\Bunt_G)$ by $c\in\mu_2$ then $\varepsilon$ acts on $\tilde\nu_B^*K$ as $(-1)^{\chi(\cE\otimes\Omega)}c$. Here $\chi$ stands for the Euler caracteristic. 
  
  As the degree 2 covering $\phi: Y\to X$ varies in $\RCov^d$ the group schemes  $U_{\phi}$ form a group scheme $U_R$ over $\RCov^d\times X$. We define the stack $\Bun_{U_R}$ and a diagram 
$$
\RCov^d\getsup{\gp_R}\Bun_{U_R}\toup{\gq_U}\Bun_H
$$ 
(see Section~\ref{Sect_6.1.4}). The fibre of $(\gp_R)_!\gq_{U}^*K$ at a $k$-point of $\RCov^d$ is the Waldspurger period of $K$ with respect to the trivial local system, this complex organizes these Waldspurger periods into a family over $\RCov^d$. It is remarkable that the Hecke property itself already implies the following acyclicity result (proved in Section~\ref{Sect_6.1.9}). 

\begin{Pp} 
\label{Pp_mainresults_3}
Let $E$ be any $\SL_2$-local system on $X$. Let $K\in\D(\Bun_H)$ be a $E$-Hecke eigensheaf and $d\ge 0$. Then $\gq_{U}^*K$ is ULA with respect to $\gp_R: \Bun_{U_R}\to\RCov^d$. Each cohomology sheaf of the complex $(\gp_R)_!\gq_{U}^*K$ is a local system on $\RCov^d$.
\end{Pp}
  
\sssec{} This paper (and Conjecture~4 iii) from \cite{Ly6}) is mostly concerned with the complex on $\RCov^d$
\begin{equation}
\label{complex_forCon4_on_RCov^d}
\cF_E:=\Four_{\psi}\tilde\nu_B^*F_G(\Aut_E)[\dimrel(\tilde\nu_B)-\dim\RCov^d]\mid_{\RCov^d}  
\end{equation}  
for $E$ an irreducible $\SL_2$-local system on $X$. We show in Corollary~\ref{Cor1_for_Section5.2} that for any $K\in\D(\Bun_H)$ the complex
\begin{equation}
\label{complex_forCor1_on_RCov^d}
\Four_{\psi}\tilde\nu_B^*F_G(K)[\dimrel(\tilde\nu_B)-\dim\RCov^d]\mid_{\RCov^d} 
\end{equation}
identifies with $(\gp_R)_!\gq_{U}^*K$ up to a shift. This is the main motivation for introducing Waldspurger periods, this is how they appeared in the work of Waldspurger \cite{W}, as a tool for calculating the Whittaker coefficients of $F_G(K)$. 

  In \cite[Section~0.3.5]{Ly6} to any local system $V$ on $X$ we associated a complex $CL^d_V$ on $\RCov^d$ whose fibre at $(\cE, s_2)$ equals 
$$
\oplus_{\theta\ge 0} \RG(X^{(\theta)}, (V\otimes\cE_0)^{(\theta)})[\theta],
$$ 
here $\phi: Y\to X$ and $\cE_0$ are associated to $(\cE, s_2)$ as above. Namely, $Y=\Spec(\cO\oplus\cE)$, $\sigma$ is the nontrivial automorphism of $Y$ over $X$, and $\cE_0$ is the sheaf of $\sigma$-anti-invariants in $\phi_*\Qlb$. So, the fibre of $CL^d_V$ over a $k$-point of $\RCov^d$ is the geometric central $L$-value of the constructible sheaf $V\otimes \cE_0$ on  $X$. We consider it as $\ZZ/2\ZZ$-graded by the parity of $\theta$. 
 
  The following is the main result of this paper (it is derived from Theorem~\ref{Th_3} in Section~\ref{Sect_6.1.11}). 
\begin{Cor}  
\label{Cor_main_for_SL_2_SO_3}
Let $E$ be an irreducible $\SL_2$-local system on $X$. There is an isomorphism on $\RCov^d$
\begin{equation}
\label{iso_maybe_of_algebras!!}
\cF_E^{\otimes 2}\,\iso\, CL_{E}(X)\otimes CL^d_{E}
\end{equation}
This isomorphism is $\ZZ/2\ZZ$-graded, where the grading on (\ref{complex_forCon4_on_RCov^d}) is given by the contribution of the stacks $\Bun_{U_R}^a$ for $a\in\ZZ/2\ZZ$, and the complexes $CL^d_{E}$ and $CL_{E}(X)$ are graded by the parity of $\theta$. One has $\DD(\cF_E[d])\,\iso\, \cF_E[d]$. If $d>0$ then $\cF_E$ is a local system on $\RCov^d$ placed in cohomological degree zero. 
\end{Cor}
  
\sssec{} Recall the following \cite[Conjecture~3]{Ly6}.
\begin{Con} 
\label{Con_existence_cS}
For any $\SL_2$-local system $E$ on $X$ and any $d\ge 0$ there is a complex $\cS_E^d\in \D(\RCov^d)$ equipped with an isomorphism 
$$
(\cS_E^d)^{\otimes 2}\,\iso\, CL^d_E
$$
Moreover, $\cS_E^d[\dim\RCov^d]$ is Verdier self-dual. The complex 
$$
\Four_{\psi}\tilde\nu_B^*(\cF^+\oplus \cF^-)[\dimrel(\tilde\nu_B)-\dim\RCov^d]\mid_{\RCov^d} 
$$
identifies canonically with $\cS_E^d$. Here $\cF^+\oplus\cF^-$ is the perverse sheaf on $\Bunt_G$ associated to $E$ in Conjecture~\ref{Con_theta_lift_to_Mp2}. In addition, $SQ_E$ of Conjecture~\ref{Con_theta_lift_to_Mp2} is the fibre of $\cS_E^0$ at the $k$-point $\cE=\cO$ of $\RCov^0$. 
\end{Con}
 
\sssec{} If $E$ is an $\SL_2$-local system on $X$ with $\H^0(X,E)=0$ then $\H^1(X, E)$ is equipped with a nondegenerate symmetric form 
$$
\H^1(X, E)\otimes \H^1(X, E)\to \H^2(X, E\otimes E)\to \H^2(X, \Qlb)\,\iso\, \Qlb
$$
and a compatible trivialization $\det\H^1(X, E)\,\iso\, \Qlb$. So, $\H^1(X, E)$
can be seen as a $\SO_{4g-4}$-torsor over $\Spec \Qlb$. View $CL_E(X)\,\iso\,\oplus_{r\ge 0} \wedge^i \H^1(X, E)$ as the corresponding ($\ZZ/2\ZZ$-graded) Clifford algebra. Let us insist on this point, the central value of the $L$-function $L(E, \frac{1}{2})$, which in the geometric setting is not just a number, but a vector space with the Frobenius operator on it, is actually an algebra.

 Let now $E$ be an $\SL_2$-local system such that for any $d\ge 0$ and a covering $\phi: Y\to X$ in $\RCov^d$ we have $\H^0(Y, \phi^*E)=0$. Recall that the fibre of $CL_E(X)\otimes CL^d_E$ at $(\cE, s_2)$ is 
$$
\otimes_{\theta\ge 0}\RG(Y^{(\theta)}, (\phi^*E)^{(\theta)})[\theta]
$$ 
As above, $\H^1(Y, \phi^*E)$ is equipped with a nondegenerate symmetric form. View $CL_E(X)\otimes CL^d_E$ as the corresponding sheaf of $\ZZ/2\ZZ$-graded Clifford algebras on $\RCov^d$. Assume in addition $E$ irreducible. By Corollary~\ref{Cor_main_for_SL_2_SO_3}, $\cF_E$ is a local system over $\RCov^d$ for all $d\ge 0$. The sheaf of endomorphisms $\cF^{\otimes 2}_E\,\iso\, \END(\cF_E)$ of $\cF_E$ is naturally a $\ZZ/2\ZZ$-graded sheaf of algebras on $\RCov^d$.
 
\begin{Con} 
\label{Con_iso_of_algebras}
 Let $E$ be an $\SL_2$-local system such that for any $d\ge 0$ and a covering $\phi: Y\to X$ in $\RCov^d$ we have $\H^0(Y, \phi^*E)=0$.
Then the isomorphism (\ref{iso_maybe_of_algebras!!}) is actually an isomorphism of $\ZZ/2\ZZ$-graded sheaves of algebras on $\RCov^d$. 
\end{Con}

 Put another way, under the assumptions of Conjecture~\ref{Con_iso_of_algebras}, the local system $CL_{E}(X)\otimes CL^d_{E}$ on $\RCov^d$ has two differerent structures of a $\ZZ/2\ZZ$-graded sheaf of algebras, and we expect these structures to coincide.

\begin{Rem} 
\label{Rem_1.2.13}
i) Let $E$ be an irreducible $\SL_2$-local system on $X$. If we were working over complex numbers with $\cD$-modules, then $\H^1(X,E)$ would carry a pure Hodge structure of weight 1, this structure yields a natural candidate for $SQ_E$. We hope our construction could have applications for the theory of motives (when $E$ is of motivic origin, the corresponding half-spin representations of the Clifford algebras appearing via the Waldspurger periods could be motivic as well).

\smallskip\noindent
ii) If $E$ is an $\SL_2$-local system on $X$ with $\H^0(X,E)=0$ then $\H^1(X, E)$ can be seen as a torsor over $\Spec \Qlb$ under $\SO_{4g-4}$. The datum of $SQ_E$ together with a $\ZZ/2\ZZ$-graded isomorphism (\ref{iso_SQ_E_defining}) is then equivalent to a datum of  lifting of this torsor to a $\Spin_{4g-4}$-torsor over $\Spec\Qlb$. 

\smallskip\noindent
iii) In the setting of Conjecture~\ref{Con_existence_cS} assume $E$ irreducible and $d>0$. For each point of $\RCov^d$ the $*$-fibre of the Clifford algebra $CL^d_E$ as this point splits, that is, identifies as $\ZZ/2\ZZ$-graded with the algebra of endomorphisms of some $\ZZ/2\ZZ$-graded vector space. It is not clear if this is true globally over $\RCov^d$. 
\end{Rem}
 
\section{Strategy of the proof and other results}
\label{Sect_Strategy}

 Essentially all our results about $\wt\SL_2$ formulated in Section~\ref{Sect_Quantum_intro} are an application of Theorem~\ref{Th_3}. Their proofs are given in Section~\ref{Sect_dual_pair_Mp2_SO3}. In the rest of Section~\ref{Sect_Strategy} we explain our approach to the proof of Theorem~\ref{Th_3}.

 We follow the strategy of \cite[Theorem~5]{Ly1}, though the proof of \select{loc.cit} does not generalize `as is' to the case of ramified coverings. The new ingredients needed for the proof of Theorem~\ref{Th_3} are
\begin{itemize}
\item ramified theta-lifting functors for the dual pair $(\GL_2,\GO_2)$ and ramified geometric Eisenstein series for $\GL_2$ described in Section~\ref{Sect_the_dual_pair_GL_GO_1.5}. 
\item Proposition~\ref{Pp_local_nature_gamma_forgetting}, which is a version of the Hecke property of the theta kernel $\Aut_{G,\tilde H}$ for the dual pair $(G=\GL_2, \tilde H=\GO_{2m}^0)$ adopted to our particular type of ramifications.  
\end{itemize} 
We first describe these ingredients and then explain our strategy of the proof of Theorem~\ref{Th_3} in Section~\ref{Sect_Strategy_Th3}. 

\ssec{The dual pair $(\GL_2, \GO_2)$} 
\label{Sect_the_dual_pair_GL_GO_1.5}

\sssec{} 
\label{Sect_2.1.1_now}
Pick a point of $\RCov^d$ given by $s_2: \cE\,\iso\, \cO(-D)$, here $D\in {^{rss}X^{(d)}}$. We have the corresponding degree 2 covering $\phi: Y\to X$.  Let $G=\GL_2$ and $H=\phi_*\Gm$. 
We think of $\D(\Pic Y)$ as some ramified version of $\D(\Pic X\times\Pic X)$. 

 Let $\Bun_{2,D}$ be the stack classifying $M\in\Bun_2$ and a subsheaf $M(-D)\subset \bar M\subset M$ with $\div(M/\bar M)=D$. In Section~\ref{Sect_3.3.5} we define a $\mu_2$-gerbe $\Bunt_{2, D}\to\Bun_{2,D}$ classifying $(\bar M\subset M)\in\Bun_{2,D}$, a $\ZZ/2\ZZ$-graded line $\cU$ of parity zero together with a $\ZZ/2\ZZ$-graded isomorphism 
$$
\cU^2\,\iso\, \det\RG(X, M/\bar M)\otimes\det\RG(X, \cO_D)
$$ 
We think of $\D(\Bunt_{2,D})$ as a ramified version of $\D(\Bun_2)$ taking into account some particular type of ramification (cf. Section~\ref{Sect_3.1_should_be}).
   
 Consider the stack $\Bun_{2,D,H}=\Bun_{2,D}\times_{\Pic X} \Pic Y$, where the map $\Pic Y\to\Pic X$ sends $\cB$ to $\Omega(D)\otimes N(\cB)^{-1}$, and $\Bun_{2,D}\to\Pic X$ sends $(\bar M\subset M)$ to $\det M$. Let $\Bunt_{2,D,H}$ be obtained from $\Bun_{2,D,H}$ by the base change $\Bunt_{2,D}\to\Bun_{2,D}$. We define an object $\Aut_{G,H}\in \D^{\prec}(\Bunt_{2,D,H})$, which is the kernel of the theta-lifting functors
\begin{equation}
\label{functor_F_G_F_H_Sect_2.1.1}
F_G:  \D^-(\Pic Y)\to \D^{\prec}(\Bunt_{2,D})\;\;\;\; F_H: \D^-(\Bunt_{2,D})_!\to \D^{\prec}(\Pic Y)
\end{equation}
(cf. Definition~\ref{Def_theta_lifting_functors_GO_2_GL_2_ram}). One has a decomposition $\Aut_{G,H}\,\iso\, \Aut_{G,H,g}\oplus \Aut_{G,H,s}$. 
   
   The following is an analog of \cite[Proposition~5]{Ly1} in our ramified setting.
\begin{Pp} 
\label{Pp_mainresults_1}
i) Both $\Aut_{G,H, g}$ and $\Aut_{G,H, s}$ are  perverse sheaves irreducible over each connected component of $\Bunt_{2, D, H}$. Moreover, $\DD(\Aut_{G,H})\,\iso\, \Aut_{G, H}$ canonically.  \\
ii) The perverse sheaf $\Aut_{G, H}$ is ULA with respect to $\gq: \Bunt_{2, D, H}\to\Pic Y$. So, $F_G$ commutes with the Verdier duality functors.
\end{Pp}

\sssec{} 
\label{Sect_2.1.3_now}
We define geometric Eisenstein series $\Eis(\Qlb\oplus \cE_0)\in \D(\Bunt_{2,D})$ in the spirit of \cite{BG} (cf. Definition~\ref{Def_3.4.2}).
The second part of the following result geometrizes a particular case of the Siegel-Weil formula (compare with \cite[Proposition~6]{Ly1} in the nonramified setting). 

\begin{Pp}
\label{Pp_mainresults_2}
i)  Let $E$ be a rank one local system on $Y$ that does not descend to $X$. Then $F_G(AE[\dim\Pic Y])$ is an irreducible perverse sheaf over each connected component of $\Bunt_{2,D}$.\\
ii) There is an isomorphism $F_G(\Qlb[\dim\Pic Y])\,\iso\, \Eis(\Qlb\oplus\cE_0)$ on $\Bunt_{2,D}$.   
\end{Pp}

\sssec{} The complexes appearing in Proposition~\ref{Pp_mainresults_2} should be ramified Hecke eigen-sheaves. Indeed, we expect that $F_G$ commutes with Hecke functors, but we do not need this fact in the present paper. 

 In our opinion, it is a little miracle in Proposition~\ref{Pp_mainresults_2} ii), as usually in the Siegel-Weil formula, that given a point $(\bar M\subset M, \cU)\in \Bunt_{2,D}$, certain sum over $\cB\in\Pic Y$ with $N(\cB)\,\iso\, (\det M)^{-1}\otimes\Omega(D)$ rewrites as a sum over some subsheaves of $M$. 
 
 We also describe the Whittaker coefficients of the sheaves appearing in Proposition~\ref{Pp_mainresults_2} (cf.  Lemma~\ref{Lm_restriction_to_Bun_PD} and \ref{Lm_irreducibility_over_cS_PD}). 

\sssec{} 
\label{Sect_2.1.5_should_be}
Let $\epsilon$ be the 2-automorphism of $\Bunt_{2,D}$ acting as $-1$ on $\cU$ and trivially on $(\bar M\subset M)$ for $(\bar M\subset M, \cU)\in \Bunt_{2,D}$. Let $\D_-(\Bunt_{2,D})\subset\D(\Bunt_{2,D})$ be the full subcategory of objects, on which $\epsilon$ acts as $-1$. 

\sssec{}  We introduce the Hecke functor 
$$
\H_G: \D_-(\Bunt_{2,D})\to \D_-(X\times\Bunt_{2,D})
$$
in Section~\ref{Sect_Hecke_functors_2.3.8}. This is some ad hoc extension of the Hecke functor corresponding to the standard representation of the Langlands dual group of $\GL_2$ to the points of ramification (sufficient for our purpose to prove Theorem~\ref{Th_3}). One of the difficulties is that $\H_G$ is known to commute with the Verdier duality only over $(X-D)\times \Bunt_{2,D}$, not over the whole of $X\times \Bunt_{2,D}$.
  
  For the convenience of the reader we explain the definition of $\H_G$ here. Compare with the Hecke functors in the ramified setting used in \cite{Y}. 

 Let $\Mod_{2,D}^r$ be the stack classifying $(\bar M\subset M)\in\Bun_{2,D}$ and an upper modification $M\subset M'$ with $\deg(M'/M)=r$ and $M'\cap \bar M(D)=M$. We have the diagram of projections
\begin{equation}
\label{diag_Hecke_introduc}
\Bun_{2,D}\;\getsup{p_M}\;\Mod_{2,D}^r\;\toup{p'_M} \;\Bun_{2,D},
\end{equation}
where $p_M$ sends the above point to $(\bar M\subset M)$, and $p'_M$ sends the above point to $(\bar M'\subset M')$. Here $M'(-D)\subset \bar M'$ and $\bar M'/M'(-D)\subset M'/M'(-D)$ is the image of the natural map $\bar M/M(-D)\to M'/M'(-D)$. 

Both $p_M, p'_M$ are smooth of relative dimension $2r$. Let $\tilde p_M: \Modt_{2,D}^r\to \Bunt_{2,D}$ be obtained from $p_M$ by the base change $\Bunt_{2,D}\to \Bun_{2,D}$. In Section~\ref{Sect_3.5.3_averaging} 
we extend (\ref{diag_Hecke_introduc}) to a diagram
$$
\Bunt_{2,D}\;\getsup{\tilde p_M}\; \Modt_{2,D}^r\; \toup{\tilde p'_M}\; \Bunt_{2,D},
$$
whose construction does depend on our choice of the covering $\phi: Y\to X$. For $r=1$ this gives the diagram
$$
\Bunt_{2,D} \getsup{\tilde p_M} \Modt^1_{2,D} \toup{\supp\times\tilde p'_M} X\times \Bunt_{2,D},
$$
where $\supp$ sends a point of $\Modt^1_{2,D}$ as above to $\div(M'/M)$. 
For this diagram and $K\in \D_-(\Bunt_{2,D})$ we set
$$
\H_G(K)=(\supp\times \tilde p'_M)_!(\tilde p_M)^*K[2]
$$
The map $\supp\times \tilde p'_M$ is proper over $(X-D)\times\Bunt_{2,D}$, but not over $X\times\Bunt_{2,D}$. The fibre of $\supp\times \tilde p'_M$ over $(x, \bar M\subset M, \cU)$ is $\PP^1$ (resp., $\AA^1$) for $x\notin D$ (resp., $x\in D$). 
  
\sssec{}  For any local system $E$ on $X$ and $r\ge 0$ we introduce in Section~\ref{Sect_3.5.3_averaging} the averaging functors 
$$
\Av^r_E: \D_-(\Bunt_{2,D})\to \D_-(\Bunt_{2,D})\;\;\;\mbox{and}\;\;\; \Av^r_E: \D(\Pic Y)\to \D(\Pic Y)
$$ 
similar to the averaging functors from \cite{G}. The following is an analog of \cite[Propostion~7]{Ly1} in our ramified setting.

\begin{Thm} 
\label{Th_1}
For any local system $E$ on $X$ one has a canonical isomorphism of functors 
$$
F_H\comp \Av^r_E\,\iso\, \Av^r_E\comp F_H
$$
from $\D_-^-(\Bunt_{2,D})_!$ to $\D(\Pic Y)$. 
\end{Thm} 
  
\sssec{}  Theorem~\ref{Th_1} is derived from the following Theorem~\ref{Th_2}. In Section~\ref{Sect_3.6.1_notations} we introduce the stack $\wt\cW_{D,H}$ classifying $x\in X, (\bar M\subset M, \,\cU)\in \Bunt_{2,D}, \cB'\in\Pic Y$ together with an isomorphism $N(\cB')\,\iso\, \cC(D-x)$, where $\cC=\Omega\otimes(\det M)^{-1}$. We also introduce two Hecke functors
$$
\H_{GH}: \D(\Bunt_{2,D,H})\to \D(\wt\cW_{D,H})\;\;\;\mbox{and}\;\;\;
\H_{HG}: \D(\Bunt_{2,D,H})\to \D(\wt\cW_{D,H}),
$$
here $\H_{GH}$ corresponds to the standard representation of the Langlands dual to $G$. The following is an analog of a special case of \cite[Theorem~1]{Ly1} in our ramified setting. 

\begin{Thm}
\label{Th_2}  There is a canonical isomorphism in $\D(\wt\cW_{D,H})$
\begin{equation}
\label{iso_for_Th_2}
\H_{GH}(\Aut_{G,H})\,\iso\, \H_{HG}(\Aut_{G,H})
\end{equation}
\end{Thm}

\sssec{} 
\label{Sect_2.1.10_intro}
While $\H_{HG}$ commutes with the Verdier duality, this property is known for $\H_{GH}$ only over the open substack of $\wt\cW_{D,H}$ given by $x\notin D$. This is one of the technical difficulties in the proof of Theorem~\ref{Th_2}. To derive Theorem~\ref{Th_1} from Theorem~\ref{Th_2}, it is essential to have the isomorphism (\ref{iso_for_Th_2}) over the whole of $\wt\cW_{D,H}$.
 
\ssec{A version of the Hecke property of $\Aut_{G,\tilde H}$}
 
Our Proposition~\ref{Pp_local_nature_gamma_forgetting} is a version of the Hecke property of the theta kernel $\Aut_{G,\tilde H}$ for the dual pair $(G=\GL_2, \tilde H=\GO_{2m}^0)$ adopted to our particular type of ramifications. It says that some averaging of $\Aut$ along $G$ is isomorphic to another averaging of $\Aut$ along $\tilde H$. The proof of Proposition~\ref{Pp_local_nature_gamma_forgetting} is the purpose of Section~\ref{Sect_Restriction_ramifications}, which is independent of the rest of the paper. 


\ssec{Strategy of the proof of Theorem~\ref{Th_3}}
\label{Sect_Strategy_Th3} Our proof consists of the following three steps.
 
\sssec{Step 1} Let $\tilde H$ be given by the exact sequence 
$$
1\to \Gm\to \GL_2\times\GL_2\to \tilde H\to 1,
$$ 
where the first map is $z\mapsto (z, z^{-1})$. Let $\rho_{\tilde H}: \Bun_2\times\Bun_2\to\Bun_{\tilde H}$ be the extension of scalars. Let $E$ be a rank 2 irreducible local system on $X$. We define a perverse sheaf $K_{\pi^*E, \det E, \tilde H}$ on $\Bun_{\tilde H}$, which is the descent of $\Aut_E\boxtimes\Aut_E$ under $\rho_{\tilde H}$. Here we denoted by $\pi: X\sqcup X\to X$ the trivial cover.
 
 As in \cite[Section~6]{Ly1}, define the group scheme $R_{\phi}$ on $X$ by the exact sequence 
$$
1\to \Gm\to \phi_*\Gm\times \phi_*\Gm\to R_{\phi}\to 1,
$$ 
where the first map is $z\mapsto (z, z^{-1})$. 
 
  Let $_{\phi}\GL_2$ be the group scheme on $X$ of automorphisms of $\phi_*\cO_Y$. Define $_{\phi}\tilde H$ by the exact sequence 
$$
1\to\Gm\to (_{\phi}\GL_2)\times(_{\phi}\GL_2)\to {_{\phi}\tilde H}\to 1,
$$ 
where the first map is $z\mapsto (z, z^{-1})$. By \cite[Section~6.1.1]{Ly1}, $\Bun_{_{\phi}\tilde H}\,\iso\,\Bun_{\tilde H}$ naturally. The natural map $R_{\phi}\to {_{\phi}\tilde H}$ induces a morphism denoted $\gq_{R_{\phi}}: \Bun_{R_{\phi}}\to\Bun_{\tilde H}$ in (\select{loc.cit.}, p. 414). 
 
 The product $\phi_*\Gm\times\phi_*\Gm\to\phi_*\Gm$ factors through $R_{\phi}\to \phi_*\Gm$. Let $p_{\phi}: \Bun_{R_{\phi}}\to\Pic Y$ be the extension of scalars with respect to the latter map. 
Unwinding the definitions, the calculation of $\mult_!(\cK\boxtimes\cK)$ is reduced to that of 
\begin{equation}
\label{complex_for_2.0.1}
(p_{\phi})_!\gq_{R_{\phi}}^*K_{\pi^*E, \det E, \tilde H}
\end{equation}
 
Let $G=\GL_2$, let $F_{\tilde H}: \D^-(\Bun_G)_!\to \D^{\prec}(\Bun_{\tilde H})$ be the theta-lifting functor for the dual pair $(G,\tilde H)$ from \cite[Definition~1]{Ly1}. Use \cite[Proposition~8]{Ly1}, which provides an isomorphism
$$
F_{\tilde H}(\Aut_{E^*})\,\iso\, (A(\det E))_{\Omega}\otimes K_{\pi^*E, \det E,\tilde H}
$$

\sssec{Step 2} 
\label{Sect_2.3.2_Step2}
Set $H=\phi_*\Gm$. Our Proposition~\ref{Pp_11_great} identifies (\ref{complex_for_2.0.1}) up to a shift with 
\begin{equation}
\label{complex_essential_for_introd} 
F_H(\Eis(\Qlb\oplus\cE_0)\otimes\delta_{\tilde D}^*\Aut_E)
\end{equation}
Here $\delta_{\tilde D}: \Bunt_{2,D}\to\Bun_2$ sends $(\bar M\subset M, \cU)$ to $M$, and the stack $\Bunt_{2,D}$ is that of Section~\ref{Sect_2.1.1_now}. Here $F_H$ denotes the ramified theta-lifting functor (\ref{functor_F_G_F_H_Sect_2.1.1}), and the Eisenstein series $\Eis(\Qlb\oplus\cE_0)$ is that of (Section~\ref{Sect_2.1.3_now} and Definition~\ref{Def_3.4.2}). This is our motivation for a study of the ramified theta-lifts and Eisenstein series summarized in  Section~\ref{Sect_the_dual_pair_GL_GO_1.5}. 
 
 The proof of Proposition~\ref{Pp_11_great} is a combination of two ingredients: Proposition~\ref{Pp_local_nature_gamma_forgetting}, which was missing in the non-ramified case, and a version of the Siegel-Weil formula given by Proposition~\ref{Pp_mainresults_2} ii).
 
\sssec{Step 3}  We establish Proposition~\ref{Pp_RS_convolution}, which is an analog of \cite[Theorem~2]{Ly1} for our ramified setting. Namely, in Lemma~\ref{Lm_for_RS_convolution_preparatory} we show that 
$$
\Eis(\Qlb\oplus\cE_0)\otimes\delta_{\tilde D}^*\Aut_E
$$ 
is of the form $\oplus_{r\ge 0} \Av^r_E(\La_{\det E})$ for some complex $\La_{\det E}\in \D_-(\Bunt_{2,D})$. The shorthand $\La$ refers to Laumon, as this sheaf is related with Laumon's construction from \cite{Laum}. Applying Theorem~\ref{Th_1}, we get 
$$
F_H\Av^r_E(\La_{\det E})\,\iso\, \Av^r_E F_H(\La_{\det E})
$$ 
The end of the proof is an easy explicit calculation of the right hand side. The answer is of local nature, that is, makes sense for not necessarily irreducible local system $E$ on $X$. This is why, as in \cite{Ly1}, we think of (\ref{complex_essential_for_introd}) as some Rankin-Selberg type convolution.
 
  Combining the above steps, one derives Theorem~\ref{Th_3}.

\ssec{} In Appendix~\ref{Section_appendix} for a point of $\RCov^d$ given by $\phi: Y\to X$ we relate the theta-sheaves for $X$ with those for $Y$ (cf. Proposition~\ref{Pp_last}).

\section{The dual pair $(\GL_2, \GO_2)$ in the ramified setting}
\label{Sect_dual_pair_GL_2_GO_2_ram}

In this section we perform the contructions and establish the results formulated in Section~\ref{Sect_the_dual_pair_GL_GO_1.5}.

\ssec{} 
\label{Sect_3.1_should_be}
The geometric Langlands program in the ramified setting was outlined in \cite{FG}. We will encounter some ramified automorphic sheaves of the following type. Let $G$ be a connected reductive group over $k$. Let $T\subset B\subset G$ be a maximal torus and Borel subgroup. Pick a multiplicity free effective divisor $D$ on $X$. Let $\Bunt_{G,D}$ denote the stack classifying a $G$-torsor $\cF_G$ on $X$, a $B$-torsor $\cF_B$ on $D$, an isomorphism $\cF_G\mid_D\,\iso\, \cF_B\times_B G$, and a trivialization $\gamma: \cF_B\times_B T\,\iso\, \cF^0_T$. Recall that we denoted by $\cF^0_T$ the trivial $T$-torsor. For $x\in D$ pick a rank one character local system $\chi_x$ on $T$, in our case it will be of finite order. The local systems $\chi_x$ will be regular, that is, the stabilizer of $\chi_x$ in the Weyl group is trivial. The group $\prod_{x\in D} T$ acts on $\Bunt_{G,D}$ changing $\gamma$, and we will have to consider $K\in \D(\Bunt_{G,D})$ that change under the action of $\prod_{x\in D} T$ by $\boxtimes_{x\in D} \;\chi_x$ and are automorphic. The Hecke algebra for the corresponding local geometric setting has been studied in \cite{KS}.

\ssec{} 
\label{Sect_2.1.2}
For $n\ge 1$ let $G_n$ be defined as in \cite[Section~2.1]{Ly1}. Namely, $G_n$ is the sheaf on $X$ of automorphisms of $\cO^n_X\oplus\Omega^n$ preserving the natural symplectic form $\wedge^2(\cO^n_X\oplus\Omega^n)\to\Omega$. The stack $\Bun_{G_n}$ classifies $M\in\Bun_{2n}$ with a symplectic form $\wedge^2 M\to\Omega$. Let $\cA_{G_n}$ be the line bundle on $\Bun_{G_n}$ with fibre $\det\RG(X, M)$ at $M$. Let $\Bunt_{G_n}\to\Bun_{G_n}$ be the $\mu_2$-gerbe of square roots of $\cA_{G_n}$. Write $\Aut$ for the theta-sheaf on $\Bunt_{G_n}$, one has a direct sum of perverse sheaves $\Aut\,\iso\, \Aut_g\oplus\Aut_s$ (see \select{loc.cit.}). 

 We will use the notation $G_2$ below, it always refers in this paper to the above sheaf of groups (and should not be confused with the notation $G_2$ used in Cartan's classification of simple algebraic groups). 

 Let $P_n\subset G_n$ be the parabolic group subscheme preserving $\cO^n_X$. Let $\nu_{P_n}: \Bun_{P_n}\to\Bun_{G_n}$ be the natural map. We may view $\Bun_{P_n}$ as the stack classifying $L\in\Bun_n$ and an exact sequence $0\to \Sym^2 L\to ?\to \Omega\to 0$ on $X$, it yields an exact sequence $0\to L\to M\to L^*\otimes\Omega\to 0$. Denote by $\tilde\nu_{P_n}: \Bun_{P_n}\to\Bunt_{G_n}$ the map sending the above point to $(M, \cB=\det\RG(X, L))$ with the induced isomorphism $\cB^2\,\iso\, \det\RG(X, M)$. 

\sssec{} 
\label{Sect_3.2.1}
For $d\ge 0$ let $X^{(d)}$ denote the $d$-th symmetric power of $X$, let $^{rss}X^{(d)}\subset X^{(d)}$ be the open subscheme classifying multiplicity free divisors. As in \cite[Section~6]{Ly1}, let $\RCov^d$ be the stack classifying $D\in {^{rss}X^{(d)}}$, $\cE\in\Bun_1$ and an isomorphism $\cE^2\,\iso\, \cO(-D)$ on $X$. According to \cite[Section7.7.2]{Ly2}, we think of $\RCov^d$ as the stack classifying a degree 2 covering $\phi: Y\to X$ ramified exactly over $D$, where $Y$ is smooth and projective. Namely, $Y=\Spec(\cO\oplus \cE)$ for a point of $\RCov^d$ as above.

 For a constructible $\Qlb$-sheaf $E$ on $X$ let 
$$
E^{(d)}=(\sym_* E^{\boxtimes d})^{S_d}
$$ 
for the map $\sym: X^d\to X^{(d)}$ sending $(x_i)$ to $\sum_i x_i$. 

\sssec{} \select{For the rest of Section~\ref{Sect_dual_pair_GL_2_GO_2_ram} we fix $d\ge 0$ and a $k$-point of $\RCov^d$ given by a degree 2 covering $\phi: Y\to X$ ramified exactly over some $D\in {^{rss}X^{(d)}}$}. 

 Let $D_Y$ be the ramification divisor on $Y$, so $\phi$ induces an isomorphism $D_Y\,\iso\, D$. Let $\sigma$ be the nontrivial automorphism of $Y$ over $X$, let $\cE_{\phi}$ be the $\sigma$-anti-invariants in $\phi_*\cO$. Let $\cE_0$ be the $\sigma$-anti-invariants in $\phi_*\Qlb$. This is a local system on $X-D$ extended by zero to $X$. Similarly, if $r\ge 1$ then $(\cE_0)^{(r)}$ is a local system on $(X-D)^{(r)}$ extended by zero to $X^{(r)}$. 
 
\sssec{}  Let $\Bunt_1$ be the stack classifying $L_1\in \Bun_1$, a $\ZZ/2\ZZ$-graded line $\cU$ of parity zero equipped with a $\ZZ/2\ZZ$-graded isomorphism
\begin{equation}
\label{iso_eta_for_Bunt_1}
\eta: \cU^2\,\iso\, \det\RG(X, L_1\mid_D)\otimes\det\RG(X, \cO_D)
\end{equation} 

\begin{Rem} 
\label{Rem_detRG_O_D}
i) The vector space $\det\RG(X,\cO_D)$ is not canonically trivialized, though one has canonically $\det\RG(X,\cO_D)^2\,\iso\, k$. For $\cA,\cA'\in\Bun_1$ one has a canonical $\ZZ/2\ZZ$-graded isomorphism
$$
\det\RG(X, \cA_D)\otimes\det\RG(X, \cA'_D)\,\iso\, \det\RG(X, \cA\otimes\cA'\mid_D)\otimes\det\RG(X, \cO_D)
$$
ii) For $L_1\in\Bun_1$ the RHS of (\ref{iso_eta_for_Bunt_1}) is $\otimes_{x\in D} \; L_{1,x}$, where each $L_{1,x}$ is of parity zero, so the order of points does not matter.
\end{Rem}

\begin{Lm} 
\label{Lm_detRG_L_mid_D_is_a_square}
For any $L_1,\cC\in\Bun_1$ there is a canonical $\ZZ/2\ZZ$-graded isomorphism
\begin{equation}
\label{detRG_L_mid_D_is_a_square}
\frac{\det\RG(X, L_1\mid_D)\otimes\det\RG(X, L_1\otimes \cE_{\phi}\otimes\cC)^2}{\det\RG(X, L_1\otimes\cC)^2}\,\iso\, 
\frac{\det\RG(X, \cO_D)\otimes\det\RG(X, \cE_{\phi}\otimes\cC)^2}{\det\RG(X, \cC)^2}
\end{equation}
\end{Lm} 
\begin{Prf}
For $\cA,\cB\in\Bun_1$ set 
$$
K(\cA,\cB)=\frac{\det\RG(X, \cA\otimes\cB)\otimes\det\RG(X, \cO)}{\det\RG(X, \cA)\otimes\det\RG(X,\cB)}
$$ 
Then $K(\cA,\cB)$ is bilinear in each variable up to a canonical isomorphism, cf. \cite{De} or \cite[Section~4.2.1-4.2.2]{Ly8}. The isomorphism $K(\cA_1\cB)\otimes K(\cA_2,\cB)\,\iso\, K(\cA_1\otimes\cA_2, \cB)$ shows that
$$
\frac{\det\RG(X, \cA_1\otimes\cA_2\otimes\cB)\otimes\det\RG(X,\cB)}{\det\RG(X, \cA_1\otimes\cB)\otimes\det\RG(X, \cA_2\otimes\cB)}
$$
is independent of $\cB\in\Bun_1$ up to a canonical isomorphism.
This shows that
$$
\frac{\det\RG(X, \cE_{\phi}\otimes\cC)\otimes\det\RG(X, L_1\otimes\cC)}{\det\RG(X, L_1\otimes \cE_{\phi}\otimes\cC)\otimes \det\RG(X,\cC)}\,\iso\, K(\cE_{\phi}, L_1)^{-1}
$$
canonically. Since $K(\cdot, L_1)$ is linear, we get using Remark~\ref{Rem_detRG_O_D}
$$
K(\cE_{\phi}, L_1)^{-2}\,\iso\, K(\cO(-D), L_1)^{-1}\,\iso\, \det\RG(X, L_1\mid_D)\otimes\det\RG(X,\cO_D)
$$
\end{Prf}

\sssec{} Let $\gr: \Bun_1\to\Bunt_1$ be the map sending $L_1$ to $(L_1, \cU,\eta)$, where
\begin{equation}
\label{def_cU_Section_233}
\cU=\frac{\det\RG(X, L_1)\otimes\det\RG(X, \cE_{\phi})}{\det\RG(X, \cO)\otimes\det\RG(X, L_1\otimes\cE_{\phi})}
\end{equation}
is equipped with the isomorphism (\ref{iso_eta_for_Bunt_1}) given by Lemma~\ref{Lm_detRG_L_mid_D_is_a_square}. The map $\gr$ is a $\mu_2$-covering. Let $A\cE_0$ denote the $\mu_2$-anti-invariants in $\gr_!\Qlb$, this is a local system of rank one and order 2 on $\Bunt_1$. The map $\gr$ is a section of the gerbe $\Bunt_1\to\Bun_1$. 

 Note that $A\cE_0$ is trivialized at $(\cO, \cU=k, \eta=\id)$. The stack $\Bunt_1$ is naturally a commutative group stack, the product $m:\Bunt_1\times\Bunt_1\to\Bunt_1$ sends $(L_1, \cU_1,\eta_1)$, $(L_2, \cU_2,\eta_2)$ to $({L_1\otimes L_2}, \cU_1\otimes\cU_2,\eta)$, where $\eta$ is the composition
$$
\cU_1\otimes\cU_2\toup{\eta_1\otimes\eta_2} \frac{\det\RG(X, L_1\mid_D)\otimes\det\RG(X, L_2\mid_D)}{\det\RG(X, \cO_D)^2}\;\iso\; \frac{\det\RG(X, L_1\otimes L_2\mid_D)}{\det\RG(X, \cO_D)},
$$
the second isomorphism is that of Remark~\ref{Rem_detRG_O_D}. One checks that $m^*A\cE_0\,\iso\, A\cE_0\boxtimes A\cE_0$ naturally. 

One has the Abel-Jacobi map 
$$
AJ: (X-D)^{(r)}\to \Bunt_1
$$ 
sending $D_1$ to $(\cO(D_1), \cU=k, \eta)$, where $\eta=\id$. One has $AJ^*(A\cE_0)\,\iso\, \cE_0^{(r)}$ naturally. Indeed, the pull-back of the cover $\gr$ by $AJ$ is a degree two cover over $(X-D)^{(r)}$, whose fibre at $D_1\in (X-D)^{(r)}$ is the set of $z\in \mathop{\otimes}\limits_{x\in D_1} (\cE_{\phi})_x$ with $z^2=1$. We used here the isomorphism $\cE_{\phi}^2\,\iso\, \cO(-D)$ fixed above.

\begin{Rem} i) For a rank one local system $\cV$ on $X$ denote by $A\cV$ the corresponding automorphic local system on $\Bun_1$. For $r\ge 0$ and the Abel-Jacobi map $AJ: X^{(r)}\to\Bun_1$, $D_1\mapsto \cO(D_1)$ one has $AJ^*A\cV\,\iso\, \cV^{(r)}$. For $d=0$ there is an ambiguity for our notation $A\cE_0$. In this case the gerbe $\Bunt_1\to \Bun_1$ admits a section $s_Y: \Bun_1\to\Bunt_1$
sending $L_1$ to $(L_1, \cU=k, \eta=\id)$. One has $s_Y^*A\cE_0\,\iso\, A\cE_0$. We hope the precise meaning for $A\cE_0$ is clear from the context. 

\smallskip\noindent
iii) If $d=0$ then (\ref{def_cU_Section_233}) is trivialized for $L_1=\Omega$, so $(A\cE_0)_{\Omega}\,\iso\,\Qlb$ in this case. This is why Theorem~\ref{Th_3} for $d=0$ coincides with \cite[Theorem~5]{Ly1}. 
\end{Rem}


\ssec{Theta-lifting functors for $(\GL_2, \GO_2)$}

\sssec{} 
\label{Sect_3.3.1_should_be}
Let $\Bun_{2,D}$ be the stack classifying $M\in\Bun_2$ together with a subsheaf $M(-D)\subset \bar M\subset M$ such that $\div(M/\bar M)=D$. This is the stack of vector bundles with parabolic structure at $D$. Let $\Sh_0^d$ denote the stack of torsion sheaves on $X$ of length $d$. We have used the map $\div: \Sh_0^d\to X^{(d)}$ from \cite{Ly3}. 

 Let $H=\phi_*\Gm$, the stack $\Bun_H$ of $H$-torsors on $X$ identifies naturally with $\Pic Y$. Set 
$$
\Bun_{2,D,H}=\Bun_{2,D}\times_{\Pic X} \Pic Y,
$$ 
where the map $\Pic Y\to\Pic X$ sends $\cB$ to $\Omega(D)\otimes N(\cB)^{-1}$, and $\Bun_{2,D}\to\Pic X$ sends $(\bar M\subset M)$ to $\cA:=\det M$. We have used the norm map $N: \Pic Y\to\Pic X$ from \cite[Appendix~A]{Ly1} given by $N(\cB)=\cE_{\phi}^{-1}\otimes\det(\phi_*\cB)$.  
 
  For $\cB\in\Pic Y$ the vector bundle $L=\phi_*\cB$ is equipped with a canonical symmetric form $\Sym^2 L\to \cC(D)$, where $\cC=N(\cB)(-D)$. Let also $\bar L=\phi_*(\cB(-D_Y))$, the symmetric form on $L$ induces an isomorphism $L\,\iso\, \bar L^*\otimes\cC$. The symmetric form $\Sym^2 \bar L\to \cC$ has a canonical section $s_c$ that fits into the commutative diagram
\begin{equation}
\label{diag_def_of_s_c}
\begin{array}{ccc}
\Sym^2 \bar L & \to & \cC\\
\uparrow\lefteqn{\scriptstyle s_c} & \nearrow\lefteqn{\scriptstyle\id}\\
\cC(-D)
\end{array}
\end{equation}   
  
\begin{Lm} 
\label{Lm_detRG_M_otimes_L}
For $(\bar M\subset M, \cB)\in \Bun_{2,D,H}$ with $\cA=\det M$, $L=\phi_*\cB$ one has a canonical $\ZZ/2\ZZ$-graded isomorphism
$$
\det\RG(X, M\otimes L)\,\iso\, \frac{\det\RG(X, M)^2\otimes\det\RG(X, L)^2}
{\det\RG(X, \cO)^2\otimes\det\RG(X, \cA)^2\otimes \det\RG(Y, \cB/\cB(-D_Y))}
$$
\end{Lm} 
\begin{Prf}
Let $\bar L=\phi_*(\cB(-D_Y))$. By \cite[Lemma~1]{Ly1}, one gets
$$
\det\RG(X, M\otimes L)\,\iso\, \frac{\det\RG(X, M)^2\otimes\det\RG(X, L)^2\otimes\det\RG(X, \cE_{\phi})}{\det\RG(X,\cO)^3\otimes\det\RG(X, \cA)\otimes\det\RG(X, \cA\otimes\cE_{\phi})}
$$
Applying this formula for $M=\cO\oplus\cA$ in particular and using $L\,\iso\, \bar L^*\otimes\cC$, one gets
$$
\det\RG(X, \bar L)\,\iso\, \det\RG(X, \cA\otimes L)\,\iso\, \frac{\det\RG(X,\cA)\otimes\det\RG(X, L)\otimes\det\RG(X, \cE_{\phi})}{\det\RG(X, \cE_{\phi}\otimes\cA)\otimes \det\RG(X,\cO)}
$$
Since $\det\RG(X, \bar L)\otimes\det\RG(Y, \cB/\cB(-D_Y))\,\iso\, \det\RG(X, L)$, our claim follows.  
\end{Prf}

\sssec{} Denote by
\begin{equation}
\label{map_tau_Section_221}
\tau: \Bun_{2, D, H}\to \Bun_{G_2}
\end{equation}
the following map. For a point of $\Bun_{2, D, H}$ on $M\otimes L$ one gets a form $\wedge^2(M\otimes L)\to \Omega(D)$ such that $M\otimes \bar L$ is the orthogonal complement of $M\otimes L$ with respect to the form with values in $\Omega$. Here $\bar L=\phi_*\cB(-D_Y)$. Note that $L/\bar L\,\iso\, \cB/\cB(-D_Y)$ as a $\cO_X$-module via the isomorphism $D_Y\,\iso\, D$ induced by $\phi$. Let $M'$ be defined by the cartesian square
$$
\begin{array}{ccc}
M\otimes L & \to & M\otimes (L/\bar L)\\
\uparrow && \uparrow\\
M' & \to & (\bar M/M(-D))\otimes (L/\bar L)
\end{array}
$$
The map $\tau$ sends the above collection to $M'$.    

\begin{Lm}
\label{Lm_calculation_detRG_M'}
 For a point $(\bar M\subset M, \cB)\in \Bun_{2, D, H}$ let $M'=\tau(\bar M\subset M, \cB)$. The $\ZZ/2\ZZ$-graded line $\det\RG(X, M')$ identifies canonically with
$$
\frac{\det\RG(X, M)^2\otimes\det\RG(X, \bar L)^2\otimes\det\RG(X, \cO_D)}
{\det\RG(X, \cO)^2\otimes\det\RG(X, \cA)^2\otimes \det\RG(X, (M/\bar M))}
$$
Here $\cA=\det M$, $\bar L=\phi_*(\cB(-D_Y))$. 
\end{Lm}
\begin{Prf} One has $\det\RG(X, M')\otimes\det\RG(X, (M/\bar M)\otimes(L/\bar L))\,\iso\, \det\RG(X, M\otimes L)$ canonically. Applying Lemma~\ref{Lm_detRG_M_otimes_L}, one gets canonically
$$
\det\RG(X, (M/\bar M)\otimes(L/\bar L))\otimes\det\RG(X, \cO_D)\,\iso\, \det\RG(X, M/\bar M))\otimes\det\RG(X, L/\bar L)
$$  
\end{Prf} 

\sssec{} 
\label{Sect_3.3.5}
Let $\Bunt_{2, D}\to \Bun_{2, D}$ be the $\mu_2$-gerbe classifying a point $(\bar M\subset M)\in \Bun_{2, D}$, a $\ZZ/2\ZZ$-graded line $\cU$ of parity zero and a $\ZZ/2\ZZ$-graded isomorphism
\begin{equation}
\label{iso_cU_squared_for_Bunt_2D} 
 \cU^2\,\iso\, \det\RG(X, M/\bar M)\otimes\det\RG(X, \cO_D)
\end{equation} 
Let $\Bunt_{2, D, H}$ be obtained from $\Bun_{2, D, H}$ by the base change $\Bunt_{2, D}\to \Bun_{2, D}$. We get a morphism
\begin{equation}
\label{map_tilde_tau_Section_221}
\tilde\tau: \Bunt_{2, D, H}\to\Bunt_{G_2}
\end{equation}
over $\tau$ sending $(\bar M\subset M, \cU, \cB)$ to $(M', \cU')$, where 
\begin{equation}
\label{def_cU'_for_tilde_tau}
\cU'=\frac{\det\RG(X, M)\otimes\det\RG(X, \bar L)}
{\det\RG(X, \cO)\otimes\det\RG(X, \cA)\otimes \cU}
\end{equation}
is equipped with the isomorphism $\cU'^2\,\iso\, \det\RG(X, M')$ given by Lemma~\ref{Lm_calculation_detRG_M'}. Set $G=\GL_2$.

\begin{Def} 
\label{Def_theta_lifting_functors_GO_2_GL_2_ram}
Recall the perverse sheaf $\Aut$ from Section~\ref{Sect_2.1.2}. Set $\Aut_{G,H}=\tilde\tau^*\Aut[\dimrel(\tilde\tau)]\in \D^{\prec}(\Bunt_{2, D, H})$ and similarly for $\Aut_{G, H, g}$, $\Aut_{G, H, s}$. As in \cite{Ly1} for the diagram of projections
$$
\Pic Y \,\getsup{\gq} \,\Bunt_{2, D, H} \,\toup{\gp} \,\Bunt_{2, D}
$$
define $F_G: \D^-(\Pic Y)\to \D^{\prec}(\Bunt_{2, D})$ and $F_H: \D^-(\Bunt_{2, D})_!\to \D^{\prec}(\Pic Y)$ by
$$
F_G(K)=\gp_!(\Aut_{G, H}\otimes\gq^*K)[-\dim\Pic Y]
$$
$$
F_H(K)=\gq_!(\Aut_{G, H}\otimes\gp^*K)[-\dim\Bunt_{2,D}]
$$
\end{Def}

\begin{Rem} i) The group stack $\Bunt_1$ acts on $\Bunt_{2,D}$. The action map $\act: \Bunt_1\times\Bunt_{2,D}\to\Bunt_{2,D}$ sends $(L,\cU)\in\Bunt_1, (\bar M\subset M,\cU_1)\in\Bunt_{2,D}$ to $(L\otimes\bar M\subset L\otimes M, \cU\otimes\cU_1)$ with the induced isomorphism 
\begin{multline*}
(\cU\otimes\cU_1)^2\,\iso\, \det\RG(X, L\mid_D)\otimes\det\RG(X, M/\bar M)\,\iso\\
 \det\RG(X, L\otimes M/L\otimes \bar M)\otimes\det\RG(X, \cO_D)
\end{multline*}
Given a rank one local system $\cV$ on $X$, say that $K\in\D(\Bunt_{2,D})$ has central character $\cV\otimes\cE_0^m$ if it is equipped with $(\Bunt_1, A\cV\otimes (A\cE_0)^m)$-equivariant structure as in \cite[Definition~3]{Ly1}. This means, in particular, that $\act^*K\,\iso\, (A\cV\otimes (A\cE_0)^m)\boxtimes K$.
By abuse of notation we denoted the restriction of $A\cV$ under $\Bunt_1\to\Bun_1$ also by $A\cV$. 

 The notion of a central character for $K\in\D(\Pic Y)$ is defined similarly for the action map $\act: \Pic X\times\Pic Y\to\Pic Y$, $(\cA,\cB)\mapsto \phi^*\cA\otimes\cB$.

\smallskip\noindent
ii) For $d=0$ we have a section $\Bun_{2,D}\to \Bunt_{2,D}$ given by $\cU=k$ with (\ref{iso_cU_squared_for_Bunt_2D}) being the identity. The inverse image with respect to this section gives an equivalence $\D_-(\Bunt_{2,D})\,\iso\, \D(\Bun_{2,D})$. Under this equivalence the functors $F_G,F_H$ identify with those of \cite[Definition~1]{Ly1} in this case. 

\smallskip\noindent
iii) One has the involution of $\Pic Y$ given by $\cB\mapsto\sigma^*\cB$, it induces an involution of $\Bunt_{2,D,H}$, and $\tilde\tau$ is invariant under this involution.
\end{Rem}

\sssec{} Let $P\subset G=\GL_2$ be the Borel subgroup of upper-triangular matrices. Write $\Bun_{P,D}$ for the stack classifying $L_1,\cA\in\Bun_1$,
an exact sequence 
\begin{equation}
\label{seq_L*otimescA_by_L}
0\to L_1\to M\to L_1^*\otimes\cA\to 0
\end{equation}
on $X$, and a lower modification $M(-D)\subset \bar M\subset M$ with $\div(M/\bar M)=D$ such that $L_{1,x}\cap (\bar M/M(-x))=0$ for any $x\in D$. Here $L_{1,x}$ is the geometric fibre of $L_1$ at $x$. The datum of $\bar M$ can be seen as a splitting of the restriction of (\ref{seq_L*otimescA_by_L}) to $D$. For a point of $\Bun_{P,D}$ the inclusion $L_1(-D)\subset \bar M$ is a subbundle.

 The stack $\Bun_{P,D}$ can be seen as a stack classifying $L_1,\cA\in\Bun_1$ and an exact sequence 
\begin{equation}
\label{seq_L*otimescA_by_L(-D)} 
0\to L_1(-D)\to \bar M\to L_1^*\otimes\cA\to 0
\end{equation}
Namely, (\ref{seq_L*otimescA_by_L}) is the push-forward of (\ref{seq_L*otimescA_by_L(-D)}) by $L_1(-D)\hook{} L_1$. 
The map $\Bun_{P,D}\to\Bun_1\times\Bun_1$ sending the above point to $(L_1,\cA)$ is a generalized vector bundle in the sense of \cite{Ly3}. 

 Let $\nu_{P,D}: \Bun_{P,D}\to\Bun_{2,D}$ be the map sending the above point to $(\bar M\subset M)$. Set 
$$
\Bun_{P,D,H}=\Bun_{2,D,H}\times_{\Bun_{2,D}}\Bun_{P,D}
$$
For a point of $\Bun_{P,D}$ one has $M/\bar M\,\iso\, L_1\mid_D$ naturally. We extend $\nu_{P,D}$ to the morphism 
$
\tilde\nu_{P,D}: \Bun_{P,D}\to\Bunt_{2,D}
$ 
sending (\ref{seq_L*otimescA_by_L(-D)}) to $(\bar M\subset M, \cU)$, where $\cU$ is given by (\ref{def_cU_Section_233}) and equipped with the isomorphism 
$$
\cU^2\,\iso\, \det\RG(X, L_1\mid_D)\otimes\det\RG(X, \cO_D)\,\iso\, \det\RG(X, M/\bar M)\otimes\det\RG(X, \cO_D)
$$ 
given by Lemma~\ref{Lm_detRG_L_mid_D_is_a_square}. 
 
  Recall the group scheme $P_2$ on $X$ defined in Section~\ref{Sect_2.1.2} right above \ref{Sect_3.2.1}. Let $\tau_P: \Bun_{P,D, H}\to\Bun_{P_2}$ be the following map. For a point of $\Bun_{P,D, H}$ given by (\ref{seq_L*otimescA_by_L}) together with $(\bar M\subset M)$, $\cB\in\Pic Y$ with $\cC=N(\cB)(-D)$
let $M'=\tau(\bar M\subset M, \cB)$. Then
$$
(L_1\otimes L)\cap M'=L_1\otimes\bar L,
$$ 
so $L_1\otimes\bar L\subset M'$ is a lagrangian subbundle. Here $L=\phi_*\cB, \bar L=\phi_*(\cB(-D_Y))$. The map $\tau_P$ sends the above point to $(L_1\otimes\bar L\subset M')$.  

 Consider the diagram
\begin{equation}
\label{diag_tilde_tau_lifted_to_P}
\begin{array}{ccc}
\Bunt_{2,D,H} & \toup{\tilde\tau} & \Bunt_{G_2}\\
\uparrow\lefteqn{\scriptstyle \tilde\nu_{P,H}}&& \uparrow\lefteqn{\scriptstyle \tilde\nu_{P_2}}\\
\Bun_{P,D,H} & \toup{\tau_P} & \Bun_{P_2},
\end{array}
\end{equation} 
where we have set $\tilde\nu_{P,H}=\tilde\nu_{P,D}\times\id$, and $\tilde\nu_{P_2}$ was defined in Section~\ref{Sect_2.1.2} right above \ref{Sect_3.2.1}.

\begin{Lm} 
\label{Lm_root_of_MoverbarM}
The diagram (\ref{diag_tilde_tau_lifted_to_P}) is naturally 2-commutative.
\end{Lm}  
\begin{Prf}
By \cite[Lemma~1]{Ly1}, one has
$$
\det\RG(X, L_1\otimes\bar L)\,\iso\, \frac{\det\RG(X, L_1)\otimes\det\RG(X, \bar L)\otimes\det\RG(X, L_1\otimes\cE_{\phi}\otimes\cC)}{\det\RG(X, \cE_{\phi}\otimes\cC)\otimes\det\RG(X, \cO)}
$$
This implies that for a point of $\Bun_{P,D,H}$ we have a natural $\ZZ/2\ZZ$-graded isomorphism $\det\RG(X, L_1\otimes\bar L)\,\iso\, \cU'$ , where $\cU'$ is given by (\ref{def_cU'_for_tilde_tau}), and $\cU$ is given by (\ref{def_cU_Section_233}).  
\end{Prf}

\begin{Lm} 
\label{Lm_tau_P_explicitly}
Consider a point of $\Bun_{P,D, H}$ as above. Its image under $\tau_P$ is given by the exact sequence $0\to \Sym^2(L_1\otimes\bar L)\to ?\to \Omega\to 0$, which is the push-forward of (\ref{seq_L*otimescA_by_L(-D)}) by the composition
$$
L_1^2\otimes \cA^{-1}(-D)\,\iso\, L_1^2\otimes\Omega^{-1}\otimes \cC(-D)\toup{\id\otimes s_c} L_1^2\otimes\Omega^{-1}\otimes\Sym^2 \bar L,
$$
where $s_c:\cC(-D)\to \Sym^2 \bar L$ is the canonical section of the symmetric form (\ref{diag_def_of_s_c}). 
\QED
\end{Lm}

\sssec{} Let $\cS_{P,D}$ be the stack classifying $L_1,\cA\in\Bun_1$ and $v: L_1^2\otimes\cA^{-1}(-D)\to\Omega$. We have a diagram of dual generalized vector bundles 
\begin{equation}
\label{diag_dual_vector_bundles_Bun_PD}
\cS_{P,D}\to \Bun_1\times\Bun_1Ê\gets \Bun_{P,D},
\end{equation}
where each map sends a point as above to $(L_1, \cA)\in\Bun_1\times\Bun_1$. 
 
 Let $\cS_{P,D,H}$ be the stack classifying $(L_1,\cA, v)\in \cS_{P,D}$, $\cB\in\Pic Y$ with $\cC=N(\cB)(-D)$ and an isomorphism $\cA\otimes\cC\,\iso\,\Omega$. By base change from (\ref{diag_dual_vector_bundles_Bun_PD}) one gets a diagram
$$
\cS_{P, D,H}\,\iso\, \cS_{P,D}\times_{\Bun_1}\Pic Y \to \Bun_1\times\Pic Y \gets \Bun_{P,D}\times_{\Bun_1}\Pic Y \,\iso\, \Bun_{P, D,H},
$$
and we write 
$
\Four_{\psi}: \D^{\prec}(\cS_{P,D,H})\to \D^{\prec}(\Bun_{P,D,H})
$
for the corresponding Fourier transform. 

 Let $\cV_{P,D,H}$ be the stack classifying $L_1\in\Bun_1$, $\cB\in \Pic Y$ for which we set $\bar L=\phi_*\cB(-D_Y)$, $\cC=N(\cB)(-D)$, $\cA=\Omega\otimes \cC^{-1}$, and a section $t: L_1\otimes\bar L\to\Omega$. Let
$$
\gp_{\cV,H}: \cV_{P,D,H}\to \cS_{P, D,H}
$$
be the map sending this point to $(L_1,\cB, v)$, where $v$ is the composition
\begin{equation}
\label{first_def_of_v}
L_1^2\otimes\cC(-D)\toup{\id\otimes s_c} L_1^2\otimes\Sym^2 \bar L\,\iso\, \Sym^2(L_1\otimes\bar L)\toup{t\otimes t} \Omega^2
\end{equation}
Recall the complex $S_{P,\psi}$ on $\Bun_{P_2}$ from \cite[Definition~3]{Ly4}. 
The following is a straightforward consequence of Lemma~\ref{Lm_tau_P_explicitly}. 
\begin{Pp} 
\label{Pp_description_tau_P_S_Ppsi}
There is a canonical isomorphism in $\D(\Bun_{P,D, H})$
$$
\tau_P^* S_{P,\psi}[\dimrel(\tau_P)]\,\iso\, \Four_{\psi}((\gp_{\cV, H})_!\Qlb[\dim \cV_{P,D,H}])
$$
\QED
\end{Pp} 

\sssec{} The stack $\cV_{P,D, H}$ can be seen as the one classifying $L_1\in\Bun_1, \cB\in\Pic Y$ and 
$$
t: \phi^*L_1\to \cB^*\otimes\Omega_Y(D_Y)
$$ 
In these terms the map $v$ given by (\ref{first_def_of_v}) is nothing but the norm $Nt: L_1^2\to \cC^{-1}\otimes\Omega^2(D)$. We have used the fact that $\Omega_Y(-D_Y)\,\iso\, \phi^*\Omega$ canonically, here $\Omega_Y$ denotes the canonical line bundle on $Y$.  

\sssec{} For the rest of Section~\ref{Sect_dual_pair_GL_2_GO_2_ram} 
we assume $d=\deg D>0$, that is, $\phi: Y\to X$ is indeed ramified. The case $d=0$ is established in \cite{Ly1}. 

Let us generalize \cite[Lemma~6]{Ly1} to our ramified case. Let $\Pic^r Y$ be the connected component of $\Pic Y$ classifying $\cB\in\Pic Y$ with $\deg\cB=r$. Write $\Pic'^r Y$ for the stack classifying $\cB\in\Pic^r Y$ together with $t: \cO\to \cB$. One defines $\Pic'^r X$ similarly. Let
$$
\phi'_{ex}: \Pic'^r Y\to \Pic'^r X\times_{\Pic X} \Pic Y
$$
be the map sending $(t: \cO\to\cB)$ to $(\cB, Nt: \cO\to N(\cB))$. The group $S_2$ acts on $\phi'_{ex}$ sending $(\cB, t)$ to $(\cB, -t)$. We have the open immersion $X^{(r)}\hook{} \Pic'^r X$ sending $D_1$ to $(t: \cO\to \cO(D_1))$. Over this open substack the map $\phi'_{ex}$ restricts to a morphism
$$
\phi': Y^{(r)}\to X^{(r)}\times_{\Pic X} \Pic Y
$$

\begin{Lm} 
\label{Lm_perversity_for_phi'_ex}
i) Let $\cF$ be an irreducible perverse sheaf on $X$ such that both $V=\cF[-1]$ and $(\DD\cF)[-1]$ are constructible sheaves, that is, placed in usual cohomological degree zero. Then $V^{(r)}[r]$ is also an irreducible perverse sheaf.

\smallskip\noindent
ii) For any $r\ge 0$ both the $S_2$-invariants and anti-invariants in $\phi'_!\Qlb[r]$ are irreducible perverse sheaves. If $r>2g_Y-2$ then the same holds for $(\phi'_{ex})!\Qlb[r]$.
\end{Lm}
\begin{Prf}
i) similar to \cite[Theorem~3.3.10]{Laum}.

\smallskip\noindent
ii) The morphism $\phi'_{ex}$ is finite. Write $\uPic Y$ for the Picard scheme of $Y$, similarly for $X$. We have a $\mu_2$-gerbe
$$
\und{\gr}: X^{(r)}\times_{\Pic X}\Pic Y\to X^{(r)}\times_{\uPic X} \uPic Y
$$
The proof of \cite[Lemma~6]{Ly1} still holds if in Step 1 we take $E$ to be a rank one local system on $Y$ that does not descend to $X$. In this case $\phi_*E[1]$ is an irreducible perverse sheaf on $X$, so $(\phi_*E)^{(r)}[r]$ is also an irreducible perverse sheaf on $X^{(r)}$ by i).
\end{Prf}

\sssec{} Let $^0\Bun_{P,D}\subset \Bun_{P,D}$ be the open substack given by the following condition

\begin{itemize}
\item[(C)] \  $2\deg L_1-\deg\cA+2g-2<0$, and if $g=0$ then moreover $2\deg L_1-\deg\cA<1$. 
\end{itemize}
\noindent
The restriction $^0\nu_{P,D}:{^0\Bun_{P,D}}\to\Bun_{2,D}$ of $\nu_{P,D}$ is smooth and surjective. Let $^0\cS_{P,D}\subset \cS_{P,D}$ be the open substack given by (C), and similarly for $^0\cS_{P,D,H}$. Set 
$$
^0\cV_{P,D,H}=\gp_{\cV, H}^{-1}(^0\cS_{P, D,H}),
$$ 
this is a smooth stack.

\smallskip\noindent
\begin{Prf}\select{of Proposition~\ref{Pp_mainresults_1}} \\
i) Let 
$
^0\Bun_{P,D,H}\subset \Bun_{P,D,H}
$ 
be the open substack given by (C). By the above, the restriction ${^0\Bun_{P,D,H}}\to \Bunt_{2,D,H}$ of $\tilde\nu_{P,H}$ is smooth and surjective. So, it suffices to show that both 
$$
\tau_P^*S_{P,\psi, g}[\dimrel(\tau_P)]\;\;\;\mbox{and}\;\;\; \tau_P^*S_{P,\psi, s}[\dimrel(\tau_P)]
$$ 
are irreducible perverse sheaves over each connected component of $^0\Bun_{P,D,H}$. 

  By ii) of Lemma~\ref{Lm_perversity_for_phi'_ex}, 
$(\gp_{\cV, H})_!\Qlb[\dim {^0\cV_{P,D,H}}]$ over any connected component of $^0\cS_{P,D,H}$ is a direct sum of two irreducible perverse sheaves. Since $\gp_{\cV, H}$ is finite, this perverse sheaf is self-dual. Since $\Four_{\psi}$ preserves perversity and irreducibility, part i) follows.

\smallskip\noindent
ii) Let $^0(\Bun_1\times\Pic Y)\subset \Bun_1\times\Pic Y$ be the open substack classifying $L_1\in\Bun_1,\cB\in\Pic Y$ with $2\deg L_1+\deg\cB-\deg D<0$. The projection $^0\cV_{P,D,H}\to {^0(\Bun_1\times\Pic Y)}$ forgetting $t$ is a vector bundle, so $(\gp_{\cV,H})_!\Qlb[\dim {^0\cV_{P,D,H}}]$ is ULA over $^0(\Bun_1\times\Pic Y)$. 
As in the proof of \cite[Proposition~5 ii)]{Ly1}, the Fourier transform preserves the ULA property over a smooth base, so $\tau_P^*S_{P,\psi}$ is ULA with respect to the projection $^0\Bun_{P,D,H}\to {^0(\Bun_1\times\Pic Y)}$, hence also ULA over $\Pic Y$. Since $\tilde\nu_{P,H}: {^0\Bun_{P,D,H}}\to \Bunt_{2,D,H}$ is smooth and surjective, the desired ULA property follows. 

 Since $\gp: \Bunt_{2,D,H}\to\Bunt_{2,D}$ is proper, $F_G$ commutes with the Verdier duality as in \cite[Section~5.1.2]{BG}.
\end{Prf}
 
\sssec{} Write $\gp_{\cS}$ for the composition $\cV_{P,D, H}\toup{\gp_{\cV, H}} \cS_{P,D,H}\to \cS_{P,D}$, where the second map is the projection forgetting $\cB$. Consider the diagram
\begin{equation}
\label{diag_for_F_cS}
\cS_{P,D}\getsup{\gp_{\cS}}\cV_{P,D, H} \toup{\gq_{\cS}} \Pic Y,
\end{equation}
where $\gq_{\cS}$ is the projection sending $(L_1,\cB, t)$ to $\cB$. From Proposition~\ref{Pp_description_tau_P_S_Ppsi} one gets the following.
 
\begin{Lm} 
\label{Lm_restriction_to_Bun_PD}
For $K\in\D(\Pic Y)$ one has naturally in $\D(\Bun_{P,D})$
$$
\tilde\nu_{P,D}^*F_G(K)\,\iso\, \Four_{\psi}\gp_{\cS !}\gq_{\cS}^*K[\dimrel(\gq_{\cS})]
$$ 
\QED
\end{Lm}

\sssec{} For $r\ge 0$ we define an open immersion $j_r: \Bun_1\times X^{(r)}\to \cS_{P,D}$ sending $(L_1, D_1)$ to $(L_1, \cA, v)$, where $\cA=\Omega^{-1}\otimes L_1^2(D_1-D)$ with the canonical inclusion 
$$
v: L_1^2\otimes\cA^{-1}(-D)\,\iso\,\Omega(-D_1)\hook{}\Omega
$$
Let $j_{r,Y}: \Bun_1\times Y^{(r)}\to \cV_{P,D,H}$ be the map sending $(L_1, \bar D_1)$ to $L_1, \cB=\Omega_Y(D_Y-\bar D_1)\otimes\phi^*L_1^*$ with the canonical inclusion $t: \phi^*L_1\to \cB^*\otimes\Omega_Y(D_Y)$. The following square is cartesian
\begin{equation}
\label{diag_j_r_restriction}
\begin{array}{ccc}
\cV_{P,D,H} & \getsup{j_{r,Y}} & \Bun_1\times Y^{(r)}\\
\downarrow\lefteqn{\scriptstyle \gp_{\cS}} && \downarrow\lefteqn{\scriptstyle \id\times \phi^r}\\
\cS_{P,D} & \getsup{j_r} & \Bun_1\times X^{(r)},
\end{array}
\end{equation}
where $\phi^r: Y^{(r)}\to X^{(r)}$ is the map $\bar D_1\mapsto \phi_*\bar D_1$. 

 For a rank one local system $E$ on $Y$ write $AE$ for the corresponding automorphic local system on $\Pic Y$. For $r\ge 0$ and the Abel-Jacobi map $AJ: Y^{(r)}\to\Pic Y$ sending $\bar D_1\to \cO_Y(\bar D_1)$, one has $AJ^*AE\,\iso\, E^{(r)}$ canonically. 
 
 Let $^r\phi: X^{(r)}\to Y^{(2r)}$ be the map sending $D'$ to $\phi^*D'$. In the ramified case one may define $N(E)$ by $N(E)=(^1\phi)^*E^{(2)}$. Then for any $r\ge 1$ one gets $(^r\phi)^*E^{(2r)}\,\iso\, (N(E))^{(r)}$ naturally. 

\begin{Lm} 
\label{Lm_irreducibility_over_cS_PD}
i) For a rank one local system $E$ on $Y$ one has naturally
$$
j_r^*\gp_{\cS !}\gq_{\cS}^*AE\,\iso\, (AN(E))_{\Omega(D)}\otimes
AN(E)^{-1}\boxtimes (\phi_*E^*)^{(r)}
$$
ii) Assume that $E$ does not descend to $X$. Then 
$\gp_{\cS !}\gq_{\cS}^*AE$ is the extension by zero under $j_r$. Besides, we have the cleanness
\begin{equation}
\label{iso_ext_by_zero_is_by_star}
(j_r)_{!}j_r^*\gp_{\cS !}\gq_{\cS}^*AE\,\iso\, (j_r)_{*}j_r^*\gp_{\cS !}\gq_{\cS}^*AE
\end{equation}
over $^0\cS_{P,D}$, and $\gp_{\cS !}\gq_{\cS}^*AE[\dim \cV_{P,D,H}]$ is an irreducible perverse sheaf over any connected component of $^0\cS_{P,D}$.

\smallskip\noindent
iii) $\gp_{\cS !}\Qlb[\dim\cV_{P,D,H}]$ is a perverse sheaf 
over $^0\cS_{P,D}$, the intermediate extension under $j_r: \Bun_1\times X^{(r)}\to \cS_{P,D}$ over the corresponding connected components of $^0\cS_{P,D}$.
\end{Lm}
\begin{Prf}
i) This follows from (\ref{diag_j_r_restriction}) and the fact that $(AE)_{\phi^*L_1}\,\iso\, (AN(E))_{L_1}$ naturally for $L_1\in \Bun_1$. We also used the fact that for any local system $V$ on $Y$ one has $\phi^r_*(V^{(r)})\,\iso\, (\phi_*V)^{(r)}$ canonically. 

\smallskip\noindent
ii) The first claim follows from $N_! (AE)=0$ for $N: \Pic Y\to \Pic X$. Since $^0\cV_{P,D,H}$ is smooth, (\ref{iso_ext_by_zero_is_by_star}) follows. To finish the proof, note that $(\phi_*E^*)^{(r)}[r]$ is an irreducible perverse sheaf on $X^{(r)}$ for any $r\ge 0$.

\smallskip\noindent
iii) The stack $^0\cV_{P,D,H}$ is smooth, and the map $\gp_{\cS}: {^0\cV_{P,D,H}}\to {^0\cS_{P,D}}$ is small. 
\end{Prf}

\medskip\noindent
\begin{Prf}\select{of Proposition~\ref{Pp_mainresults_2} i)} 
The map $\tilde\nu_{P,D}: {^0\Bun_{P,D}}\to \Bunt_{2,D}$ is smooth and surjective. By Lemmas~\ref{Lm_restriction_to_Bun_PD} and \ref{Lm_irreducibility_over_cS_PD}, $\tilde\nu_{P,D}^*F_G(AE[\dim\Pic Y])$ is an irreducible perverse sheaf over each connected component of $^0\Bun_{P,D}$. Our claim follows.
\end{Prf}

\ssec{Some ramified Eisenstein series}

\sssec{}
Let $\Bunb^{P}_D$ be the stack classifying $(\bar M\subset M)\in \Bun_{2,D}$ and a subsheaf $M_1\subset \bar M$ with $M_1\in\Bun_1$. We have a projection $\gp_D: \Bunb^{P}_D\to\Bun_{2,D}$ forgetting $M_1$. Let $\gp_{\tilde D}: \Bunb^P_{\tilde D}\to \Bunt_{2,D}$ be obtained from $\gp_D$ by the base change $\Bunt_{2,D}\to \Bun_{2,D}$. 

 Let $_0\Bunb^P_D\subset \Bunb^P_D$ be the open substack given by the property that the inclusion $M_1\hook{} M$ is maximal at all points of $D$. Let $j^P: {_0\Bunb^P_{\tilde D}}\hook{} \Bunb^P_{\tilde D}$ be defined by the same property. We get a diagram
$$
\begin{array}{cc}
_0\Bunb^P_{\tilde D}& \hook{j^P} \;\Bunb^P_{\tilde D} \;\toup{\gp_{\tilde D}}\; \Bunt_{2,D}\\
\downarrow\lefteqn{\scriptstyle \gq_{\tilde D}}\\ 
\Bunt_1,
\end{array}
$$
where $\gq_{\tilde D}$ sends $(M_1\subset \bar M\subset M, \cU)$ with an isomorphism (\ref{iso_cU_squared_for_Bunt_2D}) to $(L_1, \cU, \eta)$, where $L_1=M_1^*\otimes\cA$ and $\eta$ is the induced isomorphism (\ref{iso_eta_for_Bunt_1}). 

 Let $_0\gp_{\tilde D}: {_0\Bunb^P_{\tilde D}}\to \Bunt_{2,D}$ be the restriction of $\gp_{\tilde D}$.
\begin{Def} 
\label{Def_3.4.2}
The Eisenstein series on $\Bunt_{2,D}$ corresponding to $\phi_*\Qlb$ is defined as
$$
\Eis(\Qlb\oplus\cE_0)=(_0\gp_{\tilde D})_!\gq_{\tilde D}^*A\cE_0[\dim \Bunb^P_D]
$$
\end{Def}

  We will show below in Lemma~\ref{Lm_Eis_series_is_clean} that $j^P_{!}(\gq_{\tilde D}^*A\cE_0)\,\iso\, j^P_{!*}(\gq_{\tilde D}^*A\cE_0)$ naturally. Since $\gp_{\tilde D}$ is proper, $\Eis(\Qlb\oplus\cE_0)$ is Verdier self-dual. 
  
\begin{proof}[Proof of Proposition~\ref{Pp_mainresults_2} ii)]  Let $^0\Bunt_{P,D}={^0\Bun_{P,D}}\times_{\Bun_{2,D}}\Bunt_{2,D}$, write $\nu_{P,\tilde D}: {^0\Bunt_{P,D}}\to\Bunt_{2,D}$ for the projection to the second factor. Since the intersection of any connected component of $^0\Bunt_{P,D}$ with a fibre of $\nu_{P,\tilde D}$ is connected, by \cite[Lemma~4.8]{G} it suffices to establish an isomorphism
$$
\nu_{P,\tilde D}^*\Eis(\Qlb\oplus \cE_0)\,\iso\, \nu_{P,\tilde D}^*F_G(\Qlb[\dim\Pic Y])
$$
Since the 2-automorphism $\epsilon$ acts on both sides as $-1$, it suffices to
establish an isomorphism of shifted perverse sheaves
\begin{equation}
\label{iso_for_Pp_mainresults_2}
\Four_{\psi}^{-1}(^0\tilde\nu_{P,D})^*\Eis(\Qlb\oplus \cE_0)\,\iso\, \Four_{\psi}^{-1}(^0\tilde\nu_{P,D})^*F_G(\Qlb[\dim\Pic Y])
\end{equation}
over $^0\cS_{P,D}$. Let us do this.

Consider the stack $^0\Bun_{P,D}\times_{\Bunt_{2,D}} {_0\Bunb^P_{\tilde D}}$, its point is given by $L_1, \cA, M_1\in \Bun_1$, an exact sequence (\ref{seq_L*otimescA_by_L(-D)}), an inclusion $M_1\hook{} \bar M$ such that the induced map $M_1\to L_1^*\otimes\cA$ is also an inclusion. Indeed, if we had $M_1\subset L_1(-D)$ then $M_1$ would not be a subbundle of $M$ in a neighbourhood of $D$. Write $\cZ^r$ for the substack of the latter stack given by $\deg(L_1^*\otimes\cA)-\deg(M_1)=r$.

 Then $\cZ^r$ is the stack classifying $L_1, \cA, \in \Bun_1$, $D_1\in (X-D)^{(r)}$, and an exact sequence $0\to L_1(-D)\to ?\to L_1^*\otimes\cA\mid_{D_1}\to 0$. In this notation $M_1=L_1^*\otimes\cA(-D_1)$. The composition 
$$
\cZ^r\toup{\pr_2} {_0\Bunb^P_{\tilde D}}\;\toup{\gq_{\tilde D}} \;\Bunt_1
$$
sends the above point to $(L_1(D_1), \cU,\eta)$, where $\cU$ is given by (\ref{def_cU_Section_233}) with the induced isomorphism $\eta$. Let $s_1: \cZ^r\to (X-D)^{(r)}$ be the map sending the above point to $D_1$. By the multiplicativity of $A\cE_0$, one gets over $\cZ^r$ an isomorphism
$$
\pr_2^*\gq_{\tilde D}^* A\cE_0\,\iso\, s_1^* (\cE_0)^{(r)}
$$
 
 Let $\cY^r$ be the stack classifying $D_1\in (X-D)^{(r)}$, $L_1, \cA\in \Bun_1$ satisfying (C), and a section $v_1: L_1^2\to \cA\otimes\Omega(D-D_1)$. Consider the diagram 
$$
(X-D)^{(r)}\getsup{q_{\cY}} \cY^r\toup{p_{\cY}} {^0\cS_{P,D}},
$$ 
where $q_{\cY}$ is the projection sending the above point to $D_1$. The map
$p_{\cY}$ sends the above point to $(L_1,\cA, v)$, where $v$ is the composition 
$$
L_1^2\toup{v_1} \cA\otimes\Omega(D-D_1)\hook{} \cA\otimes\Omega(D)
$$ 
Then $p_{\cY !} q_{\cY}^* (\cE_0)^{(r)}[\dim\cY^r]$ is the contribution of $\cZ^r$ to 
$$
\Four_{\psi}^{-1}(^0\tilde\nu_{P,D})^*\Eis(\Qlb\oplus \cE_0)[\dimrel(^0\tilde\nu_{P,D})]
$$
over $^0\cS_{P,D}$. Let us show that $p_{\cY !} q_{\cY}^* (\cE_0)^{(r)}[\dim\cY^r]$ is perverse, the intermediate extension from the locus given by $v\ne 0$.
The map $p_{\cY}$ is proper, over the locus of $^0\cS_{P,D}$ given by $v\ne 0$ it is finite. 

One has $\dim\H^1(X, \cE_0)=2g-2+d$. So, $\RG((X-D)^{(r)}, (\cE_0)^{(r)})\,\iso\, \wedge^r \H^1(X, \cE_0)[-r]$ vanishes for $r>2g-2+d$. So, we assume $r\le 2g-2+d$. The codimension in $^0\cS_{P,D}$ of the locus given by $v=0$ equals
$-2\deg L_1+\deg \cA+d+(g-1)$. So, it suffices to show that
$$
r<-2\deg L_1+\deg \cA+d+(g-1)
$$
For $g\ge 1$ this follows from $2\deg L_1-\deg\cA+2g-2<0$. For $g=0$ this follows from $2\deg L_1-\deg\cA <1$ in (C). Applying Lemma~\ref{Lm_irreducibility_over_cS_PD}, one gets the desired isomorphism.
\end{proof}

\begin{Rem} 
\label{Rem_iota_very_first}
i) Since $\Qlb$ is invariant under the involution $\cB\mapsto \sigma^*\cB$ of $\Pic Y$, we get an involution on $F_G(\Qlb[\dim\Pic Y])$ that we denote by $\iota$. Its invariants on $\Eis(\Qlb\oplus\cE_0)$ is the contribution of those connected components of $\Bunb^P_{\tilde D}$ for which $\deg(M_1^*\otimes\cA)$ is even, that is, $\oplus_{r,d_1} \Eis(\Qlb\oplus\cE_0)^{r, d_1}$ for $r-d_1$ even.\\
ii) As a byproduct of our proof of Proposition~\ref{Pp_mainresults_2}, one may calculate the cohomology of the Prym variety $\Ker N$, where $N: \uPic Y\to \uPic X$ is the norm map, and $\uPic X$ denotes the Picard scheme of $X$. 
\end{Rem}

\begin{Lm} 
\label{Lm_Eis_series_is_clean}
The natural map $j^P_{!}(\gq_{\tilde D}^*A\cE_0)\,\iso\, j^P_{!*}(\gq_{\tilde D}^*A\cE_0)$ is an isomorphism.
\end{Lm}
\begin{Prf}
Consider the stack $\bar\cZ$ classifying $L_1,\cA\in\Bun_1$, an effective divisor $D_1$ on $X$, and an exact sequence $0\to L_1(-D)\to ?\to L_1^*\otimes\cA\mid_{D_1}\to 0$. A point of $\bar\cZ$ gives rise to a diagram
$$
\begin{array}{cccc}
0\to L_1(-D)\to \bar M & \to & L_1^*\otimes\cA & \to 0\\
& \nwarrow & \uparrow \\
&& M_1,
\end{array}
$$
where we have set $M_1=L_1^*\otimes\cA(-D_1)$. Let $^0\bar\cZ\subset\bar\cZ$ be the open substack  given by (C). We have a smooth and surjective map $^0\tilde\nu_{P,D}: {^0\bar\cZ}\to \Bunb^P_{\tilde D}$ sending the above point to $(M_1\subset \bar M\subset M, \cU)$ defined as above. 
The preimage of $_0\Bunb^P_{\tilde D}$ is the open substack $^0\cZ\subset {^0\bar\cZ}$ given by $D_1\cap D=\emptyset$. For $r\ge 0$ let $\bar \cZ^r$ be the locus given by $\deg D_1=r$. Then for any $r$ we have a cartesian square
$$
\begin{array}{ccc}
^0\cZ^r&\hook{} & {^0\bar\cZ^r}\\
\downarrow\lefteqn{\scriptstyle \zeta} && \downarrow\\
(X-D)^{(r)} & \hook{j_X} & X^{(r)},
\end{array}
$$
where the vertical maps are the smooth projections sending the corresponding point to $D_1$. The local system $(^0\tilde\nu_{P,D})^*\gq_{\tilde D}^*A\cE_0$ over $^0\cZ^r$ identifies with $\zeta^*(\cE_0)^{(r)}$. Since the natural map $(j_X)_!(\cE_0)^{(r)}\,\iso\, (j_X)_{!*}(\cE_0)^{(r)}$ is an isomorphism, our claim follows.
\end{Prf}

\ssec{Hecke functors}
\label{Sect_Hecke_functors_2.3.8}

\sssec{} 
We need the following Hecke functor on $\D(\Bunt_{2,D})$. Let $\Mod_{2,D}^r$ be the stack classifying $(\bar M\subset M)\in\Bun_{2,D}$ and an upper modification $M\subset M'$ with $\deg(M'/M)=r$ and $M'\cap \bar M(D)=M$. We have the diagram of projections
$$
\Bun_{2,D}\;\getsup{p_M}\;\Mod_{2,D}^r\;\toup{p'_M} \;\Bun_{2,D},
$$
where $p_M$ sends the above point to $(\bar M\subset M)$, and $p'_M$ sends the above point to $M'(-D)\subset \bar M'\subset M'$. Here $\bar M'/M'(-D)\subset M'/M'(-D)$ is the image of the natural map $\bar M/M(-D)\to M'/M'(-D)$. In other words, $\bar M'=\bar M+M'(-D)$. 

 We may also view $\Mod^r_{2,D}$ as the stack classifying $(\bar M'\subset M')\in\Bun_{2,D}$ together with a lower modification $\bar M\subset \bar M'$ such that $M'(-D)+\bar M=\bar M'$. For such a point we get $M=M'\cap \bar M(D)$. 
Both $p_M, p'_M$ are smooth of relative dimension $2r$. 

 Let $\tilde p_M: \Modt_{2,D}^r\to \Bunt_{2,D}$ be obtained from $p_M$ by the base change $\Bunt_{2,D}\to \Bun_{2,D}$. 
 
\begin{Lm} For $(\bar M\subset M\subset M')\in\Mod_{2,D}^r$ with $\bar D=\div(M'/M)$ one has a canonical $\ZZ/2\ZZ$-graded isomorphism
\begin{equation}
\label{iso_for_stack_Mod_2D}
\det\RG(M'/\bar M')\,\iso\, \det\RG(M/\bar M)\otimes \det\RG(X, \cO(\bar D)\mid_D)\otimes\det\RG(X, \cO_D)
\end{equation}
\end{Lm}
\begin{proof}
For the diagram
$$
\begin{array}{ccc}
\bar M' & \subset & M'\\
\cup && \cup\\
\bar M& \subset & M
\end{array}
$$
one gets a canonical $\ZZ/2\ZZ$-graded isomorphism
$$
\det\RG(X, \bar M'/\bar M)\otimes \det\RG(M'/\bar M')\,\iso\, \det\RG(M/\bar M)\otimes \det\RG(X, M'/M)
$$ 
The projection $M'(-D)\to \bar M'/\bar M$ induces an isomorphism $M'(-D)/M(-D)\,\iso\, \bar M'/\bar M$. Further, 
$$
\det\RG(X, (M'/M)\otimes\cO(-D))\,\iso\, \frac{\det\RG(X, M'/M)}{\det\RG(X, \cO(\bar D)\mid_D)\otimes\det\RG(X, \cO_D)}
$$ 
canonically. Indeed, by \cite[Lemma~1]{Ly1}, one has 
$$
\det\RG(X, M(-D))\,\iso\, \frac{\det\RG(X, M)}{\det\RG(X, \cA\mid_D)\otimes\det\RG(X, \cO_D)}
$$ 
and 
$$
\det\RG(X, M'(-D))\,\iso\, \frac{\det\RG(X, M')}{\det\RG(X, \cA(\bar D)\mid_D)\otimes\det\RG(X, \cO_D)}
$$
\end{proof}

\sssec{} 
\label{Sect_3.5.3_averaging}
We get the diagram
\begin{equation}
\label{diag_Modt_r_2D}
\Bunt_{2,D}\;\getsup{\tilde p_M}\; \Modt_{2,D}^r\; \toup{\tilde p'_M}\; \Bunt_{2,D},
\end{equation}
where $\tilde p'_M$ sends $(\bar M\subset M\subset M',\cU)$ and (\ref{iso_cU_squared_for_Bunt_2D}) to 
$(\bar M'\subset M', \cU')$, where 
$$
\cU'=\cU\otimes \det\RG(X,\cE_{\phi}(\bar D)/\cE_{\phi})^{-1}\otimes\det\RG(X, \cO(\bar D)/\cO),
$$
$\bar D=\div(M'/M)$, and $\cU'$ is equipped with the isomorphism
$
\cU'^2\,\iso\, \det\RG(M'/\bar M')\otimes\det\RG(X, \cO_D)
$ 
given by Lemma~\ref{Lm_detRG_L_mid_D_is_a_square} and (\ref{iso_for_stack_Mod_2D}). 

 Let $\Sh_0$ denote the stack of torsion sheaves on $X$. Let $\Sh^r_0$ be its connected component classifying torsion sheaves of length $r$. Let $\gs: \Modt_{2,D}^r\to \Sh_0$ be the map sending the above point to $M'/M$. For a local system $E$ on $X$ write $\cL_E$ for Laumon's sheaf on $\Sh_0$ defined in \cite[Section~1]{Ly3} (and originally in \cite{Laum}). Define the averaging functor 
$$
\Av^r_E: \D(\Bunt_{2,D})\to\D(\Bunt_{2,D})
$$ 
by $\Av^r_E(K)=(\tilde p'_M)_!(\tilde p_M^*K\otimes \gs^*\cL_E)[2r]$.

 In Section~\ref{Sect_2.1.5_should_be} we have introduced the full subcategory $\D_-(\Bunt_{2,D})\subset \D(\Bunt_{2,D})$. Clearly, $\Av^r_E$ preserves this subcategory. 
 
\sssec{}  Define the Hecke functor $\H_G: \D(\Bunt_{2,D})\to \D(X\times\Bunt_{2,D})$ as follows. Consider the diagram
$$
\Bunt_{2,D} \getsup{\tilde p_M} \Modt^1_{2,D} \toup{\supp\times\tilde p'_M} X\times \Bunt_{2,D},
$$
where $\supp$ sends a point of $\Modt^1_{2,D}$ as above to $\div(M'/M)$. 
For this diagram set 
\begin{equation}
\label{functor_H_G}
\H_G(K)=(\supp\times \tilde p'_M)_!(\tilde p_M)^*K[2]
\end{equation}
Over $(X-D)\times\Bunt_{2,D}$ one has $\DD\H_G(K)\,\iso\, \H_G(\DD(K))$ functorially for $K\in \D(\Bunt_{2,D})$. However, this is not clear over the whole of $X\times\Bunt_{2,D}$. 

\sssec{} For a local system $E$ on $X$ define the averaging functor 
$$
\Av^r_E: \D(\Pic Y)\to \D(\Pic Y)
$$ 
as follows. Consider the diagram $\Pic Y\getsup{\pr_2}Y^{(r)}\times \Pic Y\toup{m} \Pic Y$, where $m$ is the map $(D_1, \cB)\mapsto \cB(D_1)$, and $\pr_2$ is the projection. For $K\in\D(\Pic Y)$ set
$$
\Av^r_E(K)=\pr_{2 !}(\pr_1^*(\phi^*E)^{(r)}\otimes m^*K)[r]
$$  
This is consistent with the definition of the averaging functors from \cite[Section~4.4.2]{Ly1}.  

\ssec{} The purpose of this section is to derive Theorem~\ref{Th_1} from Theorem~\ref{Th_2}.  

\sssec{} 
\label{Sect_3.6.1_notations}
Set 
$$
\Modt_{2,D, H}^1=\Modt^1_{2,D}\times_{\Bunt_{2,D}} \Bunt_{2,D,H},
$$
where we used $\tilde p'_M: \Modt^1_{2,D}\to \Bunt_{2,D}$ to define the fibre product.

 Let $\cW_{D,H}$ be the stack classifying $x\in X, (\bar M\subset M)\in \Bun_{2,D}, \cB'\in\Pic Y$ together with an isomorphism $N(\cB')\,\iso\, \cC(D-x)$, where $\cC=\Omega\otimes\cA^{-1}$ and $\cA=\det M$. Set 
$$
\wt\cW_{D,H}=\cW_{D,H}\times_{\Bun_{2,D}}\Bunt_{2,D}
$$ 
 
  A point of $\Modt_{2,D,H}^1$ is given by a collection: $(\bar M\subset M\subset M',\cU)\in \Modt^1_{2,D}$, $\cB'\in\Pic Y$ equipped with $N(\cB')=\cC(D-x)$ with $x=\div(M'/M)$. Here $\cC=\Omega\otimes\cA^{-1}$ with $\cA=\det M$. 
 
 Consider the diagram
$$
\wt\cW_{D,H} \;\getsup{h_G}\; \Modt_{2,D,H}^1\;\toup{h'_G}\; \Bunt_{2, D,H},
$$
where $h'_G$ is the projection to the second factor, and $h_G$ sends the above point to $x, (\bar M\subset M,\cU)\in\Bunt_{2,D}, \cB'$ equipped with $N(\cB')\,\iso\, \cC(D-x)$. Define 
$$
\H_{GH}: \D(\Bunt_{2,D,H})\to \D(\wt\cW_{D,H})
$$ 
as
$\H_{GH}(K)=(h_G)_!(h'_G)^*K[2]$. 

\sssec{} Define also the functor 
$$
\H_{HG}: \D(\Bunt_{2,D,H})\to \D(\wt\cW_{D,H})
$$ 
as follows. Consider the diagram
$$
\wt\cW_{D,H}\getsup{h_H} Y\times_X \wt\cW_{D,H} \toup{h'_H} \Bunt_{2,D,H},
$$
where $h_H$ is the projection to the second factor, and $h'_H$ sends $y\in Y, (\bar M\subset M, \cU)\in\Bunt_{2,D}, \cB'$ together with $N(\cB')\,\iso\, \cC(D-\phi(y))$ to $(\bar M\subset M, \cU, \cB)$, where $\cB=\cB'(y)$. Here $\cB$ is equipped with the induced isomorphism $N(\cB)\,\iso\, \cC(D)$.
Set 
$$ 
\H_{HG}(K)=h_{H !} (h'_H)^*K[1]
$$
Since $h'_H$ is smooth of relative dimension 1, $\H_{HG}$ commutes with the Verdier duality. 
 
  Define the functor $\H_{H,X}: \D(\Pic Y)\to \D(X\times\Pic Y)$ by
$$
\H_{H, X}(K)=(\phi\times\id)_!m^*K[1] 
$$ 
for the diagram $X\times\Pic Y \getsup{\phi\times\id} Y\times\Pic Y\toup{m} \Pic Y$, here $m$ sends $(y,\cB)$ to $\cB(y)$. More generally, for $r\ge 1$ consider the equivariant derived category $ \D^{S_r}(X^r\times \Pic Y)$, where $S_r$ acts naturally on $X^r$ and trivially on $\Pic Y$. Define the functor 
$$
(\H_{H, X})^{\boxtimes r}: \D(\Pic Y)\to \D^{S_r}(X^r\times \Pic Y)
$$ 
by 
$$
(\H_{H, X})^{\boxtimes r}(K)=(\phi^r\times\id)_!(m^r)^*K[r]
$$
for the diagram $X^r\times\Pic Y \getsup{\phi^r\times\id} Y^r\times\Pic Y\toup{m^r}\Pic Y$. Here $m^r$ sends $(y_1,\ldots, y_r,\cB)$ to $\cB(y_1+\ldots+y_r)$. 
Consider the diagram of projections
\begin{equation}
\label{diag_projections_X^r_times_Pic Y}
X^r \;\getsup{\pr_1}\; X^r\times\Pic Y \;\toup{\pr_2} \;\Pic Y
\end{equation}

\begin{Lm} Let $E$ be a local system on $X$. For the diagram (\ref{diag_projections_X^r_times_Pic Y}) there is an  isomorphism 
$$
\Av^r_E(K)\,\iso\, \Hom_{S_r}(\triv, \pr_{2 !}(\pr_1^*E^{\boxtimes r}\otimes (\H_{H, X})^{\boxtimes r}(K)))
$$ 
functorial in $K\in \D(\Pic Y)$. \QED
\end{Lm}

 Write $\H^{\boxtimes r}_G$ for the $r$-fold iteration of (\ref{functor_H_G}). 

\begin{Lm} The functor $\H^{\boxtimes r}_G$ maps $\D(\Bunt_{2,D})$ to the equivariant derived category $\D^{S_r}(X^r\times\Bunt_{2,D})$. Moreover, there is an isomorphism 
$$
\Av^r_E(K)\,\iso\, \Hom_{S_r}(\triv, \pr_{2 !}(
\pr_1^*E^{\boxtimes r}\otimes\H^{\boxtimes r}_G(K)))
$$
functorial in $K\in \D(\Bunt_{2,D})$. Here the maps $\pr_i$ are the projections in the diagram $X^r\getsup{\pr_1} X^r\times\Bunt_{2,D}\toup{\pr_2}\Bunt_{2,D}$.
\end{Lm}
\begin{Prf} We argue essentially as in \cite[Propositions~1.8, 1.14 and 1.15]{G}, however we need to take into account the ramifications over $D$. 

\smallskip\noindent
{\bf Step 1} Let $\ItMod^r_{2,D}$ be the stack classifying $(\bar M\subset M\subset M')\in\Mod^r_{2,D}$ together with a complete flag $M=M_0\subset M_1\subset\ldots\subset M_r=M'$ with $\deg(M_i/M_{i-1})=1$ for all $i$. The shorthand $\ItMod$ refers to \select{iterated modifications}. Consider the diagram
$$
\begin{array}{ccccccccccccc}
\Mod^1_{2,D} &&&& \Mod^1_{2,D} &&&& \ldots &&&& \Mod^1_{2,D}\\
&\searrow\lefteqn{\scriptstyle p'_M}
 && \swarrow\lefteqn{\scriptstyle p_M} && \searrow\lefteqn{\scriptstyle p'_M} && \swarrow\lefteqn{\scriptstyle p_M} && \searrow\lefteqn{\scriptstyle p'_M} && \swarrow\lefteqn{\scriptstyle p_M}\\
&& \Bun_{2,D} &&&& \Bun_{2,D} &&&& \Bun_{2,D}
\end{array}
$$
Let us check that the limit of this diagram in the category of stacks identifies with $\ItMod^r_{2,D}$. Consider a diagram of locally free $\cO_X$-modules of rank 2 on $X$
$$
\begin{array}{cccccccc}
M_0 & \subset & M_1 &\subset &\ldots&\subset &M_r\\
\cup & & \cup &&&\ldots &\cup\\
\bar M_0 & \subset & \bar M_1 &\subset &\ldots&\subset &\bar M_r\\
\cup & & \cup &&&\ldots &\cup\\
M_0(-D) & \subset & M_1(-D) &\subset &\ldots&\subset &M_r(-D),
\end{array}
$$
such that $\div(M_i/\bar M_i)=D$ for all $i$. Then the condition $\bar M_r=\bar M_0+M_r(-D)$ is equivalent to the system of equations
$$
\bar M_i=\bar M_{i-1}+M_i(-D)\;\;\;\mbox{for all}\;\; 0<i\le r
$$
Indeed, the latter condition is equivalent to requiring that all the maps $\bar M_{i-1}/M_{i-1}(-D)\to \bar M_i/M_i(-D)$ are surjective. Since all these torsion sheaves are of the same length, this is equivalent to requiring the surjectivity of $\bar M_{0}/M_{0}(-D)\to \bar M_r/M_r(-D)$.

 For a point of $\ItMod^r_{2,D}$ as above we let $\bar M_i=\bar M+M_i(-D)$. Our claim follows.
 
\smallskip\noindent
{\bf Step 2} We have the projections
$$
\Bun_{2,D}\getsup{p_{I}} \ItMod^r_{2,D} \toup{p'_{I}}\Bun_{2,D},
$$
where $p_I$ sends a point of $\ItMod^r_{2,D}$ as above to $(\bar M\subset M)$, $p'_I$ sends it to $(\bar M'\subset M')$, which is the image of $(\bar M\subset M\subset M')$ under $p'_M$. Set 
$$
\wt\ItMod^r_{2,D}=\ItMod^r_{2,D}\times_{\Mod^r_{2,D}} \wt\Mod^r_{2,D}
$$ 
Define the diagram of projections
$$
\Bunt_{2,D}\getsup{\tilde p_{I}} \wt\ItMod^r_{2,D} \toup{\tilde p'_{I}}\Bunt_{2,D}
$$
by composing (\ref{diag_Modt_r_2D}) with $\ItMod^r_{2,D}\to \Mod^r_{2,D}$.

Let $\supp: \wt\ItMod^r_{2,D}\to X^r$ be the map sending the above point to the collection $\div(M_i/M_{i-1})$, $1\le i\le r$. From Step 1 it follows that  
$$
\H^{\boxtimes r}_G(K)\,\iso\, (\supp\times \tilde p'_I)_!(\tilde p_I^*K)[2r]
$$
functorially for $K\in\D(\Bunt_{2,D})$. Set 
$
\wt\Int^r_{2,D}=\Modt^r_{2,D}\times_{X^{(r)}} X^r,
$ 
where the map $\Modt^r_{2,D}\to X^{(r)}$ sends a point $(\bar M\subset M\subset M', \cU)$ to $\div(M'/M)$. The notations $\Int$ refers to \select{intermediate stack} and agrees with the notation of  \cite[Proposition~8]{G}. Let 
$$
\gr_{Int}: \wt\ItMod^r_{2,D}\to\wt\Int^r_{2,D},
$$
be the map whose first component is the projection $ \wt\ItMod^r_{2,D}\to \wt\Mod^r_{2,D}$, and the second component is $\supp$. The stack $\wt\ItMod^r_{2,D}$ is smooth. By \cite[Lemma~1.9]{G}, $\gr_{Int}$ is a small resolution of singularities. So, $(\gr_{Int})_!\IC$ is canonically isomorphic to the intersection cohomology sheaf on $\wt\Int^r_{2,D}$. This yields $S_r$-equivariant structure on $\H^{\boxtimes r}_G$. Recall that $\cL^r_E$ is the sheaf of $S_r$-invariants in the corresponding Springer sheaf. Now one concludes the proof as in \cite[Proposition~1.15]{G}. 
\end{Prf}
 
\sssec{} The following is an analog of \cite[Corollary~5]{Ly1}. 
 
\begin{Pp} There is a $S_r$-equivariant isomorphism of functors
\begin{equation}
\label{iso_wanted_for_Con2}
(\id\boxtimes F_H)\comp (\H_G)^{\boxtimes r}\,\iso\, (\H_{H, X})^{\boxtimes r}\comp F_H
\end{equation}  
from $\D_-(\Bunt_{2,D})$ to $\D(X^r\times\Pic Y)$. It is understood that the shift in the definition of $\id\boxtimes F_H: \D(X^r\times \Bunt_{2,D})\to \D(X^r\times\Pic Y)$ is the same as for $F_H$.
\end{Pp}
\begin{Prf}
1) Case $r=1$. Consider the diagram
$$
\begin{array}{cccccccc}
\Bunt_{2,D} &\getsup{\tilde p_M} & \Modt_{2,D}^1 & \gets & \Modt^1_{2,D,H} & \toup{h_G} & \wt\cW_{D,H} & \toup{p_{\tilde\cW}} \Bunt_{2,D}\\
 && \downarrow\lefteqn{\scriptstyle\supp\times \tilde p'_M} && \downarrow\lefteqn{\scriptstyle \supp\times h'_G}&& \downarrow\lefteqn{\scriptstyle q_{\tilde\cW}}\\
 && X\times\Bunt_{2,D} & \getsup{\id\times\gp} & X\times \Bunt_{2,D,H} & \toup{\id\times\gq} & X\times\Pic Y,
\end{array}
$$
where the left square is cartesian. The map $q_{\tilde\cW}$ sends a point $(x, \bar M\subset M, \cU, \cB')\in \wt\cW_{D,H}$ to $(x,\cB')$, and $p_{\tilde\cW}$ sends this point to $(\bar M\subset M, \cU)$.
 
  We get for $K\in \D(\Bunt_{2,D})$
$$
(\id\boxtimes F_H)\H_G(K)\,\iso\, q_{\tilde\cW !}(p_{\tilde\cW}^*K\otimes\H_{GH}(\Aut_{G,H}))[-\dim\Bunt_{2,D}]
$$
The commutative diagram 
$$
\begin{array}{ccc}
Y\times_X \wt\cW_{D,H} & \toup{h'_H} & \Bunt_{2,D,H}\\
\downarrow\lefteqn{\scriptstyle h_H} && \downarrow\lefteqn{\scriptstyle\gp}\\
\wt\cW_{D,H} & \toup{p_{\tilde\cW}} & \Bunt_{2,D},
\end{array}
$$
and Theorem~\ref{Th_2} show that the latter complex identifies with
$$
q_{\tilde\cW !}h_{H !}(h'_H)^*(\gp^*K\otimes\Aut_{G, H})[1-\dim\Bunt_{2,D}]\,\iso\, \H_{H, X}F_H(K)
$$
2) For any $r$ apply 1) $r$ times. One checks that the $S_r$-equivariant structures are identified on both sides.
\end{Prf} 

\begin{proof}[Proof of Theorem~\ref{Th_1}] Consider the diagram of projections
(\ref{diag_projections_X^r_times_Pic Y}). Tensoring (\ref{iso_wanted_for_Con2}) by $\pr_1^*(E^{\boxtimes r})$ and taking the direct image $\pr_{2 !}$ and further the $S_r$-invariants, one gets the desired isomorphism. 
\end{proof}

\ssec{} In this section we prove Theorem~\ref{Th_2}. 

\sssec{} Let $^0\wt\cW_{D,H}\subset \wt\cW_{D,H}$ be the open substack given by $x\notin D$. The desired isomorphism (\ref{iso_for_Th_2}) over $^0\wt\cW_{D,H}$ can be established as in 
\cite[Theorem~1]{Ly1}, the argument from \cite[Theorem~4]{Ly5} can also be adopted to our situation. So, it suffices to prove the following.

\begin{Pp} 
\label{Pp_interm_extension_for_Th_2}
Both sides of (\ref{iso_for_Th_2}) are perverse sheaves, the intermediate extensions with respect to the open immersion $^0\wt\cW_{D,H}\hook{}\wt\cW_{D,H}$.
\end{Pp}

\sssec{} One has an isomorphism $Y\times \Bunt_{2,D,H}\,\iso\, Y\times_X \wt\cW_{D,H}$ sending $(y, \bar M\subset M, \cU,\cB)$ to $(y, \bar M\subset M, \cU,\cB')$, where $\cB'=\cB(-y)$ is equipped with the induced isomoprhism $N(\cB')\,\iso\, \cC(D-x)$, here $x=\phi(y)$ and $\cC=\Omega\otimes(\det M)^{-1}$. 
Now Proposition~\ref{Pp_mainresults_1} shows that $(h'_H)^*\Aut_{G,H}[1]$ is perverse, the intermediate extension from the open substack given by $\phi(y)\notin D$. Since $h_H$ is finite, $H_{HG}(\Aut_{G,H})$ is perverse, the intermediate extension under $^0\wt\cW_{D,H}\hook{}\wt\cW_{D,H}$.

\sssec{} Consider the stack $\Bun^P_D$ classifying $(\bar M\subset M)\in\Bun_{2,D}$ together with a subbundle $M_1\subset \bar M$ of rank one such that $M_1\subset M$ is also a subbundle. Let $\nu^P_D: \Bun^P_D\to \Bun_{2,D}$ be the map forgetting $M_1$.
We extend $\nu^P_D$ to a morphism $\tilde\nu^P_D: \Bun^P_D\to\Bunt_{2,D}$ 
sending the above point to $(\bar M\subset M, \cU)$, where
$$
\cU=\frac{\det\RG(X, M_1^*\otimes\cA)\otimes\det\RG(X, \cE_{\phi})}{\det\RG(X, \cO)\otimes\det\RG(X, M_1^*\otimes\cA\otimes\cE_{\phi})}
$$
is equipped with the isomorphism 
$$
\cU^2\,\iso\, \det\RG(X, M_1^*\otimes\cA\mid_D)\otimes\det\RG(X , \cO_D)\,\iso\, \det\RG(X, M/\bar M)\otimes\det\RG(X, \cO_D)
$$
given by Lemma~\ref{Lm_detRG_L_mid_D_is_a_square}. A point of $\Bun^P_D$ is given by $M_1,\cA\in\Bun_1$ and an exact sequence 
\begin{equation} 
\label{ext_M_1*otimescA_by_M_1}
0\to M_1\to M\to M_1^*\otimes\cA\to 0
\end{equation}
For this point $\bar M/M(-D)=M_1\mid_D$. 

\sssec{} Let $^0\Bun^P_D\subset\Bun^P_D$ be the open substack given by 
\begin{equation}
\label{def_of^0Bun^P_D} 
 2\deg M_1-\deg\cA<1-2g-d,
\end{equation} 
here $d=\deg(D)$ and $\cA=\det M$. This condition on the degrees implies $\H^1(X, M_1^{-2}\otimes \cA(-D))=0$. The restriction $^0\tilde\nu^P_D: {^0\Bun^P_D}\to\Bunt_{2,D}$ of $\tilde\nu^P_D$ is smooth and surjective. Set 
$$
\Bun^P_{D,H}=\Bun^P_D\times_{\Bun_{2,D}}\Bun_{2,D,H}\;\;\;\;\;\mbox{and}\;\;\;\;\; 
{^0\Bun^P_{D,H}}={^0\Bun^P_D}\times_{\Bun_{2,D}}\Bun_{2,D,H}
$$ 

\sssec{} Let $\tau^P: \Bun^P_{D,H}\to \Bun_{P_2}$ be the following map. For a point of $\Bun^P_{D,H}$ given by (\ref{ext_M_1*otimescA_by_M_1}) and $\cB\in\Pic Y$ with $N(\cB)\,\iso\,\cC(D)$ and $\cC=\Omega\otimes\cA^{-1}$ let $L=\phi_*\cB$ and $M'=\tau(\bar M\subset M, \cB)$. Then
$M_1\otimes L\subset M'$ is a lagrangian subbundle. The map $\tau^P$ sends the above point to $(M_1\otimes L\subset M')\in\Bun_{P_2}$. The corresponding exact sequence 
$$
0\to \Sym^2(M_1\otimes L)\to ?\to\Omega\to 0
$$ 
is the push-forward of $0\to M_1^2\otimes\cC\to M\otimes M_1\otimes\cC\to \Omega\to 0$ under
$\id\otimes s_c: M_1^2\otimes\cC\to M_1^2\otimes\Sym^2 L$. Here $s_c: \cC\to\Sym^2 L$ is the canonical section of the symmetric form $\Sym^2 L\to \cC(D)$. 
One checks that the following diagram is canonically 2-commutative
\begin{equation}
\label{diag^P_for_2.3.10.1}
\begin{array}{ccc}
\Bunt_{2,D,H} & \toup{\tilde\tau} & \Bunt_{G_2}\\
\uparrow\lefteqn{\scriptstyle \tilde\nu^P_H} && \uparrow\lefteqn{\scriptstyle \tilde\nu_{P_2}}\\
\Bun^P_{D,H} & \toup{\tau^P} & \Bun_{P_2},
\end{array}
\end{equation}
where we have set $\tilde\nu^P_H=\tilde\nu^P_D\times\id$. 

\sssec{} 
\label{Sect_3.7.7_should_be}
Let $\cZ^P$ be the stack classifying a point of $\Bun^P_{D,H}$ given by (\ref{ext_M_1*otimescA_by_M_1}), $\cB\in\Pic Y$ equipped with $N(\cB)\,\iso\, \cC(D)$ with $\cC=\Omega\otimes\cA^{-1}$, and a section $s: \phi^*M_1\to \cB^*\otimes\Omega_Y$. Let $\ev_{\cZ}: \cZ^P\to\A^1$ be the map sending this point to the pairing of $N(s): M_1^2\to \cA\otimes\Omega$ with (\ref{ext_M_1*otimescA_by_M_1}). Let $p_{\cZ}: \cZ^P\to \Bun^P_{D,H}$ be the map forgetting $s$. We will need the following.
 
\begin{Lm} 
\label{Lm_Aut_GH_over_Bun^P_DH}
There is an isomorphism 
$$
(\tilde\nu^P_H)^*\Aut_{G,H}[\dimrel(\tilde\nu^P_H)]\,\iso\, p_{\cZ !}\ev_{\cZ}^*\cL_{\psi}[a],
$$
where $a$ is a function of a connected component of $\Bun^P_{D,H}$ given by $a=\dim(\Bun^P_{D,H})+\chi(Y, \phi^*M_1^*\otimes \cB^*\otimes\Omega_Y)$ for a point of $\Bun^P_{D,H}$ as above. 
\end{Lm}
\begin{Prf} 
The section $s$ can be rewritten as $s: M_1\otimes L\to \Omega$, where $L=\phi_*\cB$. Combining the above description of $\tau^P$ with \cite[Proposition~1]{Ly1}, one gets the desired result.
\end{Prf}
 
\sssec{} Set $\cW^P=\wt\cW_{D,H}\times_{\Bunt_{2,D}}{^0\Bun^P_D}$, let $\tilde\nu_{\cW}: \cW^P\to \wt\cW_{D,H}$ be the projection to the first factor. To prove Proposition~\ref{Pp_interm_extension_for_Th_2}, we will consider the restrictions of both sides of (\ref{iso_for_Th_2}) under $\tilde\nu_{\cW}$. 

 Let $h^P: \Mod^P\to\cW^P$ be obtained from $h_G: \Modt^1_{2,D,H}\to \wt\cW_{D,H}$ by the base change $\tilde\nu_{\cW}: \cW^P\to \wt\cW_{D, H}$. 
 
 A point of $\cW^P$ is given by a collection: $x\in X$, $M_1,\cA\in\Bun_1$ satisfying (\ref{def_of^0Bun^P_D}), an exact sequence (\ref{ext_M_1*otimescA_by_M_1}), $\cB'\in\Pic Y$ equipped with $N(\cB')\,\iso\, \cC(D-x)$, here $\cC=\Omega\otimes\cA^{-1}$. 
 
 A point of $\Mod^P$ is given by a point of $\cW^P$ as above together with an upper modification $M\subset M'$ with $x=\div(M'/M)$ such that $M'\cap \bar M(D)=M$. The latter condition is equivalent to requiring that the inclusion $M_1\hook{} M'$ is maximal in a neighbourhood of $D$. 

\sssec{} Let $^0\Mod^P\subset \Mod^P$ be the open substack given by the property that $M_1\subset M'$ is a subbundle. Let $^1\Mod^P$ be its complement in $\Mod^P$. Let $^iK$ denote the contribution of $^i\Mod^P$ to the complex 
$$
\tilde\nu_{\cW}^*\H_{GH}(\Aut_{G,H})[\dimrel(\tilde\nu_{\cW})],
$$ 

 Let $h^{'P}: {^0\Mod^P}\to \Bun^P_{D,H}$ be the map sending a point of $^0\Mod^P$ as above to the collection: the exact sequence
\begin{equation}
\label{ext_M_1otimescA'_by_M_1} 
0\to M_1\to M'\to M_1^*\otimes \cA(x)\to 0,
\end{equation}
$\cB'$ equipped with $N(\cB')\,\iso\, \cC(D-x)$, where $\cC=\Omega\otimes\cA^{-1}$. We get a diagram
$$
\begin{array}{ccccc}
\cW^P & \getsup{^0h^P} & 
^0\Mod^P &\toup{h^{'P}} & \Bun^P_{D,H} \\
\downarrow\lefteqn{\scriptstyle \tilde\nu_{\cW}} && \downarrow && \downarrow\lefteqn{\scriptstyle \tilde\nu^P_H}\\
\wt\cW_{D,H} & \getsup{h_G} & \Modt^1_{2,D,H} & \toup{h'_G} & \Bunt_{2,D,H},
\end{array}
$$
where we denoted by $^0h^P$ the restriction of $h^P$. It is easy to check that the right square of this diagram is 2-commutative. By definition, 
\begin{equation}
\label{def_of_^0K}
^0K=(^0h^P)_!(\tilde\nu^P_H h^{'P})^*\Aut_{G,H}[\dimrel(\tilde\nu^P_H h^{'P})]
\end{equation}

\sssec{}  Let $\cY\cW^P$ be the stack classifying a point of $\cW^P$ as above together with $s: \phi^*M_1\to \cB'^*\otimes\Omega_Y$ such that $N(s): M_1^2\to \cA\otimes\Omega(x)$ actually lies in $\H^0(X, M_1^{-2}\otimes\cA\otimes\Omega)$. 
 
 Let $\ev_{\cY\cW}: \cY\cW^P\to\A^1$ be the map sending the above point to the pairing of $N(s)$ with (\ref{ext_M_1*otimescA_by_M_1}). Let $p_{\cY\cW}: \cY\cW^P\to \cW^P$ be the map forgetting $s$.

\begin{Lm} 
\label{Lm_first_description_^0K}
There is an isomorphism over $\cW^P$
$$
^0K\,\iso\, (p_{\cY\cW})_!\ev_{\cY\cW}^*\cL_{\psi}[b],
$$ 
where $b$ is a function of a connected component of $\cW^P$ given by 
$$
b=\chi(Y, \phi^*M_1^*\otimes\cB'^*\otimes\Omega_Y)-1+\dim(\cW^P)$$
\end{Lm} 
\begin{Prf}
Let $f: {^0\cY}\to {^0\Mod^P}$ be the stack classifying a point of $^0\Mod^P$ as above together with a section $s: \phi^*M_1\to \cB'^*\otimes\Omega_Y$. Let $^0\ev: {^0\cY}\to \A^1$ be the map sending a point of $^0\cY$ to the pairing of (\ref{ext_M_1otimescA'_by_M_1}) with $N(s): M_1^2\to \Omega\otimes\cA(x)$. 
By Lemma~\ref{Lm_Aut_GH_over_Bun^P_DH} combined with (\ref{def_of_^0K}), the complex $^0K$ rewites (up to a shift) as the direct image with compact support of $^0\ev^*\cL_{\psi}$ under the composition 
$$
^0\cY\toup{f} {^0\Mod^P}\toup{^0h^P}\cW^P
$$

 Let $i: {^0\cY'}\hook{}{^0\cY}$ be the closed substack given by requiring 
$N(s)\in \H^0(X, M_1^{-2}\otimes\cA\otimes\Omega)$. Let us show that the natural map
\begin{equation}
\label{map_for_Lm_over_cW^P}
(^0h^P)_!f_!(^0\ev)^*\cL_{\psi}\to 
(^0h^P)_!f_!i_*i^*(^0\ev)^*\cL_{\psi}
\end{equation}
is an isomorphism. Let $\eta$ be a $k$-point of $\cW^P$, write $Z_{\eta}$ for the fibre of $^0\Mod^P\to \cW^P$ over $\eta$. Note that $(M_1^2\otimes\cA^*)_x$ acts freely and transitively on $Z_{\eta}$. Since we assume $d>0$, for a point of $^0\Bun^P_D$ we get $2\deg M_1-\deg\cA<0$, so that $\H^0(X, M_1^2\otimes\cA^*)=0$, and the sequence is exact
$$
0\to (M_1^2\otimes\cA^*)_x\toup{\xi} \H^1(X, M_1^2\otimes\cA^*(-x))\to \H^1(X, M_1^2\otimes\cA^*)\to 0
$$
Under the action of $v\in (M_1^2\otimes\cA^*)_x$ on $Z_{\eta}$ the class of the exact sequence (\ref{ext_M_1otimescA'_by_M_1}) changes by $\xi(v)$.
Integrating over $Z_{\eta}$ with $s$ fixed first, we will get zero unless  $N(s)\in \H^0(X, M_1^{-2}\otimes\cA\otimes\Omega)$. So, (\ref{map_for_Lm_over_cW^P}) is an isomorphism. Our claim easily follows.
\end{Prf} 
 
\begin{Rem} 
\label{Rem_exact_sequence_over_cYcW^P}
Let $\cY$ be the stack classifying $x\in X$, $\cB\in\Pic Y$ with $s: \cO\to \cB$ such that $N(s)\in \H^0(X, N(\cB))$ lies in $\H^0(X, N(\cB)(-x))$. Let $\cX$ be the stack classifying $\cB\in \Pic Y$, $y\in Y$ and $s: \cO\to \cB(-y)$. Let $\pi_{\cX}: \cX\to \cY$ be the map sending the above point to $(x,\cB, s)$, where $x=\phi(y)$. Let $\kappa: \cY'\to \cY$ be the locally closed substack given by the property that $x\notin D$, and $s: \cO\to \cB(-\phi^*(x))$ is regular. Let $\pr': \cY'\to X$ be the projection sending the above point to $x$. Then one has an exact sequence on $\cY$
$$
0\to \Qlb\to \pi_{\cX !}\Qlb\to \kappa_!\pr'^*\cE_0\to 0
$$
\end{Rem}
 
\sssec{} Let $\kappa: {'\cY\cW^P}\subset \cY\cW^P$ be the locally closed substack given by the property that $s: \phi^*M_1\to \cB'^*\otimes\Omega_Y(-\phi^*(x))$ is regular and $x\notin D$. Let $q_{\cY\cW}: {'\cY\cW^P}\to X$ be the map sending the above point to $x$. 

 Let $^0\cZ^P$ denote the preimage of $^0\Bun^P_{D,H}$ under the map $p_{\cZ}: \cZ^P\to \Bun^P_{D,H}$ defined in Section~\ref{Sect_3.7.7_should_be}. 
 
\sssec{}  A point of $Y\times{^0\cZ^P}$ is given by $y\in Y$, an exact sequence (\ref{ext_M_1*otimescA_by_M_1}) satisfying (\ref{def_of^0Bun^P_D}), $\cB\in\Pic Y$ equipped with $N(\cB)\,\iso\, \cC(D)$, and $s: \phi^*M_1\to \cB^*\otimes\Omega_Y$. Let $h_{\cZ}: Y\times{^0\cZ^P}\to \cY\cW^P$ be the map sending a point as above to the same data with $\cB'=\cB(-y)$. The following diagram commutes
$$
\begin{array}{ccccc}
\cY\cW^P && \xleftarrow{h_{\cZ}} && Y\times {^0\cZ^P}\\
\downarrow\lefteqn{\scriptstyle p_{\cY\cW}} &&&& \downarrow\lefteqn{\scriptstyle \id\times p_{\cZ}}\\
\cW^P & \getsup{h_{\cW}} & Y\times_X \cW^P & \xleftarrow[\sim]{\xi^P} & Y\times {^0\Bun^P_{D,H}}\\
\downarrow\lefteqn{\scriptstyle \tilde\nu_{\cW}} && \downarrow\lefteqn{\scriptstyle \id\times\tilde\nu_{\cW}}&& \downarrow\lefteqn{\scriptstyle  \tilde\nu^P_H\comp\pr_2}\\
\wt\cW_{D,H} & \getsup{h_H} & Y\times_X \wt\cW_{D,H} & \toup{h'_H} & \Bunt_{2,D,H} 
\end{array}
$$
where $h_{\cW}$ is the projection. We denoted by
$\xi^P$ the isomorphism sending a collection $y\in Y$, (\ref{ext_M_1*otimescA_by_M_1}) and $\cB$ with $N(\cB)\,\iso\, \cC(D)$ to the collection (\ref{ext_M_1*otimescA_by_M_1}), $y$, $\cB'$, where $\cB'=\cB(-y)$. 

\sssec{} By Remark~\ref{Rem_exact_sequence_over_cYcW^P}, we get an exact sequence on $\cY\cW^P$
\begin{equation}
\label{seq_on_cYcW^P}
0\to \ev_{\cY\cW}^*\cL_{\psi}\to h_{\cZ !}(\Qlb\boxtimes\ev_{\cZ}^*\cL_{\psi})\to \kappa_!(q_{\cY\cW}^*\cE_0\otimes\kappa^*\ev_{\cY\cW}^*\cL_{\psi})\to 0
\end{equation}
Applying $(p_{\cY\cW})_!$ to (\ref{seq_on_cYcW^P}) and using Lemmas~\ref{Lm_Aut_GH_over_Bun^P_DH} and \ref{Lm_first_description_^0K}, we get an exact triangle
\begin{equation}
\label{triangle_first_for_^0K}
^0K\to \tilde\nu_{\cW}^*\H_{HG}(\Aut_{G,H})[\dimrel(\tilde\nu_{\cW})] \to \cK,
\end{equation}
where we have set 
$$
\cK=(p_{\cY\cW})_!\kappa_!(q_{\cY\cW}^*\cE_0\otimes\kappa^*\ev_{\cY\cW}^*\cL_{\psi})
[b],
$$ 
here $b$ is the function defined in Lemma~\ref{Lm_first_description_^0K}. 

\sssec{} A point of $^1\Mod^P$ yields an exact sequence
\begin{equation}
\label{seq_M_1^*otimescA_by_M_1(x)}
0\to M_1(x)\to M'\to M_1^*\otimes\cA\to 0,
\end{equation} 
which is the push-forward of (\ref{ext_M_1*otimescA_by_M_1}) under $M_1\to M_1(x)$. Let $h''^P: {^1\Mod^P}\to \Bun^P_{D,H}$ be the map sending a point of 
$^1\Mod^P$ as above to the collection: (\ref{seq_M_1^*otimescA_by_M_1(x)}) and $\cB'$ equipped with $N(\cB')\,\iso\, \cC(D-x)$. We get a diagram
$$
\begin{array}{ccccc}
\cW^P & \getsup{^1h^P} &
^1\Mod^P & \toup{h''^P} & \Bun^P_{D,H}\\
\downarrow\lefteqn{\scriptstyle \tilde\nu_{\cW}} &&
\downarrow\lefteqn{\scriptstyle {^1p}} && \downarrow\lefteqn{\scriptstyle \tilde\nu^P_H}\\
\wt\cW_{D,H} & \getsup{h_G} & \Modt^1_{2,D,H} & \toup{h'_G} & \Bunt_{2,D,H},
\end{array}
$$ 
where $^1p$ is the composition $^1\Mod^P\hook{}\Mod^P\to \Modt^1_{2,D,H}$. 
We denote by $^1h^P$ the composition $^1\Mod^P\hook{}\Mod^P\toup{h^P} \cW^P$. The map $^1h^P$ is an open immersion given by $x\notin D$.
The right square in the above diagram is not 2-commutative. But it is almost 2-commutative, in the sense that there is an isomorphism
\begin{equation}
\label{iso_for_non2-commutative_diag_Th2}
^1p^*(h'_G)^*\Aut_{G,H}\,\iso\, \pr_{\Mod}^*\cE_0\otimes (h''^P)^*(\tilde\nu^P_H)^*\Aut_{G,H},
\end{equation}
where $\pr_{\Mod}:{^1\Mod^P}\to X$ is the map sending a point of $^1\Mod^P$ as above to $x$. 

\begin{Lm} There is an isomorphism $\cK\,\iso\, {^1K}$ on $\cW^P$. 
\end{Lm}
\begin{Prf}
Using (\ref{iso_for_non2-commutative_diag_Th2}), one gets 
$$
^1K\,\iso\, (^1h^P)_!(\pr^*_{\Mod}\cE_0\otimes (h''^P)^*(\tilde\nu^P_H)^*\Aut_{G,H})[1+\dimrel(\tilde\nu^P_H\comp h''^P)]
$$
Applying Lemma~\ref{Lm_Aut_GH_over_Bun^P_DH}, one gets the desired isomorphism.
\end{Prf}
 
\begin{Lm} The map $h'_G\comp(^1p)$ is smooth, so $^1K[-1]$ is a perverse sheaf on $^1\Mod^P$. Its extension by zero under the map $^1h^P: {^1\Mod^P}\hook{} \cW^P$ is also perverse.
\end{Lm}
\begin{Prf} 
Given a point $(\bar M'\subset M', \cB', N(\cB')\,\iso\, \cC'(D))$ of $\Bunt_{2,D,\tilde H}$, one first picks a line subbundle $M_1(x)\subset \bar M'$ such that $M_1(x)\subset M'$ is also a subbundle. This is a smooth map because $2\deg M_1(x)-\deg \cA'<2-2g-d$ with $\cA'=\det M'$, and so $\H^1(X, (M_1(x))^{-2}\otimes \cA'(-D))=0$. Then one picks any $x\in X-D$. Finally, the exact sequence $0\to M_1(x)/M_1\to M'/M_1\to M_1^*\otimes\cA\to 0$ splits, and a choice of a splitting is also a smooth map.

 The last claim follows from the fact that the inclusion $X-D\hook{} X$ is affine.
\end{Prf}

\sssec{} Consider the distinguished triangle over $\cW^P$
\begin{equation}
\label{triangle_for_HG}
^1K[-1]\toup{\alpha} {^0K}\to \tilde\nu_{\cW}^*\H_{HG}(\Aut_{G,H})[\dimrel(\tilde\nu_{\cW})]
\end{equation}
obtained from (\ref{triangle_first_for_^0K}). Since the RHS of (\ref{triangle_for_HG}) is perverse, $^0K$ is also perverse over $\cW^P$. So, (\ref{triangle_for_HG}) is an exact sequence of perverse sheaves.

 By definition of $^iK$, we also have a distinguished triangle
\begin{equation}
\label{triangle _for_GH}
^1K[-1]\toup{\beta} {^0K}\to \tilde\nu_{\cW}^*\H_{GH}(\Aut_{G,H})[\dimrel(\tilde\nu_{\cW})]
\end{equation}

 Over $^0\wt\cW_{D,H}$ the map $h_G: \Mod^1_{2,D,H}\to \wt\cW_{D,H}$ is proper, so $\H_{GH}(\Aut_{G,H})$ is Verdier self-dual over $^0\wt\cW_{D,H}$. From (\ref{triangle _for_GH}) we see that $\H_{GH}(\Aut_{G,H})$ is perverse 
over $^0\wt\cW_{D,H}$. 

\begin{Lm} 
\label{Lm_irreducibilities_parts_of_^1K[-1]}
1) Over $\cW^P$ one has 
$$
\Hom(^1K[-1], \tilde\nu_{\cW}^*\H_{HG}(\Aut_{G,H})[\dimrel(\tilde\nu_{\cW})])=0
$$ 
in the category of perverse sheaves.\\
2) Both $S_2$-invariants and anti-invariants in $^1K[-1]$ are irreducible perverse sheaves on $^1\Mod^P$. (They correspond to the contributions of $\Aut_g$ and $\Aut_s$ to  $^1K[-1]$). 
\end{Lm}
\begin{Prf} 1) By adjointness for $^1h^P$, we have to prove this over $^1\Mod^P$. 
Recall that we have a diagram
$$
\begin{array}{ccccc}
\cW^P & \getsup{h_{\cW}} & Y\times_X \cW^P\\
\downarrow\lefteqn{\scriptstyle \tilde\nu_{\cW}} && \downarrow\lefteqn{\scriptstyle \id\times\tilde\nu_{\cW}}\\
\wt\cW_{D,H} & \getsup{h_H} & Y\times_X \wt\cW_{D,H} & \toup{h'_H} & \Bunt_{2,D,H}
\end{array}
$$
By adjointness, we have to show that over $(Y-D_Y)\times_X \cW^P$ one has
$$
\Hom(h_{\cW}^*(^1K)[-1], (\id\times\tilde\nu_{\cW})^*(h'_H)^*\Aut_{G,H}[1+\dimrel(\tilde\nu_{\cW})])=0
$$
Let $\cS$ be the stack classifying $y\in Y-D_Y$ with $x=\phi(y)$, $M_1, \cA\in \Pic X$ satisfying (\ref{def_of^0Bun^P_D}), $\cB'\in\Pic Y$ equipped with $N(\cB')\,\iso\, \cC(D-x)$, here $\cC=\Omega\otimes\cA^{-1}$. 

 Let $Y_{\cS}$ be the stack classifying a point of $\cS$ as above together with a section $\bar s: M_1^2\to \cA\otimes\Omega$. Then $Y_{\cS}$ and $(Y-D_Y)\times_X \cW^P$ are dual generalized vector bundles over $\cS$.
Both perverse sheaves are Fourier transforms under $\Four_{\psi}: \D(Y_{\cS})\to \D((Y-D_Y)\times_X \cW^P)$.

 Let $\bar Y_{\cS}$ be the stack classifying a point of $\cS$ as above together with a section $s: \phi^*M_1\to \cB'^*\otimes\Omega_Y(-y)$. We have a projection $\delta: \bar Y_{\cS}\to Y_{\cS}$ given by $\bar s=N(s)$. Let also $i: \bar Y'_{\cS}\subset \bar Y_{\cS}$ be the closed substack given by the condition that $s: \phi^*M_1\to \cB'^*\otimes\Omega_Y(-\phi^*(x))$ is regular.

For a point of $\cS$ one has an inclusion
$$
\H^0(Y, \phi^*M_1^*\otimes \cB'^*\otimes \Omega_Y(-\phi^*(x)))\to \H^0(Y, \phi^*M_1^*\otimes \cB'^*\otimes \Omega_Y(-y)),
$$
whose cokernel is a 1-dimensional space. Indeed, $\H^1(Y, \phi^*M_1^*\otimes \cB'^*\otimes \Omega_Y(-\phi^*(x)))=0$ because of (\ref{def_of^0Bun^P_D}). 
The morphism $\delta$ is finite. 

 It suffices to show that 
$$
\Hom_{Y_{\cS}}(\delta_! i^*\Qlb[\dim\bar Y_{\cS}-1],
\delta_!\Qlb[\dim\bar Y_{\cS}])=0
$$ 
in the category of perverse sheaves on $Y_{\cS}$. The condition (\ref{def_of^0Bun^P_D}) implies $\deg(M_1^{-2}\otimes\cA\otimes\Omega)>2g_Y-2$. So, by Lemma~\ref{Lm_perversity_for_phi'_ex}, both $S_2$-invariants and anti-invariants in $\delta_!\Qlb[\dim\bar Y_{\cS}]$ are perverse sheaves irreducible over each connected component.
Let $Y'_{\cS}\subset Y_{\cS}$ be the closed substack given by $\bar s\in \H^0(X, M_1^{-2}\otimes\cA\otimes\Omega(-x))$. Then $\delta_! i^*\Qlb$ is the extension by zero from $Y'_{\cS}$. The restriction to $Y_{\cS}-Y'_{\cS}$ of $S_2$-invariants or anti-invariants in $\delta_!\Qlb$ does not vanish. Our claim follows.

\medskip\noindent
2) Let $\bar\cS$ be the stack classifying $x\in X-D$, $M_1,\cA\in\Pic X$ satisfying (\ref{def_of^0Bun^P_D}), $\cB'\in\Pic Y$ equipped with $N(\cB')\,\iso\, \cC(D-x)$, here $\cC=\Omega\otimes\cA^{-1}$. 
 
  Let $Y_{\bar\cS}$ be the stack classifying a point of $\bar\cS$ as above together with $\bar s: M_1^2\to \cA\otimes\Omega(-x)$. Let $\bar Y_{\bar\cS}$ be the stack classifying a point of $\bar\cS$ as above together with $s: \phi^*M_1\to \cB'^*\otimes\Omega_Y(-\phi^*(x))$. We have a finite map $\bar\delta: \bar Y_{\bar\cS}\to Y_{\bar\cS}$ given by $\bar s=N(s)$. 

For a point of $\bar\cS$ the condition (\ref{def_of^0Bun^P_D}) implies $\deg(M_1^{-2}\otimes\cA\otimes\Omega(-x))>2g_Y-2$. By Lemma~\ref{Lm_perversity_for_phi'_ex}, both $S_2$-invariants and anti-invariants in $\bar\delta_!\Qlb[\dim\bar Y_{\bar\cS}]$ are perverse sheaves irreducible over each connected component. 

 Let $\cV$ be the stack classifying a point of $\bar\cS$ as above together with an exact sequence $0\to M_1(x)\to ?\to M_1^*\otimes\cA\to 0$ on $X$. Then $\cV$ and $Y_{\bar\cS}$ are dual generalized vector bundles over $\cS$. 
 
 We have a moprhism $f: {^1\Mod^P}\to \cV$ over $\bar\cS$ taking (\ref{ext_M_1*otimescA_by_M_1}) to its push-forward via $M_1\hook{} M_1(x)$. The map $f$ is a vector bundle of rank one, in particular its smooth with connected fibres. So, both $S_2$-invariants and anti-invariants in 
$$
^1K[-1]\,\iso\, f^*\Four_{\psi}(\bar\delta_!\Qlb[\dim\bar Y_{\bar\cS}])[1]
$$ 
are perverse sheaves irreducible over each connected component of $^1\Mod^P$. 
\end{Prf}

\sssec{} By 1) of Lemma~\ref{Lm_irreducibilities_parts_of_^1K[-1]}, there is a morphism $\gamma: {^1K}[-1]\to {^1K}[-1]$ over $\cW^P$ such that $\beta=\alpha\comp\gamma$. If the kernel of $\gamma$ is not zero then $\H_{GH}(\Aut_{G,H})$ would have a nontrivial perverse cohomology sheaf in degree $-1$ over $^0\wt\cW_{D,H}$, this is not the case. So, $\gamma$ is an isomorphism. This implies that the distinguished triangles (\ref{triangle _for_GH}) and (\ref{triangle_for_HG}) are isomorphic. Since $\tilde\nu_{\cW}$ is smooth, $\H_{GH}(\Aut_{G,H})$ is perverse, the intermediate extension under $^0\wt\cW_{D,H}\hook{} \wt\cW_{D,H}$. 
 
 Proposition~\ref{Pp_interm_extension_for_Th_2} and Theorem~\ref{Th_2} are proved.

\ssec{Local Rankin-Selberg-type convolutions}


\sssec{} Let $\cQ$ be the stack classifying $L_1\in \Bun_1$ and an exact sequence
\begin{equation}
\label{seq_cL_by_cLotimesOmega}
0\to L_1\otimes\Omega\to M_2\to L_1\to 0
\end{equation}
Let $\nu_{\cQ}: \cQ\to\Bun_{2,D}$ be the map sending the above point to $(M_2(-D)\subset \bar M_2\subset M_2)$, where $\bar M_2/M_2(-D)=L_1\otimes\Omega\mid_D$. Let $\tilde\nu_{\cQ}: \wt\cQ\to \Bunt_{2,D}$ be obtained from $\nu_{\cQ}$ by the base change $\Bunt_{2,D}\to\Bun_{2,D}$. Let $\ev_{\tilde Q}\to\A^1$ be the map sending the above point to the class of the sequence (\ref{seq_cL_by_cLotimesOmega}). A point of $\wt\cQ$ is given by (\ref{seq_cL_by_cLotimesOmega}) and a line $\cU$ equipped with 
$$
\eta: \cU^2\,\iso\, \det\RG(X, L_1\mid_D)\otimes\det\RG(X, \cO_D)
$$

 Let $\tilde q_{\cQ}: \wt\cQ\to\Bunt_1$ be the map sending the above point to $(L_1,\cU, \eta)$. Let $q_{\cQ}$ denote the composition $\wt\cQ\toup{\tilde q_{\cQ}}\Bunt_1\to\Bun_1$, where the second map is the projection forgetting $\cU$. Recall that $AE_1$ denotes the automorphic local system on $\Bun_1$ corresponding to $E_1$ (cf. Section~\ref{Sect_1.1.3_intro}). For a rank 1 local system $E_1$ on $X$ set
$$
\La_{E_1}=(\tilde\nu_{\cQ})_!(\tilde q_{\cQ}^*A\cE_0\otimes\ev_{\tilde\cQ}^*\cL_{\psi}\otimes q_{\cQ}^*AE_1)[\dim\cQ]
$$
 
Let $\tilde\cQ^{d_1}$ be the component of $\tilde\cQ$ given by $\deg(L_1\otimes\Omega)=d_1$. Write $\La_{E_1}^{d_1}$ for the contribution of $\tilde\cQ^{d_1}$ to $\La_{E_1}$. So, $\La_{E_1}^{d_1}$ is supported on the connected component of $\Bunt_{2,D}$ given by $\deg M_2=2d_1-2g+2$. 
 
\sssec{} Denote by $_0\Bunb^{P, r, d_1}_{\tilde D}$ the component of $_0\Bunb^P_{\tilde D}$ given by $\deg M_1=d_1$ and $\deg M=r+2d_1-\deg\Omega$. Let 
$\Eis(\Qlb\oplus \cE_0)^{r, d_1}$ denote the contribution of this component to $\Eis(\Qlb\oplus \cE_0)$. Recall the map $\delta_{\tilde D}: \Bunt_{2,D}\to\Bun_2$ defined in Section~\ref{Sect_2.3.2_Step2}. It sends $(\bar M\subset M, \cU)$ to $M$.

\begin{Lm} 
\label{Lm_for_RS_convolution_preparatory}
Let $E$ be an irreducible rank 2 local system on $X$. For $d_1\in\ZZ, r\ge 0$ one has canonically
$$
\Eis(\Qlb\oplus \cE_0)^{r, d_1}\otimes\delta^*_{\tilde D}\Aut_E\,\iso\, \Av^r_E(\La^{d_1}_{\det E})[\dim\Bun_2]
$$
So, $\Eis(\Qlb\oplus \cE_0)\otimes\delta^*_{\tilde D}\Aut_E\,\iso\, \oplus_{r\ge 0} \;\Av^r_E(\La_{\det E})[\dim\Bun_2]$. 
\end{Lm}
\begin{Prf}
Let $\cQ\Mod^r$ be the stack classifying a point of $\cQ$ given by (\ref{seq_cL_by_cLotimesOmega}) and an upper modification $M_2\subset M$ with $\deg(M/M_2)=r$, $M\in\Bun_2$. Let $\Bun'_2$ be the stack classifying $M\in\Bun_2$ together with a subsheaf $M_1\subset M$ such that $M_1\in \Bun_1$. 
Let $\varrho: \cQ\Mod^r\to \Bun'_2$ be the map sending the above point to $(M_1\subset M)$ with $M_1=L_1\otimes\Omega$. Let $\varrho': \Bun'_2\to\Bun_2$ be the map forgetting $M_1$. 

 Let $\gs: \cQ\Mod^r\to \Sh_0$ be the map sending the above point to $M/M_2$. Let $\ev_{\cQ\Mod}$ be the composition $\cQ\Mod^r\to \cQ\toup{\ev_{\cQ}}\A^1$, where the first map is the projection forgetting $M$. Let $q_{\cQ\Mod}:\cQ\Mod^r$ be the map sending the above point to $L_1$. By \cite[Section~7.9]{FGV}, there is an isomorphism
$$
\varrho'^*\Aut[\dimrel(\varrho')]\,\iso\, \ev_{\cQ\Mod}^*\cL_{\psi}\otimes \gs^*\cL^r_E\otimes q_{\cQ\Mod}^*(A\det E)[\dim\cQ\Mod^r]
$$
over the components of $\Bun'_2$ given by $\deg\Omega+\deg M-2\deg M_1\ge 0$.  

 For a point of $\cQ\Mod^r$ let $\bar M_2\subset M_2$ be the image of (\ref{seq_cL_by_cLotimesOmega}) under $\nu_{\cQ}$. Then the condition that the inclusion $M_1\hook{} M$ is maximal at $D$ is equivalent to requiring that $(\bar M_2\subset M_2\subset M)\in\Mod^r_{2,D}$, that is, $M\cap \bar M_2(D)=M_2$. 
The first isomorphism follows.

 To get the second one, use the cuspidality of $\Aut_E$. Namely, over the components of $\Bun'_2$ given by $\deg M-2\deg M_1+\deg\Omega<0$ one has $\varrho'^*\Aut_E=0$. The second isomorphism follows.
\end{Prf}

\sssec{} Write $\varepsilon: \Pic Y\to\Pic Y$ for the automorphism $\cB\mapsto \cB^*\otimes\Omega_Y$. For $d_1\in\ZZ, r\ge 0$ let 
$$
\mult^{r, d_1}: \Pic^{d_1}X \times Y^{(r)}\to \Pic^{2d_1+r} Y
$$
be the map sending $(M_1, D_1)$ to $(\phi^*M_1)(D_1)$.  

\begin{Pp} 
\label{Pp_RS_convolution}
Let $d_1\in\ZZ$, $r\ge 0$. For any rank 2 local system $E$ on $X$ and a rank one local system $E_1$ on $X$ one has
\begin{equation}
\label{complex_great_local}
F_H(\Av^r_E(\La^{d_1}_{E_1}))\,\iso\, (AE_1)_{\Omega^{-1}}\otimes\varepsilon_!\mult^{r, d_1}_!(AE_1\boxtimes (\phi^*E)^{(r)})[g-1+r]
\end{equation}
over $\Pic Y$.
\end{Pp}
\begin{Prf}
\Step 1 Case $r=0$. Consider the stack $\cQ\times_{\Bun_{2,D}}\Bun_{2,D,H}$, where we used the map $\nu_{\cQ}$ to define the fibre product. A point of this stack is given by an exact sequence (\ref{seq_cL_by_cLotimesOmega}) together with $\cB\in \Pic Y$ equipped with $N(\cB)\,\iso\, L_1^{-2}(D)$.

 We have a map $\cQ\times_{\Bun_{2,D}}\Bun_{2,D,H}\to \Bun^P_{D,H}$ sending the above point to $(L_1\otimes\Omega\subset M_2, \cB)$ equipped with $N(\cB)\,\iso\, L_1^{-2}(D)$. Apply Lemma~\ref{Lm_Aut_GH_over_Bun^P_DH} to describe the restriction of $\Aut_{G,H}$ to 
$\tilde\cQ\times_{\Bunt_{2,D}}\Bunt_{2,D,H}$. The complex $\pr_1^*\La_{E_1}\otimes\pr_2^*\Aut_{G,H}$ is the inverse image under the natural map
$$
\tilde\cQ\times_{\Bunt_{2,D}}\Bunt_{2,D,H}\to \cQ\times_{\Bun_{2,D}}\Bun_{2,D,H}
$$

 So, let $Z$ be the stack classifying pairs: a point of $\cQ\times_{\Bun_{2,D}}\Bun_{2,D,H}$ as above together with $s: \phi^* M_1\to \cB^*\otimes\Omega_Y$, where $M_1=L_1\otimes\Omega$. The norm of $s$ is a morphism $N(s): M_1^2\to N(\cB^*\otimes\Omega_Y)\,\iso\, L_1^2\otimes\Omega^2$. 

 Let $\ev_Z: Z\to \A^1$ be the map sending a point of $Z$ as above to the pairing of (\ref{seq_cL_by_cLotimesOmega}) with the section $N(s)\in \H^0(X, \cO)=k$. Over $Z$ we get the sheaf 
$$
\ev_Z^*\cL_{\psi}^{-1}\otimes \ev_{\cQ}^*\cL_{\psi}\otimes q_{\cQ}^* AE_1[\dim(\cQ\times_{\Bun_{2,D}}\Bun_{2,D,H})+1-g_Y],
$$ 
and $F_H(\La^{d_1}_E)$ identifies with its direct image with compact support under the map $Z\to \Pic Y$ sending the above point to $\cB$. 

 Consider the stack $\cW$ classifying $M_1=L_1\otimes\Omega\in \Pic^{d_1} X$, $\cB\in \Pic Y$ equipped with $N(\cB)\,\iso\, L_1^{-2}(D)$, and $s: \phi^*M_1\to \cB^*\otimes\Omega_Y$. We have a map $p_{Z}: Z\to \cW$ forgetting the extension (\ref{seq_cL_by_cLotimesOmega}). The complex 
$p_{Z !}(\ev^*\cL_{\psi}^{-1}\otimes \ev_{\cQ}^*\cL_{\psi}\otimes q_{\cQ}^* AE_1)$ is the extension by zero from the locus given by $N(s)=1$. Over the latter locus $s$ is an isomorphism, so this locus identifies with $\Pic ^{d_1} X$. The resulting map $\Pic^{d_1} X\to \Pic Y$ sends $M_1$ to $\varepsilon(\phi^* M_1)$. The case $r=0$ is done.

\medskip
\Step 2 By Theorem~\ref{Th_1} and Step 1, we get
$$
F_H(\Av^r_E(\La^{d_1}_{E_1}))\,\iso\, \Av^r_E((AE_1)_{\Omega^{-1}}\otimes \epsilon_!\mult_!^{0,d_1}(AE_1))[g-1]
$$
It is straightforward to identify the latter complex with the right hand side of (\ref{complex_great_local}). 
\end{Prf}

\begin{Rem} Let $r\ge 1$. For $0\le a\le r/2$ let $s^a: (X-D)^{(a)}\times Y^{(r-2a)}\to Y^{(r)}$ be the map sending $(\bar D_1, D_2)$ to $\phi^*\bar D_1+D_2$. 
For any rank 2 local system $E$ on $X$ the constructible sheaf $(\phi^*E)^{(r)}$ on $Y^{(r)}$ admits a filtration $0=F_{-1}\subset F_0\subset\ldots\subset F_1\subset...$
with successive quotients being
\begin{equation}
\label{succ_quotients_(phi^*E)^r}
F_a/F_{a-1}=s^a_!((\det E)^{(a)}\boxtimes \phi^*(E^{(r-2a)}))
\end{equation}
for $0\le a\le r/2$. This is an analog of \cite[Remark~7]{Ly1} in our ramified setting. Let $\cY^r$ be the stack classifying $L_1\in \Bun_1$, $M\in \Bun_2$, $\bar D\in X^{(r)}$, an inclusion $s_1: M_1=L_1\otimes\Omega\hook{} M$ which is maximal in a neighbourhood of $D$, and an isomorphism 
$$
s_2: L_1^2\otimes\Omega(\bar D)\,\iso\, \det M
$$ 
We have a natural map $\cQ\Mod^r\to \cY^r$ sending $(M_1\subset M_2\subset M)$ to $(M_1\subset M, s_2, \bar D)$, where $\bar D=\div(M/M_2)$ and $s_2$ is the induced trivialization. We have a map $\phi_{\cY}: \cY^r\to \Bun_{2,D}$ sending the above point to $(\bar M\subset M)$, where
$\bar M/M(-D)=M_1\mid_D$. Let $\tilde\phi_{\cY}: \wt\cY^r\to\Bunt_{2,D}$ be obtained from $\phi_{\cY}$ by the base change $\Bunt_{2,D}\to\Bun_{2,D}$. The complex $\Av^r_E(\La^{d_1}_{E_1})$ can be seen as a direct image with compact support under $\tilde\phi_{\cY}$ of the Whittaker sheaf, cf. \cite[Definition~2]{Ly3}. The stack $\tilde\cY^r$ has a stratification given by fixing the degree of the subbundle in $M$ generated by $M_1$. One checks that calculating $F_H(\Av^r_E(\La^{d_1}_{E_1}))$ with respect to this stratification, one gets a stratification of $(\phi^*E)^{(r)}$ by constructible subsheaves
with successive quotients (\ref{succ_quotients_(phi^*E)^r}). We believe but we have not checked that the two above stratifications on $(\phi^*E)^{(r)}$ coincide.
\end{Rem}

\section{Restriction of the theta-sheaf and ramifications}
\label{Sect_Restriction_ramifications}

This section is essentially independent of the others. Our purpose here is to prove Proposition~\ref{Pp_local_nature_gamma_forgetting} below, to be used in the proof of our main results in Section~\ref{Sect_Waldspurger_periods_main}. 

 Our Proposition~\ref{Pp_local_nature_gamma_forgetting} is some version of the Hecke property of the theta kernel $\Aut_{G,\tilde H}$ for the dual pair $(G=\GL_2, \tilde H=\GO_{2m}^0)$ adopted to our particular ramifications. It says that some averaging of $\Aut$ along $G$ is isomorphic to another averaging of $\Aut$ along $\tilde H$. 

\ssec{} Let $D$ be a multiplicity free divisor on $X$ of degree $d\ge 0$ (in this section $d$ need not be even). Fix $m\ge 1$. Let $G=\GL_2$. Let $\Bun_Q$ be the stack classifying $V\in\Bun_{2m},\cC\in\Bun_1$ with a lower modification $V(-D)\subset V^{\perp}\subset V$ such that $\div(V/V^{\perp})=2D$, a morphism $b: \Sym^2 V\to \cC(D)$ inducing an isomorphism $V\,\iso\, (V^{\perp})^*\otimes\cC$, and a compatible trivialization $\gamma: (\det V)(-D)\,\iso\, \cC^m$. This means that $\gamma^2$ is the isomorphism induced by $b$. 
 
  Let $\tilde H=\GO_{2m}^0$, the connected component of unity of the split orthogonal similitude group. So, $\Bun_{\tilde H}$ classifies $V\in \Bun_{2m}, \cC\in\Bun_1$ with a symmectric form $\Sym^2 V\to\cC$ inducing an isomorphism $V\,\iso\, V^*\otimes\cC$, and a compatible trivialization $\det V\,\iso\, \cC^m$. Let $\Bun_{Q,D}$ be the stack classifying a point of $\Bun_Q$ as above together with a lower modification $V^{\perp}\subset \bar V\subset V$ such that $(\bar V,b,\gamma)\in\Bun_{\tilde H}$. This implies that $\div(V/\bar V)=D$. Write $q_Q: \Bun_{Q,D}\to\Bun_Q$ for the map $q_Q$ forgetting $\bar V$. 

 Let $\Bun_{G, Q}$ be the stack classifying a point of $\Bun_Q$ as above with $M\in\Bun_2$, for which we set $\cA=\det M$, and an isomorphism $\cA\otimes\cC\,\iso\, \Omega$. Set 
$$
\Bun_{2,D,Q}=\Bun_{2,D}\times_{\Bun_G}\Bun_{G,Q},
$$  
where the map $\Bun_{2,D}\to\Bun_G$ sends $(\bar M\subset M)$ to $M$. Let 
$$
\tau_Q: \Bun_{2,D,Q}\to\Bun_{G_{2m}}
$$ 
be the map sending $(\bar M\subset M)\in\Bun_{2,D}, (V^{\perp}\subset V)\in\Bun_Q$ to $M'$ defined by the cartesian square
\begin{equation}
\label{diag_def_of_M'_forSection_311}
\begin{array}{ccc}
M\otimes V & \to& M\otimes(V/V^{\perp})\\
\uparrow && \uparrow\\
M' & \to & \bar M/M(-D)\otimes (V/V^{\perp})
\end{array}
\end{equation}
The form $\wedge^2(M\otimes V)\to \Omega(D)$ induces a symplectic form on $M'$. 

\begin{Lm} 
\label{Lm_for_modification_of_theta_detRG}
1) For a point of $\Bun_{G,Q}$ as above one has canonically
$$
\det\RG(X, M\otimes V)\,\iso\, \frac{\det\RG(X, M)^{2m}\otimes \det\RG(X, V)^2}{\det\RG(X, \cA)^{2m}\otimes\det\RG(X, \cO)^{2m}\otimes\det\RG(X, V/V^{\perp})}
$$
2) For a point of $\Bun_{2,D,Q}$ as above one has canonically
$$
\det\RG(X, M')\,\iso\, \frac{\det\RG(X, M)^{2m}\otimes \det\RG(X, V^{\perp})^2}{\det\RG(X, \cA)^{2m}\otimes\det\RG(X, \cO)^{2m}\otimes\det\RG(X, M/\bar M)^2}
$$
\end{Lm}
\begin{Prf}
1) is proved as in Lemma~\ref{Lm_detRG_M_otimes_L}. \\
2) There is a canonical $\ZZ/2\ZZ$-graded isomorphism 
$$
\det\RG(X, (M/\bar M)\otimes (V/V^{\perp}))\,\iso\, \det\RG(X, M/\bar M)^2\otimes \det\RG(X, V/V^{\perp})
$$
Our claim follows. 
\end{Prf} 

\sssec{} The map $\tau_Q$ extends to a morphism $\tilde\tau_Q: \Bun_{2,D,Q}\to\Bunt_{G_{2m}}$ sending the above point to $(M', \cB)$, where
$$
\cB=\frac{\det\RG(X, M)^{m}\otimes \det\RG(X, V^{\perp})}{\det\RG(X, \cA)^{m}\otimes\det\RG(X, \cO)^{m}\otimes\det\RG(X, M/\bar M)}
$$
is equipped with the isomorphism $\cB^2\,\iso\, \det\RG(X, M')$ given by Lemma~\ref{Lm_for_modification_of_theta_detRG}.

\sssec{} 
\label{Sect_4.1.3_G_tildeH}
As in \cite[Section~3.2.1]{Ly1}, we denote by $\Bun_{G,\tilde H}$ the stack classifying $M\in\Bun_2$ with $\cA=\det M$, $(V,\cC)\in \Bun_{\tilde H}$, and an isomorphism $\Omega\,\iso\, \cA\otimes\cC$ on $X$. We will use the morphism $\tilde\tau: \Bun_{G,\tilde H}\to\Bunt_{G_{2m}}$ introduced in \cite[Section~3.2.1]{Ly1}. It sends the above point to $(M\otimes V, \cB)$, where 
$$
\cB=\frac{\det\RG(X, V)\otimes \det\RG(X, M)^m}{\det\RG(X, \cO)^m\otimes\det\RG(X, \cA)^m}
$$
is equipped with a canonical isomorphism $\cB^2\,\iso\, \det\RG(X, M\otimes V)$.
 
\sssec{} Let $\Bun_{G,Q,D}$ be the stack obtained from $\Bun_{G,Q}$ by the base change $q_Q: \Bun_{Q,D}\to\Bun_Q$. Let $p_Q: \Bun_{G,Q,D}\to\Bun_{G,\tilde H}$ be the map sending the above point to $(M, \bar V)$. Consider the diagram
$$
\Bunt_{G_{2m}}\getsup{\tilde\tau_Q}
\Bun_{2,D,Q}\toup{\delta_D} \Bun_{G,Q}\getsup{q_Q} \Bun_{G,Q,D}\toup{p_Q} \Bun_{G,\tilde H}\toup{\tilde\tau}\Bunt_{G_{2m}},
$$
where $\delta_D: \Bun_{2,D}\to\Bun_G$ is the projection forgetting $\bar M$. 

\begin{Pp} 
\label{Pp_local_nature_gamma_forgetting}
There is an isomorphism
$$
q_{Q !}p_Q^*\tilde\tau^*\Aut[\dimrel(p_Q\comp\tilde\tau)]\,\iso\, \delta_{D !}\tilde\tau_Q^*\Aut[\dimrel(\tilde\tau_Q)]
$$
\end{Pp}

\ssec{Proof of Proposition~\ref{Pp_local_nature_gamma_forgetting}}

\sssec{} Let $X_D$ be the formal neighbourhood of $D$ in $X$, $F=k(X)$ and $F_D=\prod_{x\in D} F_x$, where $F_x$ is the completion of $F$ at $x$. Write $\cO_{X_D}$ for the structure sheaf of $X_D$. 

 We will use the following version of the theta-sheaf from \cite{LL}. Let $W$ be a rank $2n$ vector bundle on $X_D$ with a symplectic form $\wedge^2 W\to\Omega$. Write $\cL_d(W(F_D))$ for the scheme of discrete lagrangians in the Tate space $W(F_D)$. We use the notation for relative determinants from \cite[Section~5.2]{LL}. Let $\wt\cL_d(W(F_D))$ be the stack classifying $R\in \cL_d(W(F_D))$, a $\ZZ/2\ZZ$-graded line $\cB$ of parity zero, and a $\ZZ/2\ZZ$-graded isomorphism $\cB^2\,\iso\,\det(W: R)$. One has the nonramified Weil category $W(\wt\cL_d(W(F_D)))$ defined as in \cite[Theorem~2]{LL}, and the perverse theta-sheaf $S_{W(F_D)}$ defined as in \cite[Section~6.5]{LL}. It is naturally $\Sp(W(\cO_{X_D}))$-equivariant. 

 Write $\cX^D_ n$ for the stack quotient $\wt\cL_d(W(F_D))/\Sp(W(\cO_{X_D}))$, it classifies a rank $2n$ vector bundle $\cW$ on $X_D$ with a symplectic form $\wedge^2\cW\to\Omega$, $R\in \cL_d(\cW(F_D))$, a line $\cB$ equipped with $\cB^2\,\iso\,\det(W: R)$ as above. The perverse theta-sheaf on $\cX^D_ n$ will be denoted $S_{\cX}$ for brevity (notation compatible with \cite[Section~6.5]{LL}). This is an object of the category $\D_{\Sp(W(\cO_{X_D}))}(\wt\cL_d(W(F_D)))$ defined as in \cite[Section~5.3]{LL}. 
 
\sssec{} Fix a rank 2 vector bundle $M_0$ on $X_D$ with $\cA_0=\det M_0$. Fix also a line bundle $\cC_0$ and a vector bundle $V_0$ of rank $2m$ on $X_D$ with a lower modification $V_0(-D)\subset V_0^{\perp}\subset V_0$ such that $\div(V_0/V_0^{\perp})=2D$, a map $b: \Sym^2 V_0\to \cC_0(D)$ inducing $V_0\,\iso\, (V_0^{\perp})^*\otimes\cC$, and a compatible trivialization $\gamma: (\det V_0)(-D)\,\iso\, \cC_0^m$. 
We assume given an isomorphism $\cA_0\otimes\cC_0\,\iso\, \Omega$ over $X_D$. 

 Let $^D\Bun_{G,Q}$ be the stack classifying a point of $\Bun_{G,Q}$ as above together with compatible trivializations 
\begin{equation}
\label{trivialization_upperscript_D} 
 M\,\iso\, M_0\mid_{X_D}, \;\;\cA\,\iso\, \cA_0\mid_{X_D}, \;\;V\,\iso\, V_0\mid_{X_D}, \;\;\cC\,\iso\, \cC_0\mid_{X_D}
\end{equation} 
preserving all the structures, in particular inducing $V^{\perp}\,\iso\, V_0^{\perp}$. 
 
 We have the group scheme over $X_D$
$$
\cG=\{(g, h)\in \GSp(M_0)\times\GO^0(V_0)\mid (g,h)\;\;\mbox{acts trivially on} \;\;\cA_0\otimes\cC_0, \; h(V_0^{\perp})=V_0^{\perp}\}
$$ 
Here $\GO^0(V_0)$ denotes the connected component of unity of the orthogonal similitude group $\GO(V_0)$. The group ind-scheme $\cG(F_D)$ acts naturally on $^D\Bun_{G,Q}$. This is analogous to \cite[Section~2.3.3]{BD}.

 Let $\cL_d(M_0\otimes V_0(F_D))$ denote the scheme of discrete lagrangians in $M_0\otimes V_0$. Let 
$$
\xi: {^D\Bun_{G,Q}}\to \cL_d(M_0\otimes V_0(F_D))
$$ 
be the map sending $(M,V)$ to $\H^0(X-D, M\otimes V)$ which becomes a discrete lagrangian in  $M_0\otimes V_0(F_D)$ via (\ref{trivialization_upperscript_D}). The map $\xi$ is $\cG(F_D)$-equivariant. Taking the stack quotient by $\cG(\cO_{X_D})$, one gets a map still denoted 
\begin{equation}
\label{map_xi_for_Section_313}
\xi: \Bun_{G,Q}\to \cL_d(M_0\otimes V_0(F_D))/\cG(\cO_{X_D})
\end{equation}
by abuse of notation. 

\sssec{} Write $\cY_{G,Q}$ for the stack classifying $R\in \cL_d(M_0\otimes V_0(F_D))$ and a $\ZZ/2\ZZ$-graded line $\cB$ of parity zero equipped with a $\ZZ/2\ZZ$-graded isomorphism
\begin{equation}
\label{iso_cB_square_for_cY_GQ}
\cB^2\,\iso\, \det(M_0\otimes V_0: R)\otimes \det(V_0: V_0^{\perp})
\end{equation} 

 Let $\cY_{2,D,Q}$ denote the stack classifying a point $(\cB, R)\in \cY_{G,Q}$ and a lower modification $M_0(-D)\subset \bar M\subset M_0$ with $\div(M_0/\bar M)=D$. 
 
 Let $\cY_{G,Q,D}$ denote the stack classifying a point $(\cB, R)\in \cY_{G,Q}$ and a lower modification $V_0^{\perp}\subset \bar V\subset V_0$ such that $(\bar V, b,\gamma)$ is a $\tilde H$-torsor on $X_D$. 
Write 
\begin{equation}
\label{diag_over_cY_GQ}
\cY_{2,D,Q}\;\toup{\bar \delta_D}\;\cY_{G,Q}\;\getsup{\bar q_Q}\;\cY_{G,Q,D}
\end{equation} 
for the projections forgetting $\bar M$ and $\bar V$ respectively. 

\sssec{} Let $\bar\tau_Q: \cY_{2,D,Q}\to \cX^D_{2m}$ be the map sending 
$(\cB, R,\bar M)\in\cY_{2,D,Q}$ to $(M', R, \cB_1)$, where $M'$ is given by the cartesian square
$$
\begin{array}{ccc}
M_0\otimes V_0 & \to & M_0\otimes(V_0/V_0^{\perp})\\
\uparrow && \uparrow\\
M' & \to & \bar M/M_0(-D)\otimes(V_0/V_0^{\perp}),
\end{array}
$$
and 
$$
\cB_1=\frac{\cB}{\det(V_0: V_0^{\perp})\otimes\det(M_0/\bar M)}
$$
is equipped with the induced isomorphism $\cB_1^2\,\iso\, \det(M': R)$. We have used a canonical isomorphism 
$$
\det((M_0/\bar M)\otimes (V_0/V_0^{\perp}))\,\iso\, \det(M_0/\bar M)^2\otimes \det(V_0: V_0^{\perp})
$$

\sssec{} Let $\bar p: \cY_{G,Q,D}\to \cX^D_{2m}$ be the map sending $(R,\cB, \bar V)$ to $(M_0\otimes\bar V, R, \cB_2)$, where 
$$
\cB_2=\cB\otimes\det(V_0: \bar V)^{-1}
$$
is equipped with the induced isomorphism $\cB_2^2\,\iso\, \det(M_0\otimes\bar V: R)$. We have used the fact that a choice of a maximal isotropic submodule $\bar V/V_0^{\perp}\subset V_0/V_0^{\perp}$  yields a $\ZZ/2\ZZ$-graded isomorphism 
$$
\det(V_0: V_0^{\perp})\,\iso\, \det(\cC_0(D): \cC_0)\otimes \det(\cO: \cO(-D))
$$

 The group $\cG(\cO_{X_D})$ acts on the diagram (\ref{diag_over_cY_GQ}) naturally, it is understood that it acts trivially on $\cB$. By abuse of notation, we denote by
$$
\bar\cY_{2,D,Q}\;\toup{\bar \delta_D}\;\bar\cY_{G,Q}\;\getsup{\bar q_Q}\;\bar\cY_{G,Q,D}
$$
the diagram obtained from (\ref{diag_over_cY_GQ}) by taking the stack quotients by $\cG(\cO_{X_D})$.

\sssec{} Extend (\ref{map_xi_for_Section_313}) to a morphism $\tilde\xi: \Bun_{G,Q}\to \bar\cY_{G,Q}$ sending $(M,V)$ to $(R,\cB)$, where $R=\H^0(X-D, M\otimes V)$ and
$$
\cB=\frac{\det\RG(X, M)^m\otimes\det\RG(X, V)}{\det\RG(X,\cA)^m\otimes\det\RG(X,\cO)^m}
$$
is equipped with the isomorphism (\ref{iso_cB_square_for_cY_GQ}) given by 1) of Lemma~\ref{Lm_for_modification_of_theta_detRG}. 
One gets the diagram
$$
\begin{array}{ccccc}
\Bun_{2,D,Q}&\toup{\delta_D} &\Bun_{G,Q}&\getsup{q_Q}&\Bun_{G,Q,D}\\
\downarrow && \downarrow\lefteqn{\scriptstyle \tilde\xi} &&\downarrow\\
\bar\cY_{2,D,Q} &\toup{\bar \delta_D} & \bar\cY_{G,Q} &\getsup{\bar q_Q} & \bar\cY_{G,Q,D},
\end{array}
$$
where both squares are cartesian, the left and right vertical arrows are defined in a way similar to $\tilde\xi$. The above diagram extends to the following one
\begin{equation}
\label{diag_big_one_for_Section313}
\begin{array}{ccccccccc}
\Bunt_{G_{2m}} & \getsup{\tilde\tau_Q}&
\Bun_{2,D,Q}&\toup{\delta_D} &\Bun_{G,Q}&\getsup{q_Q} &\Bun_{G,Q,D}&\toup{\tilde\tau\comp p_Q}&\Bunt_{G_{2m}} \\ 
\downarrow\lefteqn{\scriptstyle\tilde\xi_{\cX}} && \downarrow &&  \downarrow\lefteqn{\scriptstyle\tilde\xi} &&  \downarrow &&  \downarrow\lefteqn{\scriptstyle\tilde\xi_{\cX}} \\
\cX^D_{2m} & \getsup{\bar\tau_Q} & \bar\cY_{2,D,Q} &\toup{\bar \delta_D} & \bar\cY_{G,Q} &\getsup{\bar q_Q} & \bar\cY_{G,Q,D} & \toup{\bar p} & \cX^D_{2m}
\end{array}  
\end{equation}
Here $\tilde\xi_{\cX}$ sends $(\cM,\cU)\in\Bunt_{G_{2m}}$, where $\cU^2\,\iso\, \det\RG(X, \cM)$, to the discrete lagrangian $L=\H^0(X-D, \cM)\subset \cM(F_D)$ together with $\cU^2\,\iso\, \det(\cM(\cO_{X_D}): L)$. 

\sssec{}  One defines the derived category $\D(\bar\cY_{G,Q}):=\D_{\cG(\cO_{X_D})}(\cY_{G,Q})$ as in \cite[Section~5.3]{LL}. One also defines a functor 
\begin{equation}  
\label{functor_tildexi^*}
\tilde\xi^*: \D(\bar\cY_{G,Q})\to \D(\Bun_{G,Q})
\end{equation}
as in \cite[Section~7.2]{LL}, it is normalized to be "as exact for the perverse t-structures as possible". Here is the precise definition. 
  
  For $r\ge 0$ let $_r\Bun_{G,Q}\subset\Bun_{G,Q}$ be the open substack given by $\H^0(X, M\otimes V^{\perp}(-rD))=0$ for a point $(M, V^{\perp}\subset V, \cC)\in \Bun_{G,Q}$ as above. The category $\D(\Bun_{G,Q})$ identifies with the projective 2-limit of $\D(_r\Bun_{G,Q})$ as $r$ goes to infinity. 

 For $N\ge r\ge 0$ and $(M, V^{\perp}\subset V, \cC,b,\gamma)\in{_r\Bun_{G,Q}}$ the finite-dimensional $k$-vector space 
\begin{equation}
\label{symplectic space Sect_4.2.7}
_N(M\otimes V):=(M\otimes V)(ND)/(M\otimes V^{\perp}(-ND))
\end{equation} 
is symplectic, and $\H^0(X, (M\otimes V)(ND))\subset {_N(M\otimes V)}$ is lagrangian. The natural map $\H^0(X, (M\otimes V)(ND))\to {_N(M\otimes V)}$ is injective, because $\H^0(X, M\otimes V^{\perp}(-ND))=0$ for a point of $_r\Bun_{G,Q}$. Define the vector space $_N(M_0\otimes V_0)$ as in (\ref{symplectic space Sect_4.2.7}). 

Let $_r\cL(_N(M\otimes V))$ be the scheme of lagrangians $L$ in $_N(M\otimes V)$ such that 
$$
L\cap (M\otimes V^{\perp}(-rD)/(M\otimes V^{\perp}(-ND)))=0
$$ 
This is a smooth scheme of finite type. The group $\cG(\cO/\cO(-(2N+1)D))$ acts naturally on $_r\cL(_N(M_0\otimes V_0))$. Define the scheme $_r\cL(_N(M_0\otimes V_0))$ similarly. 

 For $N\ge r\ge 0$ and $r_1\ge 2N+1$ we get a morphism
$$
\xi_N: {_r\Bun_{G,Q}}\to {_r\cL(_N(M_0\otimes V_0))/\cG(\cO/\cO(-r_1D))}
$$
sending the above point to $\H^0(X, (M\otimes V)(ND))\subset {_N(M\otimes V)}$. 

 We view $M_0\otimes V_0, (M_0\otimes V_0^{\perp})(-ND)$ as compact lattices in the Tate space $M_0\otimes V_0(F_D)$. The following relative determinants in the notations from \cite[Section~5.2]{LL} become
$$
\det(M_0\otimes V_0: (M_0\otimes V_0^{\perp})(-ND))\,\iso\, \det(M_0\otimes V_0/(M_0\otimes V_0^{\perp})(-ND))
$$
and $\det(V_0: V_0^{\perp})\,\iso\, \det(V_0/V_0^{\perp})$. 

Over $_r\cL(_N(M_0\otimes V_0))$ one has the line bundle with fibre
$$
\det(M_0\otimes V_0: (M_0\otimes V_0^{\perp})(-ND))\otimes\det L\otimes \det(V_0: V_0^{\perp})
$$ 
at $L\in {_r\cL(_N(M_0\otimes V_0))}$, it is naturally $\cG(\cO/\cO(-r_1D))$-equivariant. Denote by 
$$
_r\tilde\cL(_N(M_0\otimes V_0))/\cG(\cO/\cO(-r_1D))
$$ 
the gerbe of square roots of this line bundle over the stack quotient $_r\cL(_N(M_0\otimes V_0))/\cG(\cO/\cO(-r_1D))$. The map $\xi_N$ extends to
$$
\tilde\xi_N: {_r\Bun_{G,Q}}\to {_r\tilde\cL(_N(M_0\otimes V_0))/\cG(\cO/\cO(-r_1D))}
$$
Further, for $N_1\ge N\ge r\ge 0$ and $r_1\ge 2N_1+1$ one has a commutative diagram
$$
\begin{array}{ccc}
_r\Bun_{G,Q} & \toup{\tilde\xi_N} & _r\tilde\cL(_N(M_0\otimes V_0))/\cG(\cO/\cO(-r_1D))\\
& \searrow\lefteqn{\scriptstyle \tilde\xi_{N_1}} & \uparrow\lefteqn{\scriptstyle p_{\cL}}\\
&& _r\tilde\cL(_{N_1}(M_0\otimes V_0))/\cG(\cO/\cO(-r_1D))
\end{array}
$$
The functors $K\mapsto \tilde\xi_N^*K[\dimrel(\tilde\xi_N)]$ from $D_{\cG(\cO)}( _r\tilde\cL(_N(M_0\otimes V_0)))$ to $\D(_r\Bun_{G,Q})$ are compatible with the transition functors $p_{\cL}^*[\dimrel(p_{\cL})]$, so passing to the limit by $N$ they yield a functor
$$
_r\tilde\xi^*: D_{\cG(\cO)}( _r\tilde\cL(M_0\otimes V_0(F_D)))\to \D(_r\Bun_{G,Q})
$$
Passing to the limit by $r$, one defines the desired functor (\ref{functor_tildexi^*}). 

 The other restriction functors for the diagram (\ref{diag_big_one_for_Section313}) are defined similarly. Thus, Proposition~\ref{Pp_local_nature_gamma_forgetting} is reduced to the following local claim.
 
\begin{Pp} 
\label{Pp_local_nature_gamma_forgetting_second}
There is a $\cG(\cO_{X_D})$-equivariant isomorphism over $\cY_{G,Q}$
\begin{equation}
\label{iso_for_Pp_local_nature_gamma_forgetting_second}
\bar \delta_{D !}\bar\tau_Q^*\cS_{\cX}[\dimrel(\bar\tau_Q)]\,\iso\, 
\bar q_{Q !} \bar p^*\cS_{\cX}[\dimrel(\bar p)]
\end{equation}
\end{Pp}
\begin{Prf}
Assume in addition given lower modification $V_0^{\perp}\subset \bar V_0\subset V_0$ such that $(\bar V_0, b, \gamma)$ is a $\tilde H$-torsor over $X_D$. This yields a decomposition $V_0/V_0^{\perp}\,\iso\, U\oplus U^*\otimes \cC_0(D)\mid_D$, where $U$ is a line bundle on $D$ such that $\bar V_0/V_0^{\perp}=U$, and $U^*\otimes \cC_0(D)\mid_D$ is isotropic in $V_0/V_0^{\perp}$. 

 Let $\wt\cL_d(M_0\otimes V_0(F_D))$ denote the stack classifying $R\in \cL_d(M_0\otimes V_0(F_D))$, a $\ZZ/2\ZZ$-graded line $\bar\cB$ of parity zero equipped with a $\ZZ/2\ZZ$-graded isomorphism
$
\bar\cB^2\,\iso\, \det(M_0\otimes \bar V_0: R)
$.
We get an isomorphism
$$
\cY_{G,Q}\,\iso\, \wt\cL_d(M_0\otimes V_0(F_D))
$$ 
sending $(R,\cB)$ equipped with (\ref{iso_cB_square_for_cY_GQ}) to $(R,\bar\cB)$, where $\bar\cB=\cB\otimes\det(V_0: \bar V_0)^{-1}$ with the induced isomorphism 
$$
\bar\cB^2\,\iso\, \det(M_0\otimes \bar V_0: R)
$$ 
Let also $\wt\Sp(M_0\otimes V_0(F_D))$ be the metaplectic group stack classifying $g\in \Sp(M_0\otimes V_0(F_D))$ and a line $\cB$ together with $\cB^2\,\iso\, \det(M_0\otimes\bar V_0: g(M_0\otimes\bar V_0))$.

 Let $W=M_0\otimes V_0^{\perp}$, this is a compact lattice in $M_0\otimes V_0(F_D)$ and $W^{\perp}=M_0\otimes V_0$. Pick 
$$
L\in \cL_d(M_0\otimes V_0(F_D))
$$ 
such that $L\cap W=0$ and the image of $L\cap W^{\perp}$ in $W^{\perp}/W$ is $L_W:=M_0/M_0(-D)\otimes U$. 

 We use some notions and notations from \cite{LL}. Let $H_W=(W^{\perp}/W)\oplus k$ be the Heisenberg group of $W^{\perp}/W$. Let $L^0\in \wt\cL(M_0\otimes V_0(F_D))$ be the lifting of $L$ defined by the enhanced structure $(L, \cB_L)$, where $\cB_L=\det L_W$ is equipped with a canonical isomorphism 
$$
\cB_L^2\,\iso\, \det L_W\otimes\det(M_0\otimes \bar V_0: W)\,\iso\, \det(M_0\otimes\bar V_0: L)
$$  
As in \cite[Sections~6.1-6.3]{LL}, one has the equivariant derived category $\D\cH_{L_W}$ and the functor 
$$
\cF_{L^0, W}: \D\cH_{L_W}\to \D(\wt\cL_d(M_0\otimes V_0(F_D)))\,\iso\, \D(\cY_{G,Q})
$$
The lagrangian subspace $(M_0\otimes \bar V_0)/W \subset W^{\perp}/W$ yields an equivalence 
$$
\D\cH_{L_W}\,\iso\, \D(M_0\otimes (V_0/\bar V_0))
$$ 
exact for the perverse t-structures . 

 Define the following object $K_{\bar V_0}\in \D\cH_{L_W}$. Let $Z$ be the scheme classifying the pairs: a lower modification $M_0(-D)\subset \bar M_0\subset M_0$ with $\div(M_0/\bar M_0)=D$, and $y\in (\bar M_0/M_0(-D))\otimes (V_0/\bar V_0)$. Let $p_Z: Z\to M_0\otimes (V_0/\bar V_0)$ be the map sending $(\bar M_0, y)$ to $y$. Set
$$
K_{\bar V_0}=p_{Z !}\Qlb[\dim Z]
$$
Denote by $J$ the set of all lower modifications $V_0^{\perp}\subset \bar V\subset V_0$ such that $(\bar V, b,\gamma)$ is a $\tilde H$-torsor over $X_D$. The set $J$ is naturally in bijection with the set of effective divisors $D'$ satisfying $D'\le D$. 

 For $\bar V\in J$ let $I_{\bar V}$ denote the constant perverse sheaf on $M_0\otimes((\bar V+\bar V_0)/\bar V_0)$ extended by zero to $M_0\otimes(V_0/\bar V_0)$. 
One easily gets a canonical isomorphism
\begin{equation}
\label{iso_in_DH_LW}
K_{\bar V_0}\,\iso\, \oplus_{\bar V\in J} I_{\bar V}
\end{equation}
Applying $\cF_{L^0, W}$ to (\ref{iso_in_DH_LW}) one gets the isomorphism (\ref{iso_for_Pp_local_nature_gamma_forgetting_second}). In particular, $\cF_{L^0, W}(I_{\bar V})$ is isomorphic to the theta-sheaf on $\cY_{G,Q}$ corresponding to the lagrangian compact lattice $M_0\otimes \bar V$ in $M_0\otimes V_0(F_D)$. 

 Now let $\bar V\in J$. Pick $N\in \cL_d(M_0\otimes V_0(F_D))$ such that $N\cap W=0$, and the image of $N\cap W^{\perp}$ in $W^{\perp}/W$ is $M_0\otimes (\bar V/V_0^{\perp})$. One similarly defines $N^0\in \wt\cL_d(M_0\otimes V_0(F_D))$, the functor $\cF_{N^0, W}: \D\cH_{N_W}\to \D(\cY_{G,Q})$, and the object $K_{\bar V}\in \D\cH_{N_W}$. 
 
  As in \cite[Section~6.2]{LL}, one gets a canonical interwining functor
$\cF_{N^0_W, L^0_W}: \D\cH_{L_W}\to \D\cH_{N_W}$. By \cite[Theorem~2]{LL}, the diagram 
$$
\begin{array}{cccc}
\D\cH_{L_W} && \toup{\cF_{L^0, W}} & \D(\cY_{G,Q})\\
\downarrow\lefteqn{\scriptstyle \cF_{N^0_W, L^0_W}} && 
\nearrow\lefteqn{\scriptstyle \cF_{N^0, W}}\\
\D\cH_{N_W}
\end{array}
$$
is canonically 2-commutative. Besides, $\cF_{N^0_W, L^0_W}(K_{\bar V_0})\,\iso\, K_{\bar V}$ naturally, and $\cF_{N^0_W, L^0_W}$ transforms (\ref{iso_in_DH_LW}) to a similar isomorphism for $\bar V$. For this reason the isomorphism (\ref{iso_for_Pp_local_nature_gamma_forgetting_second}) obtained is independent of our choice of $\bar V_0, L$, so it is  $\cG(\cO_{X_D})$-equivariant. 
\end{Prf}

\medskip 

Proposition~\ref{Pp_local_nature_gamma_forgetting} is also proved.

\section{Waldspurger periods}
\label{Sect_Waldspurger_periods_main}

In this section we finish the proof of Theorem~\ref{Th_3}, one of our main results. We follow the strategy explained in Section~\ref{Sect_Strategy}.

\ssec{The set-up} 

\sssec{}
\label{Sect_5.1.1_set-up}
Keep notations of Section~\ref{Sect_Geometric_Waldspurger_intro}. Recall the group scheme $U_{\phi}$ on $X$ given by the cokernel of the natural map $\Gm\to \phi_*\Gm$. As in \cite{Ly1}, we denote by $e_{\phi}: \Pic Y\to \Bun_{U_{\phi}}$ the extension of scalars map. Recall that $\Bun_{U_{\phi}}$ identifies with the stack classifying $\cB\in\Pic Y$ equipped with $N(\cB)\,\iso\, \cO_X$ and a compatible isomorphism $\gamma: \cB\mid_{D_Y}\,\iso\, \cO_{D_Y}$. However, this identification is not canonical as there is an automorphism $\cB\mapsto \cB^{-1}$ of $\Bun_{U_{\phi}}$. Our convention is that for the latter identification $e_{\phi}$ sends $\cB$ to $\sigma^*\cB\otimes\cB^{-1}$, this agrees with the map $e_{\phi}$ of \cite[Section~6.3.3]{Ly1}. 

 For $a\in\ZZ/2\ZZ$ we write $\Bun_{U_{\phi}}^a$ for the
connected component of $\Bun_{U_{\phi}}$ corresponding to $a$, so that $\Bun^0_{U_{\phi}}$ is the connected component of unity (cf. \select{loc.cit.})

\sssec{} Let $G=\GL_2$ and $\tilde H=\GO^0_4$, it is included into an exact sequence 
\begin{equation}
\label{seq_def_GO^0_4}
1\to\Gm\to \GL_2\times\GL_2\to \tilde H\to 1,
\end{equation}
where the fist map is $z\mapsto (z, z^{-1})$. Recall the notation $H=\phi_*\Gm$ from Section~\ref{Sect_the_dual_pair_GL_GO_1.5}, so $H$ is not to be confused with $\tilde H$.

Recall the stack $\Bun_{G,\tilde H}$ and the map $\tilde\tau: \Bun_{G,\tilde H}\to\Bunt_{G_4}$ from Section~\ref{Sect_4.1.3_G_tildeH}. We use some notations from \cite[Section~6]{Ly1}. In particular, the group scheme $R_{\phi}$ on $X$ is defined by the exact sequence $1\to \Gm\to \phi_*\Gm\times \phi_*\Gm\to R_{\phi}\to 1$, where the first map is $z\mapsto (z, z^{-1})$. 
 
  Let $_{\phi}\GL_2$ be the group scheme on $X$ of automorphisms of $\phi_*\cO_Y$. Define $_{\phi}\tilde H$ by the exact sequence $1\to\Gm\to (_{\phi}\GL_2)\times(_{\phi}\GL_2)\to {_{\phi}\tilde H}\to 1$, where the first map is $z\mapsto (z, z^{-1})$. Recall from \cite[Section~6.1.1]{Ly1} that $\Bun_{_{\phi}\tilde H}\,\iso\,\Bun_{\tilde H}$ naturally. The natural map $R_{\phi}\to {_{\phi}\tilde H}$ induces a morphism denoted $\gq_{R_{\phi}}: \Bun_{R_{\phi}}\to\Bun_{\tilde H}$ in (\select{loc.cit.}, p. 414). 
  
\sssec{}  By \cite[Lemma~16]{Ly1}, the stack $\Bun_{R_{\phi}}$ can be seen as the stack classifying $\cB_1,\cB_2\in\Pic Y$, an isomorphism $\gamma: N(\cB_1)\,\iso\,N(\cB_2)$, and its refinement $\gamma_{12}: \cB_1\mid_{D_Y}\,\iso\, \cB_2\mid_{D_Y}$. This means that $\gamma_{12}^2$ is the restriction of $\gamma$ to $D$. In these terms, the map $\gq_{R_{\phi}}$ sends $(\cB_1,\cB_2,\gamma, \gamma_{12})$ to 
$$
(V, \cC, \Sym^2 V\to \cC, \und{\gamma}: \det V\,\iso\, \cC^2),
$$ 
where $V\subset (\phi_*\cB_1)\oplus(\phi_*\cB_2)$ is given by the cartesian square
$$
\begin{array}{ccc}
(\phi_*\cB_1)\oplus(\beta_*\cB_2) & \to & \cB_1\mid_{D_Y}\oplus \cB_2\mid_{D_Z}\\
\uparrow && \uparrow\lefteqn{\scriptstyle \id+\gamma_{12}}\\
V & \to &\cB_1\mid_{D_Y},
\end{array}
$$  
$\cC=N(\cB_1)(-D)$, and the quadratic form $\Sym^2 V\to\cC$ is the restriction of the difference of the forms on $\phi_*\cB_i$ (this is explained in \cite[Section~6.2]{Ly1}). Tensoring the canonical isomorphisms $\det(\phi_*\cB_i)\,\iso\, \cE_{\phi}\otimes\cC(D)$, one gets an isomorphism 
$$
\det(\phi_*\cB_1\oplus\phi_*\cB_2)\,\iso\, \cC^2(D)
$$ 
Since $\div((\phi_*\cB_1\oplus\phi_*\cB_2)/V)=D$, it yields an isomorphism $\und{\gamma}: \det V\,\iso\, \cC^2$. 

\sssec{} Let $p_{\phi}: \Bun_{R_{\phi}}\to\Pic Y$ be the map sending the above point to $\cB_1$. The product $\phi_*\Gm\times\phi_*\Gm\to\phi_*\Gm$ factors through $R_{\phi}\to \phi_*\Gm$, and $p_{\phi}$ is the extension of scalars with respect to the latter map. 

 Define $\Bun_{G, R_{\phi}}$ by the cartesian square
$$
\begin{array}{cccc}
\Bun_{G, R_{\phi}} & \toup{\id\times\gq_{R_{\phi}}} & \Bun_{G,\tilde H} & \toup{\gp_{\tilde H}} \Bun_G\\ 
\downarrow && \downarrow\lefteqn{\scriptstyle \gq_{\tilde H}}\\
\Bun_{R_{\phi}} & \toup{\gq_{R_{\phi}}} & \Bun_{\tilde H},
\end{array}
$$ 
where $\gq_{\tilde H}, \gp_{\tilde H}$ denote the pojections sending $(M\in\Bun_G, V\in\Bun_{\tilde H})$ to $V$ and $M$ respectively. 

\sssec{} The stack $\Bun_{G,R_{\phi}}$ classifies $M\in \Bun_2$ with $\cA=\det M$ and $\cC=\Omega\otimes\cA^{-1}$, $\cB_1\in\Pic Y$, $\cB_2\in\Pic Y$, isomorphisms $N(\cB_1)\,\iso\, N(\cB_2)\,\iso\, \cC(D)$, and its refinement 
$$
\gamma_{12}: \cB_1\mid_{D_Y}\,\iso\, \cB_2\mid_{D_Y}
$$ 

 Write $\Bun_{G, H, H}$ for the stack classifying $M\in \Bun_2, \cB_1,\cB_2\in\Pic Y$ with isomorphisms $N(\cB_1)\,\iso\, N(\cB_2)\,\iso\, \cC(D)$, where $\cC=\Omega\otimes\cA^{-1}$ and $\cA=\det M$.

 Let $\Bun_{2,D,H,H}$ (resp., $\Bunt_{2,D,H,H}$)
be obtained from $\Bun_{G, H, H}$ by the base change $\delta_D: \Bun_{2,D}\to \Bun_G$ (resp., $\delta_{\tilde D}: \Bunt_{2,D}\to \Bun_G$). Recall that $\delta_{\tilde D}$ sends $(\bar M\subset M, \cU)\in\Bunt_{2,D}$ to $M$. Let 
$$
\tilde\tau_i: \Bunt_{2, D, H, H}\to \Bunt_{G_2}
$$ 
be the map sending $(\cB_1,\cB_2)$, $(\bar M\subset M, \cU)$ to
$\tilde\tau(\bar M\subset M, \cU, \cB_i)$, where $\tilde\tau$ is given by (\ref{map_tilde_tau_Section_221}). 

\sssec{} Consider the diagram
$$
\Bunt_{G_2}\times\Bunt_{G_2} \getsup{\tilde\tau_1\times\tilde\tau_2}
\Bunt_{2,D, H,H}\; \toup{\delta_{\tilde D}}\; \Bun_{G, H, H}\;\getsup{\xi_{\phi}}\; \Bun_{G, R_{\phi}}\toup{\id\times \gq_{R_{\phi}}}\Bun_{G, \tilde H}\toup{\tilde\tau} \Bunt_{G_4}
$$
where $\xi_{\phi}$ is the map forgetting $\gamma_{12}$.

 We need the following corollary of Proposition~\ref{Pp_local_nature_gamma_forgetting}.
\begin{Pp} 
\label{Pp_xi_phi_direct_image_of_theta-sheaf}
There is an isomorphism
\begin{equation}
\label{iso_Pp_10_over_Bun_GHH}
\delta_{\tilde D !}(\tilde\tau_1^*\Aut\otimes\tilde\tau_2^*\Aut)[\dimrel(\tilde\tau_1\times\tilde\tau_2)]\,\iso\, \xi_{\phi !}(\id\times \gq_{R_{\phi}})^*\tilde\tau^*\Aut[\dimrel\tilde\tau\comp(\id\times \gq_{R_{\phi}})]
\end{equation}
\end{Pp}
\begin{proof} Consider the map $\Bun_{G,H,H}\to \Bun_{G,Q}$ sending $(M, \cB_1,\cB_2)$ to $(M, V^{\perp}\subset V, \cC, \gamma)$ with $V=\phi_*(\cB_1\oplus\cB_2)$, $V^{\perp}=\phi_*(\cB_1(-D_Y)\oplus \cB_2(-D_Y))$, $\cC=N(\cB_i)(-D)$. The morphism $\gamma: \det V\,\iso\,\cC^2(D)$ is the product of the corresponding isomorphisms $\det(\phi_*\cB_i)\,\iso\, \cE_{\phi}\otimes N(\cB_i)$. 

For $(M, \cB_1,\cB_2)\in\Bun_{G,H,H}$ a choice of a lower modification 
$$
\phi_*(\cB_1(-D_Y)\oplus\cB_2(-D_Y))\subset \bar V\subset \phi_*(\cB_1\oplus\cB_2)
$$ 
such that $\bar V$ (with the induced maps $\Sym^2\bar V\to\cC$ and $\gamma$) is a $\tilde H$-torsor on $X$, is equivalent to a choice of $\gamma_{12}: \cB_1\mid_{D_Y}\,\iso\, \cB_2\mid_{D_Y}$ refining $N(\cB_1)\mid_D\,\iso\, N(\cB_2)\mid_D$. So, our claim follows from Proposition~\ref{Pp_local_nature_gamma_forgetting}. We have used here \cite[Proposition~3]{Ly1} for the trivial cover $X\sqcup X\to X$. 
\end{proof}
\begin{Rem} 
\label{Rem_iota}
The permutation of two factors on $\GL_2\times\GL_2$ induces an involution of $\tilde H$, write $\iota$ for the induced involution on $\Bun_{\tilde H}$. It sends $(V, \cC, \Sym^2 V\to \cC, \gamma: \det V\,\iso\, \cC^2)$ to $(V, \cC, \Sym^2 V\to \cC, -\gamma: \det V\,\iso\, \cC^2)$. The permutation of two factors on $\phi_*\Gm\times\phi_*\Gm$ induces an involution on $R_{\phi}$. Write also $\iota$ for the corresponding involution on $\Bun_{R_{\phi}}$. It sends $$
(\cB_1,\cB_2, N(\cB_1)\,\iso\, N(\cB_2), \cB_1\mid_{D_Y}\,\iso\, \cB_2\mid_{D_Y})
$$ 
to $(\cB_1,\sigma^*\cB_2, N(\cB_1)\,\iso\, N(\sigma^*\cB_2), \cB_1\mid_{D_Y}\,\iso\, \sigma^*\cB_2\mid_{D_Y})$. The map $\gq_{R_{\phi}}: \Bun_{R_{\phi}}\to\Bun_{\tilde H}$ is equivariant for these involutions. Let $\iota$ act on $\Bun_{G,H,H}$ sending $(M,\cB_1,\cB_2)$ to $(M, \cB_1, \sigma^*\cB_2)$. We will use the fact that (\ref{iso_Pp_10_over_Bun_GHH}) is naturally equivariant with respect to these involutions $\iota$. 
\end{Rem}
 
\sssec{} Denote by $F_{\tilde H}: \D(\Bun_2)\to \D(\Bun_{\tilde H})$ the theta-lifting functor introduced in \cite[Definition~1]{Ly1}. 

 Let $E$ be a rank 2 irreducible local system on $X$. Recall that $\Aut_E$ denotes the corresponding automorphic perverse sheaf on $\Bun_2$ normalized as in \cite[Definition~8]{Ly1}. In particular, it has central character $\det E$. 
In \cite[Section~5.1]{Ly1} we introduced a perverse sheaf $K_{\pi^*E, \det E, \tilde H}$ on $\Bun_{\tilde H}$, here $\pi: X\sqcup X\to X$ is the trivial covering. The extension of scalars map corresponding to (\ref{seq_def_GO^0_4}) is denoted $\rho_{\tilde H}: \Bun_2\times\Bun_2\to\Bun_{\tilde H}$. One has  canonically
$$
\rho_{\tilde H}^*(K_{\pi^*E, \det E, \tilde H})[\dimrel(\rho_{\tilde H})]\,\iso\, \Aut_E\boxtimes\Aut_E
$$

\begin{Pp} 
\label{Pp_11_great}
For an irreducible rank 2 local system $E$ on $X$ one has
\begin{equation}
\label{iso_before_RS_convolution}
p_{\phi !}\gq_{R_{\phi}}^*F_{\tilde H}(\Aut_{E^*})[\dimrel(\gq_{R_{\phi}})]\,\iso\, F_H(\Eis(\Qlb\oplus \cE_0)\otimes \delta_{\tilde D}^*\Aut_{E^*})[-\dim\Bun_2],
\end{equation}
where $F_H$ is the functor from Definition~\ref{Def_theta_lifting_functors_GO_2_GL_2_ram}.  
\end{Pp}
\begin{Prf} Recall the stacks $\Bun_{2,D,H}$ and $\Bunt_{2,D,H}$
defined in Sections~\ref{Sect_3.3.1_should_be} and \ref{Sect_3.3.5}.
Consider the commutative diagram
$$
\begin{array}{ccccccc}
\Bunt_{G_4} & \getsup{\tilde\tau}&\Bun_{G,\tilde H} & \getsup{\id\times\gq_{R_{\phi}}} & \Bun_{G, R_{\phi}} & \toup{\pr_R} & \Bun_{R_{\phi}}\\
&&\downarrow\lefteqn{\scriptstyle\gp_{\tilde H}}
&&\downarrow\lefteqn{\scriptstyle \xi_{\phi}} && \downarrow\lefteqn{\scriptstyle p_{\phi}}\\
&&\Bun_G &\getsup{\pr_G} &\Bun_{G, H, H} & \toup{\nu_{\phi}} & \Pic Y\\
&&&&\uparrow\lefteqn{\scriptstyle \delta_{\tilde D}} && \uparrow\lefteqn{\scriptstyle \gq}\\
&&\Bunt_{G_2}\times\Bunt_{G_2}&\getsup{\tilde\tau_1\times\tilde\tau_2}
&\Bunt_{2,D, H, H} & \toup{\pr_1} & \Bunt_{2,D, H},
\end{array}
$$
where $\nu_{\phi}$ sends $(M, \cB_1,\cB_2)$ to $\cB_1$, and $\pr_1$ sends $(\cU, \bar M\subset M, \cB_1, \cB_2)$ to $(\cU, \bar M\subset M, \cB_1)$. The map $\gq$ here is that of Definition~\ref{Def_theta_lifting_functors_GO_2_GL_2_ram}, the map $\pr_R$ sends $(M, \cB_1,\cB_2,\gamma_{12})$ to $(\cB_1,\cB_2,\gamma_{12})$. By definition, the left hand side of (\ref{iso_before_RS_convolution}) identifies with
$$
p_{\phi !}\pr_{R !}(\id\times\gq_{R_{\phi}})^*(\gp_{\tilde H}^*\Aut_{E^*}\otimes\tilde\tau^*\Aut)[-\dim\Bun_{G_4}+\dimrel(\gp_{\tilde H}\comp(\id\times\gq_{R_{\phi}}))]
$$
Applying Proposition~\ref{Pp_xi_phi_direct_image_of_theta-sheaf}, the above complex is identified with
$$
\nu_{\phi !}(\pr_G^*\Aut_{E^*}\otimes (\delta_{\tilde D})_!(\tilde\tau_1^*\Aut\otimes \tilde\tau_2^*\Aut))[\dimrel(\tilde\tau_1\times\tilde\tau_2)-\dim\Bun_G]
$$
Now applying Proposition~\ref{Pp_mainresults_2} ii), one gets the desired isomorphism.
\end{Prf} 

\medskip\noindent
\begin{Rem} Since $p_{\phi}$ is invariant under the involution $\iota$ introduced in Remark~\ref{Rem_iota}, $\iota$ acts on the left hand side of (\ref{iso_before_RS_convolution}). It acts on the right hand side of (\ref{iso_before_RS_convolution}) via its action on $\Eis(\Qlb\oplus\cE_0)$ described in Remark~\ref{Rem_iota_very_first}.
\end{Rem}

\ssec{Proof of Theorem~\ref{Th_3}}

\sssec{} Since $\Gm\subset \Gm\times\Gm$ via $z\mapsto (z,z^{-1})$, we get a natural map $R_{\phi}\to U_{\phi}\times U_{\phi}$. It yields an extension of scalars map denoted  $\kappa_{\phi}: \Bun_{R_{\phi}}\to \Bun_{U_{\phi}}\times\Bun_{U_{\phi}}$.

 By \cite[Remark~10]{Ly1}, the condition ($C_W$) of Definition~\ref{Def_WP} is satisfied. So, there exists a complex $\cK=\cK_{E,\cJ}\in \D(\Bun_{U_{\phi}})$ equipped with  
$$
e_{\phi}^*\cK[\dimrel(e_{\phi})]\,\iso\, (A\cJ)^{-1}\otimes\phi_1^*\Aut_E[\dimrel(\phi_1)]
$$
By \cite[Remark~10]{Ly1}, $\cK$ is a direct sum of (possibly shifted) irreducible perverse sheaves. From the diagram
$$
\begin{array}{ccccc}
\Bun_2\times\Bun_2 & \getsup{\phi_1\times\phi_1} & \Pic Y\times\Pic Y & \toup{e_{\phi}\times e_{\phi}} & \Bun_{U_{\phi}}\times\Bun_{U_{\phi}}\\
\downarrow\lefteqn{\scriptstyle \rho_{\tilde H}} && \downarrow & \nearrow\lefteqn{\scriptstyle\kappa_{\phi}}\\
\Bun_{\tilde H} & \getsup{\gq_{R_{\phi}}} & \Bun_{R_{\phi}}
\end{array}
$$
we get an isomorphism
$$
\kappa_{\phi}^*(\cK\boxtimes\cK)[\dimrel(\kappa_{\phi})]
\,\iso\, \gq_{R_{\phi}}^*K_{\pi^*E, \det E,\tilde H}\otimes p_{\phi}^*(A\cJ)^{-1}[\dimrel(\gq_{R_{\phi}})]
$$
\sssec{} Consider the commutative diagram of homomorphisms of groups
$$
\begin{array}{ccc}
R_{\phi} & \to & U_{\phi}\times U_{\phi}\\
\downarrow && \downarrow\lefteqn{\scriptstyle m}\\
\phi_*\Gm & \to & U_{\phi},
\end{array}
$$
where $m$ is the product, and the left vertical arrow is induced by the product map $\phi_*\Gm\times\phi_*\Gm\to \phi_*\Gm$. It yields a cartesian square of morphisms of extension of scalars
$$
\begin{array}{ccc}
\Bun_{R_{\phi}} & \toup{\kappa_{\phi}} & \Bun_{U_{\phi}}\times\Bun_{U_{\phi}}\\
\downarrow\lefteqn{\scriptstyle p_{\phi}} &&  \downarrow\lefteqn{\scriptstyle \mult}\\
\Pic Y & \toup{e_{\phi}} & \Bun_{U_{\phi}}
\end{array}
$$
So, one gets an isomorphism
\begin{equation}
\label{iso_inside_Th3_first}
e_{\phi}^*\mult_!(\cK\boxtimes\cK)[\dimrel(e_{\phi})]\,\iso\, (A\cJ)^{-1}\otimes (p_{\phi})_!\gq_{R_{\phi}}^*K_{\pi^*E, \det E, \tilde H}[\dimrel(\gq_{R_{\phi}})] 
\end{equation}
\sssec{} By \cite[Proposition~8]{Ly1}, one has
$$
F_{\tilde H}(\Aut_{E^*})\,\iso\, A(\det E)_{\Omega}\otimes K_{\pi^*E, \det E, \tilde H}
$$
So, using (\ref{iso_before_RS_convolution}), the right hand side of (\ref{iso_inside_Th3_first}) identifies with
$$
(A\cJ)^{-1}\otimes(A\det E)_{\Omega}^{-1}\otimes F_H(\Eis(\Qlb\oplus \cE_0)\otimes \delta_{\tilde D}^*\Aut_{E^*})[-\dim\Bun_2]
$$
Using Lemma~\ref{Lm_for_RS_convolution_preparatory} and Proposition~\ref{Pp_RS_convolution}, the latter complex is identified with 
\begin{multline*}
(A\cJ)^{-1}\otimes(A\det E)_{\Omega}^{-1}\otimes(\mathop{\oplus}\limits_{r\ge 0} F_H(\Av^r_{E^*}(\La^{d_1}_{\det E^*})))\,\iso\, \\
\mathop{\oplus}\limits_{r\ge 0} \;(A\cJ)^{-1}\otimes\varepsilon_!\mult^{r,d_1}_!((A\det E)^{-1}\boxtimes (\phi^*E^*)^{(r)})[g-1+r],
\end{multline*}
where $d_1$ is a function of a connected component of $\Pic Y$ given by $\deg(\cB^*\otimes\Omega_Y)=2d_1+r$ for $\cB\in\Pic Y$. The sum is over $r\ge 0$ such that $d_1\in\ZZ$. 

\sssec{} For any $r\ge 0, d_1\in\ZZ$ one has the cartesian square
$$
\begin{array}{ccc}
\Pic^{d_1} X\times Y^{(r)} & \toup{\mult^{r, d_1}} & \Pic^{2d_1+r} Y\\
\downarrow && \downarrow\lefteqn{\scriptstyle e_{\phi}\comp \varepsilon}\\
Y^{(r)} & \toup{m_{\phi, r}} & \Bun_{U_{\phi}}
\end{array}
$$ 
Since $\dimrel(e_{\phi})=g-1$, this finally yields isomorphisms
\begin{multline*}
(e_{\phi}\comp\varepsilon)^*\mult_!(\cK\boxtimes\cK)\,\iso\,
\mathop{\oplus}\limits_{r\ge 0} \; (A\cJ^{-1})_{\Omega_Y}\otimes (A\cJ)\otimes \mult^{r,d_1}_!((A\det E)^{-1}\boxtimes (\phi^*E^*)^{(r)})[r]\,\iso\, 
\\
\mathop{\oplus}\limits_{r\ge 0} \; (A\cJ^{-1})_{\Omega_Y}\otimes (e_{\phi}\comp\varepsilon)^*(m_{\phi, r})_!(\cJ\otimes\phi^*E^*)^{(r)}[r]
\end{multline*}
over $\Pic^m Y$ for any $m$. It is compatible with the descent data for $e_{\phi}\comp\varepsilon$. So, 
$$
\mult_!(\cK\boxtimes\cK)\,\iso\, \mathop{\oplus}\limits_{r\ge 0} \; (A\cJ^{-1})_{\Omega_Y}\otimes(m_{\phi, r})_!(\cJ\otimes\phi^*E^*)^{(r)}[r],
$$
the sum being over all $r\ge 0$. For $r$ even (resp., odd) $m_{\phi, r}$ maps $Y^{(r)}$ to $\Bun^0_{U_{\phi}}$ (resp., to $\Bun^1_{U_{\phi}}$). 

 If $\phi^*E$ is irreducible then $\RG(Y^{(r)}, (\cJ\otimes\phi^*E^*)^{(r)})[r]\,\iso\, \wedge^r V$ with $V=\H^1(Y, \cJ\otimes\phi^*E^*)$. This concludes the proof of Theorem~\ref{Th_3}.
 
\medskip\noindent
\begin{Rem} 
\label{Rem_great}
Let the involution $\iota$ act on $\Bun_{U_{\phi}}\times\Bun_{U_{\phi}}$ permuting the factors. Since $\Bun_{U_{\phi}}$ is a commutative group stack, $\iota$ acts on $\mult_!(\cK\boxtimes\cK)$. For the action of $\iota$ on $\Bun_{R_{\phi}}$ defined in Remark~\ref{Rem_iota} the map $\kappa_{\phi}: \Bun_{R_{\phi}}\to\Bun_{U_{\phi}}\times\Bun_{U_{\phi}}$ is $\iota$-equivariant. So, (\ref{iso_inside_Th3_first}) is $\iota$-equivariant. We conclude that the grading by the action of $\iota$ on 
$$
\mathop{\oplus}\limits_{r\ge 0} \; (A\cJ^{-1})_{\Omega_Y}\otimes(m_{\phi, r})_!(\cJ\otimes\phi^*E^*)^{(r)}[r]
$$ 
is given by fixing $r\!\!\!\mod \! 4$. In particular, the $\iota$-invariants of this complex over $\Bun^0_{U_{\phi}}$ are
$$
\mathop{\oplus}\limits_{r\ge 0, \, r= 0\!\!\!\!\mod\! 4} \;\; (A\cJ^{-1})_{\Omega_Y}\otimes(m_{\phi, r})_!(\cJ\otimes\phi^*E^*)^{(r)}[r]
$$
This agrees with Conjecture~\ref{Con_iso_of_algebras}. Indeed, let $n\ge 1$. If $V$ is a $2n$-dimensional $\Qlb$-vector space with a nondegenerate symmetric form, let $S^{\pm}$ denote the half-spin representations of the corresponding Spin group. Then 
$$
\Sym^2(S^+)\oplus \Sym^2(S^-)\,\iso\, \oplus_{k\in \ZZ} \wedge^{n+4k} V
$$ 
and $S^+\otimes S^-\,\iso\, \wedge^{n-1} V\oplus \wedge^{n-3} V\oplus\ldots$. In our case $V=\H^1(Y, \cJ\otimes \phi^*E^*)$ is of dimension divisible by 4. In general, if $\dim V$ is divisible by 4 then for the $\ZZ/2\ZZ$-graded isomorphism 
$$
(S^+\oplus S^-)^{\otimes 2}\,\iso\, \oplus_{r\ge 0}\wedge^r V
$$ 
the additional grading by $r\!\!\mod\! 4$ on the right hand side corresponds to the additional grading of $(S^+\oplus S^-)^{\otimes 2}$ with respect to the involution permuting the two factors in the tensor product.
\end{Rem} 

\section{The dual pair $(\wt\SL_2, \SO_3)$}
\label{Sect_dual_pair_Mp2_SO3}

In this section we perform the constructions and establish the results formulated in Section~\ref{Sect_Quantum_intro}. They are mostly derived from Theorem~\ref{Th_3}.

\ssec{} We use some notations from \cite[Section~0.3]{Ly6}. In particular, $H=\SO_3$ split, $G=\Sp(\cO_X\oplus\Omega)$ a group scheme over $X$, $B\subset G$ is the parabolic group subscheme over $X$ preserving $\Omega$. The stack $\Bun_B$ classifies $\cE\in\Bun_1$ and an exact sequence on $X$
\begin{equation}
\label{ext_cE^-1_by_cEotimesOmega}
0\to\cE\otimes\Omega\to M\to \cE^{-1}\to 0
\end{equation}
The map $\tilde\nu_B: \Bun_B\to\Bunt_G$ sends the above point to $(M,\cB)$, where $\cB=\det\RG(X, \cE\otimes\Omega)$ with the induced isomorphism $\cB^2\,\iso\, \det\RG(X,M)$. 

 The stack $\cS_B$ classifies $\cE\in\Bun_1$ and $s_2: \cE^2\to\cO_X$. We denote by $\Four_{\psi}: \D(\Bun_B)\,\iso\, \D(\cS_B)$ the Fourier transform. One has an open immersion $\RCov^d\hook{} \cS_B$. 
 
 Set $\Bun_{B,H}=\Bun_B\times\Bun_H$ and $\Bunt_{G,H}=\Bunt_G\times\Bun_H$. Recall the map $\tilde\tau: \Bunt_{G,H}\to\Bunt_{G_3}$ from \cite[Section~7]{Ly4}. It sends $(M,\cB, V)$ to $(M\otimes V, \cB')$, where 
$$
\cB'=\cB^3\otimes\det\RG(X, V)\otimes \det\RG(X, \cO)^{-3}
$$ 
with the induced isomorphism $\cB'^2\,\iso\, \det\RG(X, M\otimes V)$. 

\sssec{} Let $\tau_B: \Bun_{B,H}\to\Bun_{P_3}$ be the map sending (\ref{ext_cE^-1_by_cEotimesOmega}) and $V$ to $(\cE\otimes\Omega\otimes V\subset M\otimes V)$. The following diagram is 2-commutative
$$
\begin{array}{ccc}
\Bunt_{G,H} & \toup{\tilde\tau} & \Bunt_{G_3} \\
\uparrow\lefteqn{\scriptstyle \tilde\nu_B\times\id} && \uparrow\lefteqn{\scriptstyle \tilde\nu_{P_3}}\\
\Bun_{B,H} & \toup{\tau_B} & \Bun_{P_3},
\end{array}
$$
here $\tilde\nu_{P_3}$ is the map from Section~\ref{Sect_2.1.2}.  
Using \cite[Proposition~1]{Ly1}, from the above diagram one gets the following.

 Let $\cV_{H,B}$ be the stack classifying $\cE\in\Bun_1$, $V\in\Bun_H$ and a section $s: \cE\to V$. Let $\gp_{\cV}: \cV_{H,B}\to \cS_B$ be the map sending this point to $(\cE, s_2)$, where $s_2$ is the composition $\cE^2\toup{s\otimes s} \Sym^2 V\to \cO$. We get a diagram
$$
\Bun_H\, \getsup{\gq_{\cV}}\,\cV_{H,B}\,\toup{\gp_{\cV}}\, \cS_B,
$$
where $\gq_{\cV}$ sends $(E,V,s)$ to $V$. Let $F_G: \D^-(\Bun_H)_!\to\D^{\prec}(\Bunt_G)$ be the theta-lifting functor from \cite[Section~0.3.2]{Ly6}. 
\begin{Lm} 
\label{Lm_description_F_cS_over_RCov^d}
There is an isomorphism of functors $\D^-(\Bun_H)_!\to \D(\cS_B)$
$$
\gp_{\cV !} \gq_{\cV}^*(\cdot)[\dimrel(\gq_{\cV})]\,\iso\, \Four_{\psi}\tilde\nu_B^*F_G(\cdot)[\dimrel(\tilde\nu_B)]
$$
\end{Lm}
\begin{Rem} According to \cite[Proposition~7]{Ly4}, once $\sqrt{-1}\in k$ is fixed, the isomorphism of Lemma~\ref{Lm_description_F_cS_over_RCov^d} is defined uniquely up to a sign.
\end{Rem}

\sssec{} 
\label{Sect_6.1.4}
For $(\cE, s_2)\in \RCov^d$ we have the corresponding degree 2 covering $\phi: Y\to X$, here $Y=\Spec(\cO_X\oplus\cE)$. Recall the group scheme $U_{\phi}$ on $X$ from Section~\ref{Sect_Geometric_Waldspurger_intro}. 
As $(\cE, s_2)$ varies in $\RCov^d$ the groups $U_{\phi}$ form a group scheme $U_R$ over $\RCov^d\times X$. Let $\Bun_{U_R}$ be the stack whose $S$-points is the category of pairs: a map $\nu: S\to\RCov^d$ and a $(\nu\times\id)^*U_R$-torsor over $S\times X$. 

\begin{Lm} 
\label{Lm_another_cart_square_Section51}
There is a cartesian square, where the horizontal maps are open immersions
$$
\begin{array}{ccc}
\Bun_{U_R} & \hook{} & \cV_{H,B}\\
\downarrow\lefteqn{\scriptstyle \gp_R} && \downarrow\lefteqn{\scriptstyle \gp_{\cV}}\\
\RCov^d & \hook{} & \cS_B
\end{array}
$$
\end{Lm}
\begin{Prf}
Let $(\cE, V, s)\in \cV_{H,B}$ such that $\gp_{\cV}(\cE, V, s)\in\RCov^d$. Let $Y=\Spec(\cO_X\oplus\cE)$ and $\phi: Y\to X$ the corresponding covering ramified at a multiplicity free divisor $D\in X^{(d)}$. Clearly, $s: \cE\to V$ is a subbundle, let $L$ denote the kernel of $s^*: V\to \cE^*$. The symmetric form $\Sym^2 V\to\cO_X$ induces a generically nondegenerate symmetric form $\Sym^2 L\to\cO_X$ and a compatible isomorphism $\det L\,\iso\, \cE$. By \cite[Proposition~13]{Ly2}, this yields a line bundle $\cB$ on $Y$ and isomorphisms $L\,\iso\, \phi_*\cB$, $N(\cB)\,\iso\, \cO_X$ such that the symmetric form on $L$ is the canonical one $\Sym^2(\phi_*\cB)\to N(\cB)$. 

Let $D_Y\in Y^{(d)}$ be the ramification divisor of $\phi$, note that $\phi^*\cE\,\iso\, \cO(-D_Y)$ canonically. The upper modification $(\cE^{-1}\otimes L)\oplus\cO\hook{} \cE^{-1}\otimes V$ gives rise to the compatible trivialization $\gamma: \cB\mid_{D_Y}\,\iso\, \cO_{D_Y}$. Namely, 
\begin{equation}
\label{iso_geom_fibre_D_phi} 
((\cE^{-1}\otimes L)\oplus\cO)\mid_D\,\iso\, (\cB(D_Y)/\cB(-D_Y))\oplus \cO_D,
\end{equation}
and the image of $\cE^{-1}\otimes V(-D)$ in (\ref{iso_geom_fibre_D_phi}) is the graph of $\gamma$. So, $(\cB, \gamma, N(\cB)\,\iso\, \cO_X)\in\Bun_{U_{\phi}}$.
\end{Prf}

\sssec{} Let $\gq_U$ be the composition $\Bun_{U_R}\hook{}\cV_{H,B}\toup{\gq_{\cV}}\Bun_H$. Consider the stack $^{rss}\cS^d_P$ classifying $U\in \Bun_2, \cA\in\Bun_1$ and $\Sym^2 U\to\cA\otimes\Omega$ such that $U\hook{} U^*\otimes\cA\otimes\Omega$, and $\div(U^*\otimes\cA\otimes\Omega/U)\in {^{rss}X^{(d)}}$. We have the diagram 
\begin{equation}
\label{diag_e_U_square}
\begin{array}{ccccc}
\Bun_2 & \getsup{\gq_1} & ^{rss}\cS^d_P & \toup{\gp_1}&\RCov^d \\
\downarrow\lefteqn{\scriptstyle e_H} && \downarrow\lefteqn{\scriptstyle e_U}& \nearrow\lefteqn{\scriptstyle \gp_R}\\
\Bun_H & \getsup{\gq_U} & \Bun_{U_R},
\end{array}
\end{equation}
here the top line is the diagram considered in \cite[beginning of Section~6.1.2]{Ly1}. So, $\gq_1$ sends the above point to $U$. The map $\gp_1$ sends the above point to $\cE=(\cA\otimes\Omega)^{-1}\otimes\det U$ with the induced inclusion $\cE^2\hook{} \cO_X$. The map $e_H$ is the extension of scalars with respect to $\GL_2\to H$, so it sends $U$ to $V=\END_0(U)$ with the induced symmetric form $\Sym^2 V\to\cO_X$ and a compatible trivialization $\det V\,\iso\,\cO_X$. 

 The stack $^{rss}\cS_P^d$ can be seen as the one classifying a degree 2 covering $\phi: Y\to X$ ramified exactly at $D\in {^{rss}X^{(d)}}$ and a line bundle $\cB_0$ on $Y$. In these terms $\gq_1$ is the map sending $\cB_0$ to $\phi_*\cB_0$. The map $e_U$ preserves the covering $\phi$ and sends $\cB_0$ to $\cB=\sigma^*\cB_0\otimes\cB_0^{-1}$ with the induced isomorphisms $\gamma: \cB\mid_{D_Y}\,\iso\, \cO\mid_{D_Y}$ and $N(\cB)\,\iso\, \cO_X$. Here $\sigma$ is the nontrivial automorphism of $Y$ over $X$.

 Note that the fibre of $e_U$ over a $k$-point of $\RCov^d$ given by 
$\phi: Y\to X$ is the morphism $e_{\phi}: \Pic Y\to \Bun_{U_{\phi}}$ considered in Section~\ref{Sect_5.1.1_set-up}. 

\begin{Lm} The diagram (\ref{diag_e_U_square}) commutes.
\end{Lm}
\begin{Prf}
Consider a point of $\RCov^d$ given by $\phi: Y\to X$ ramified over $D\in {^{rss}X^{(d)}}$. As above, $\cE$ denotes the $\sigma$-antiinvariants in $\phi_*\cO_Y$, $\sigma$ is the nontrivial automorphism of $Y$ over $X$. Let $\cB_0\in\Pic Y$ and $U=\phi_*\cB_0$. 
The sequence 
$$
0\to \phi_*(\cB_0^2(-D_Y))\to\Sym^2 U\to N(\cB_0)\to 0
$$ 
is exact and $\det U\,\iso\, \cE\otimes N(\cB_0)$ canonically. This yields an exact sequence
\begin{equation}
\label{seq_for_END_0_U}
0\to L\to 
\END_0(U)\to \cE^{-1}\to 0
\end{equation}
with $L=\phi_*(\cB_0\otimes\sigma^*\cB_0^{-1})$. The symmetric form $\Sym^2 U\to N(\cB_0)$ has a canonical section $s: N(\cB_0)(-D)\to \Sym^2 U$ defined as follows. One has a cartesian square
$$
\begin{array}{ccc}
N(\cB_0)\oplus \phi_*(\cB_0^2) & \to & \cB_0^2\mid_{D_Y}\oplus \cB_0^2\mid_{D_Y}\\
\uparrow && \uparrow\lefteqn{\scriptstyle \diag}\\
\Sym^2 U & \to & \cB^2_0\mid_{D_Y}
\end{array}
$$
Then $s$ is defined uniquely by the property that the composition $N(\cB_0)(-D)\to\Sym^2 U\to N(\cB_0)\oplus \phi_*(\cB^2_0)$ equals $(i, 0)$, where $i$ is the canonical inclusion. Now $s$ can be seen as $\cE\hook{} \END_0(U)$, so we get an upper modification $(\cE^{-1}\otimes L)\oplus\cO\hook{} \cE^{-1}\otimes\END_0(U)$. This shows that $\END_0(U)\,\iso\, V$ naturally, where $V$ is the image of $(\cB, \gamma, N(\cB)\,\iso\, \cO_X)$ under $\gq_U$. To summarize, $V$ fits into the diagram, where the square is cartesian
$$
\begin{array}{ccc}
\cE^{-1}\oplus \phi_*(\cB_0(D_Y)\otimes\sigma^*\cB_0^{-1}) & \to & \cO(D_Y)/\cO\oplus \cO(D_Y)/\cO\\
\uparrow && \uparrow\lefteqn{\scriptstyle \diag}\\
V & \to & \cO(D_Y)/\cO\\
\uparrow\\
\cE\oplus \phi_*(\cB_0\otimes\sigma^*\cB_0^{-1})
\end{array}
$$
\end{Prf}

From Lemmas~\ref{Lm_another_cart_square_Section51} and \ref{Lm_description_F_cS_over_RCov^d} one immediately derives the following.

\begin{Cor} 
\label{Cor1_for_Section5.2}
There is an isomorphism of functors $\D^-(\Bun_H)_!\to \D(\RCov^d)$
$$
\gp_{R !} \gq_{U}^*(\cdot)[\dimrel(\gq_U)]\,\iso\, \Four_{\psi}\tilde\nu_B^*F_G(\cdot)[\dimrel(\tilde\nu_B)]\mid_{\RCov^d}
$$
\QED
\end{Cor}

\sssec{Proof of Proposition~\ref{Pp_mainresults_3}} 
\label{Sect_6.1.9}
The map $e_U: {^{rss}\cS^d_P}\to \Bun_{U_R}$ is smooth and surjective. Since the ULA property is local in the smooth topology of the source, the first claim follows from \cite[Theorem~3]{Ly1}. For the second claim apply Corollary~\ref{Cor1_for_Section5.2}. Since $\gp_R$ is proper, 
$p_{R !}\gq_U^*K$ is ULA with respect to the identity map $\RCov^d\to\RCov^d$. So, each cohomology sheaf of (\ref{complex_forCor1_on_RCov^d}) is a local system. \QED

\sssec{} Write $^{rss}\cS^{d,m}_P$ for the open substack of $^{rss}\cS^d_P$ given by $\deg(\cA\otimes\Omega)=m$. The image of $^{rss}\cS^{d,m}_P$ under $e_U$ depends only on $a=m\!\!\mod \! 2$ and is denoted $\Bun_{U_R}^a$ for $a\in \ZZ/2\ZZ$. So, $\Bun_{U_R}$ is a disjoint union of the open substacks $\Bun_{U_R}^a$ for  $a\in \ZZ/2\ZZ$.

\sssec{Proof of Corollary~\ref{Cor_main_for_SL_2_SO_3}} 
\label{Sect_6.1.11}
From Theorem~\ref{Th_3} and Corollary~\ref{Cor1_for_Section5.2} one gets the desired isomorphism. 

 As in \cite[Theorem~3]{Ly1}, one shows that $\gq_1^*\Aut_E[\dimrel(\gq_1)]$ is Verdier self-dual for any $d$, this heavily relies on the Hecke property of $\Aut_E$. For $d>4g-4$ this also follows from the fact that the map $\gq_1$ is smooth by \cite[Lemma~13]{Ly1}. Now from Corollary~\ref{Cor1_for_Section5.2} we learn that (\ref{complex_forCon4_on_RCov^d}) is Verder self-dual. 
 
If $d>0$ then for any $(\cE, s_2)\in \RCov^d$ and for the corresponding $\cE_0$ we get $\H^0(X, E\otimes\cE_0)=0$, so (\ref{complex_forCon4_on_RCov^d}) is a local system. \QED

\ssec{The trace map}  

\sssec{} Let $E$ be an irreducible $\SL_2$-local system on $X$ and $d\ge 0$. Over the open substack of $\RCov^d$ given by $\phi: Y\to X$ such that $\H^0(Y, \phi^*E)=0$ the sheaf $\cF_E$ given by (\ref{complex_forCon4_on_RCov^d}) is a local system. In particular, for $d>0$ this is a local system. Les us describe the trace map 
$$
\tr: CL_E(X)\otimes CL^d_E\,\iso\,\END(\cF_E)\to\Qlb
$$ 
over this open substack of $\RCov^d$. More precisely, we will describe the fibre $\tr_{\phi}$ of $\tr$ at a point $\phi: Y\to X$ of $\RCov^d$.

 Consider the map $f_{\phi}:\Bun_{U_{\phi}}\to \Bun_{R_{\phi}}$ sending $\cB$ with $N(\cB)\,\iso\, \cO_X$, $\cB\mid_{D_Y}\,\iso\, \cO_{D_Y}$ to $(\cB_1=\cB^{-1}\otimes\Omega_Y, \cB_2=\Omega_Y)$ with the induced isomorphisms $N(\cB_1)\,\iso\, N(\cB_2)$ and $\cB_1\mid_{D_Y}\,\iso\, \cB_2\mid_{D_Y}$. 
The composition $\kappa_{\phi}f_{\phi}$ sends $\cB$ to $(\cB, \cB)\in\Bun_{U_{\phi}}\times\Bun_{U_{\phi}}$. 

 Let $\mu_{\phi}: \Bun_{U_{\phi}}\to\Pic Y$ be the map $\cB\mapsto \cB^{-1}\otimes\Omega_Y$. Let $E$ be any irreducible rank 2 local system on $X$ (we will specialize later to $\SL_2$-local systems). One has a natural map over $\Bun_2$
$$
\Aut_E\otimes\Aut_{E^*}\to\Qlb[2\dim\Bun_2]
$$ 
 
\sssec{} Let $\zeta_{\tilde H}: \Bun_{\SO_3}\to\Bun_{\tilde H}$ be the map $V_1\mapsto (V_1\oplus\cO)\otimes\Omega=V$, here $V$ is equipped with the induced form $\Sym^2 V\to\Omega^2$ and trivialization $\det V\,\iso\, \Omega^4$. The above map yields a morphism
$\zeta_{\tilde H}^*K_{\pi^*E, \det E, \tilde H}\to \Qlb[1-g+2\dim\Bun_2]$ over $\Bun_{\SO_3}$. We have a commutative diagram
$$
\begin{array}{ccc}
\Bun_{U_{\phi}} & \to & \Bun_{\SO_3}\\
\downarrow\lefteqn{\scriptstyle f_{\phi}} && \downarrow\lefteqn{\scriptstyle \zeta_{\tilde H}}\\
\Bun_{R_{\phi}} & \toup{\gq_{R_{\phi}}} & \Bun_{\tilde H}
\end{array}
$$ 
This yields a map 
$$
f_{\phi}^*\gq_{R_{\phi}}^*K_{\pi^*E, \det E, \tilde H}\to \Qlb[1-g+2\dim\Bun_2]
$$ 
and in turn, a map
$$
\mu_{\phi}^*p_{\phi !}\gq_{R_{\phi}}^*K_{\pi^*E, \det E, \tilde H}\to\Qlb[1-g+2\dim\Bun_2]
$$
Since $\dimrel(\gq_{R_{\phi}})=2\dim\Pic Y-2\dim\Bun_2$, this is a map
$$
\mu_{\phi}^*p_{\phi !}\gq_{R_{\phi}}^*K_{\pi^*E, \det E, \tilde H}[\dimrel(\gq_{R_{\phi}})]\to\Qlb[2\dim\Pic Y+1-g]
$$
Since $F_{\tilde H}(\Aut_{E^*})\,\iso\, A(\det E)_{\Omega}\otimes K_{\pi^*E, \det E, \tilde H}$, we get a morphism over $\Bun_{U_{\phi}}$
\begin{multline*}
\mu_{\phi}^*F_H(\Eis(\Qlb\oplus \cE_0)\otimes\delta^*_{\tilde D}\Aut_{E^*})[-\dim\Bun_2]\,\iso\, 
\mu_{\phi}^*p_{\phi !}\gq_{R_{\phi}}^*F_{\tilde H}(\Aut_{E^*})[\dimrel(\gq_{R_{\phi}})]\to \\
A(\det E)_{\Omega}[2\dim\Pic Y+1-g]
\end{multline*}
We have
\begin{multline*}
\mu_{\phi}^*F_H(\Eis(\Qlb\oplus \cE_0)\otimes\delta^*_{\tilde D}\Aut_{E^*})[-\dim\Bun_2]\,\iso\\ \mathop{\oplus}\limits_{r\ge 0, \; even} (A\det E)_{\Omega}\otimes\mu_{\phi}^*\varepsilon_!\mult^{r, -r/2}((A\det E)^*\boxtimes (\phi^*E^*)^{(r)})[g-1+r]
\end{multline*}

\sssec{} Consider the stack $\cY_{U_{\phi}}$ classifying pairs:
$$
(\cB, N(\cB)\,\iso\, \cO_X, \cB\mid_{D_Y}\,\iso\, \cO_{D_Y})\in\Bun_{U_{\phi}}
$$ 
and a subsheaf $M_1\subset \phi_*\cB$ of rank 1. Let $\pi_{\phi}: \cY_{U_{\phi}}\to\Bun_{U_{\phi}}$ be the map forgetting $M_1$. For $r\ge 0$ even we have a substack $\cY^r_{U_{\phi}}$ given by the property that the divisor $D_1$ of zeros of $\phi^*M_1\hook{} \cB$ is of degree $r$. Let $h_{\phi}: \cY_{U_{\phi}}^r\to \Bun_1\times Y^{(r)}$ be the map sending the above point to $(M_1, D_1)$. We get a map over $\Bun_{U_{\phi}}$
$$
\mathop{\oplus}\limits_{r\ge 0, \, even}  \mu_{\phi}^*\varepsilon_!\mult^{r, -r/2}((A\det E)^*\boxtimes (\phi^*E^*)^{(r)})[r]\to \Qlb[2\dim\Pic Y+2-2g]
$$
The above map rewrites as
$$
\eta: \mathop{\oplus}\limits_{r\ge 0, \, even} \pi_{\phi !}h_{\phi}^*((A\det E)^*\boxtimes (\phi^*E^*)^{(r)}[r]\to \Qlb[d+2g-2]
$$
Note that $\dim\Bun_{U_{\phi}}=\frac{d}{2}+g-1$, so that the right hand side is the dualizing  complex on $\Bun_{U_{\phi}}$. 

\sssec{}  The map 
$
\Bun_1\times\Bun_{U_{\phi}}\to \coprod_{s \, even} \Pic^s Y, \; (M_1, \cB)\mapsto \phi^*M_1^*\otimes\cB
$ 
is \'etale and surjective. So, the projection 
$$
\cY_{U_{\phi}}\to \coprod_{r\ge 0, \; even} Y^{(r)}
$$ 
sending $(M_1, \cB, (\phi^*M_1)(D_1)\,\iso\, \cB)$ to $D_1$ is \'etale and surjective. Now 
$$
\eta\in \mathop{\oplus}\limits_{r\ge 0, \, even}\H^0(\cY^r_{U_{\phi}},  \DD(h_{\phi}^*((A\det E)^*\boxtimes (\phi^*E^*)^{(r)}[r]))
$$
The sheaf $h_{\phi}^*((A\det E)^*\boxtimes (\phi^*E^*)^{(r)}[r])$ is perverse on $\cY^r_{U_{\phi}}$ and
$$
\DD(h_{\phi}^*((A\det E)^*\boxtimes (\phi^*E^*)^{(r)}[r]))\,\iso\, h_{\phi}^*((A\det E)\boxtimes (\phi^*E)^{(r)}[r])
$$
Assume $\det E$ trivial. Using the projection $\cY^r_{U_{\phi}}\to Y^{(r)}$, which is finite and \'etale, one gets
$$
\eta\in \mathop{\oplus}\limits_{r\ge 0, \, even}\H^0(Y^{(r)}, (\phi^*E)^{(r)}[r])\subset \mathop{\oplus}\limits_{r\ge 0, \, even} \H^0(\cY^r_{U_{\phi}}, h_{\phi}^*(\Qlb \boxtimes (\phi^*E)^{(r)}[r])
$$
That is,
\begin{equation}
\label{map_trace_found_essentially}
\eta\in\mathop{\oplus}\limits_{r\ge 0, \, even} \wedge^r \H^1(Y, \phi^*E)
\end{equation}

\sssec{} Let $\cK$ be defined as in Theorem~\ref{Th_3} 
for $E$ and $\cJ=\Qlb$. Applying $\mu_{\phi}^*$ to (\ref{iso_inside_Th3_first}), from $\eta$ we get a morphism
$$
\mu_{\phi}^*e_{\phi}^*\mult_!(\cK\boxtimes\cK)[\dimrel(e_{\phi})]\to \Qlb[2\dim\Pic Y+1-g] 
$$
Note that $\dimrel(e_{\phi})=g-1$. Applying $\RG$ to the induced map
$$
\mult_!(\cK\boxtimes\cK)\to (e_{\phi}\mu_{\phi})_!\Qlb[2\dim\Pic Y+2-2g] 
$$
we get the trace map 
$$
\tr_{\phi}: \oplus_{r\ge 0} \RG(Y^{(r)}, (\phi^*E)^{(r)})[r]\,\iso\,
\RG(\Bun_{U_{\phi}}, \cK)^{\otimes 2}\to \RG(\Bun_{U_{\phi}}, \DD(\Qlb))\to \Qlb
$$
that coincides with (\ref{map_trace_found_essentially}). In particular, it vanishes on the summands with $r$ odd. 

\begin{Rem} Given a $\Qlb$-vector space $V$ with a nondegenerate symmetric form with $\dim V$ divisible by 4, let $S^{\pm}$ denote the half-spin representation of the Spin covering of $\SO(V)$. The trace map $\oplus_{i\ge 0} \wedge^i V\,\iso\, \End(S^+\oplus S^-)\to\Qlb$ vanishes on $\wedge^i V$ unless $i=0$ or $i=\dim V$. Conjecture~\ref{Con_iso_of_algebras} would imply that the composition of (\ref{iso_maybe_of_algebras!!}) with the trace map 
$\End(\cF_E)\to \Qlb$ vanishes on the corresponding direct summands of $CL_{E}(X)\otimes CL^d_{E}$. This claim could happen easier to check than Conjecture~\ref{Con_iso_of_algebras}.
\end{Rem}

\appendix
\section{Theta-sheaves and two-fold coverings}
\label{Section_appendix}

\ssec{} Here we generalize \cite[Proposition~3]{Ly1} to the case of ramified degree 2 coverings. Keep notations of Section~\ref{Sect_the_dual_pair_GL_GO_1.5}, so we pick $d\ge 0$, a $k$-point $(\cE, s_2)\in\RCov^d$ and denote by $\phi: Y\to X$ the corresponding degree 2 covering ramified over a multiplicity free effective degree $d$ divisor $D$ on $X$. Recall the group scheme $G_n$ on $X$ defined in Section~\ref{Sect_2.1.2}. 

 Let $\Bun_{G_n, Y}$ denote the stack classifying a rank $2n$ vector bundle on 
$Y$ with a symplectic form $\wedge^2 M\to\Omega_Y$. Let $\pi_n: \Bun_{G_n, Y}\to\Bun_{G_{2n}}$ be the following map. Given $M\in\Bun_{G_n, Y}$, we  equip $\phi_*M$ with the form $\wedge^2(\phi_*M)\to\Omega$ given as the composition
$$
(\phi_*M)\otimes(\phi_* M)\to \phi_*(M\otimes M)\to \phi_*(\wedge^2 M)\to \phi_*\Omega_Y\to \Omega,
$$ 
where the latter map is the projection $\phi_*\Omega_Y\,\iso\, \Omega\oplus \Omega\otimes\cE(D)\to \Omega$ on the $\sigma$-invariants in $\phi_*\Omega_Y$. One easily checks that this form is generically non degenerate. This yields an inclusion $\omega: \phi_*M\to (\phi_*M)^*\otimes\Omega$. Note that $\det(\phi_*M)\,\iso\, \cO(-nD)\otimes N(\Omega_Y^n)\,\iso\, \Omega^{2n}$. Thus, $\det\omega=\id$ and $\phi_*M\in \Bun_{G_{2n}}$.

  Let $\cA_{G_n, Y}$ be the line bundle on $\Bun_{G_n, Y}$ with fibre $\det\RG(Y, M)$ at $M\in \Bun_{G_n, Y}$. Let $\Bunt_{G_n, Y}$ be the gerbe of square roots of this line bundle. Let $\tilde\pi_n: \Bunt_{G_n, Y}\to \Bunt_{G_{2n}}$ be the map sending $(M,\cB, \cB^2\,\iso\, \det\RG(Y, M))$ to $(\phi_*M, \cB)$ with the induced isomorphism $\cB^2\,\iso\, \det\RG(X, \phi_*M)$. Recall the decomposition into irreducible perverse sheaves $\Aut\,\iso\,\Aut_g\oplus\Aut_s$ from Section~\ref{Sect_2.1.2}, we refer to these pieces as generic and special part of $\Aut$.  

\begin{Pp} 
\label{Pp_last}
There is a canonical isomorphism $\tilde\pi_n^*\Aut[\dimrel(\tilde\pi_n)]\,\iso\, {_Y\Aut}$ preserving  the generic and special parts. Here $_Y\Aut$ is the theta-sheaf on $\Bunt_{G_n, Y}$. 
\end{Pp}
\begin{proof}
Let $_i\Bun_{G_m}$ denote the stratum of $\Bun_{G_m}$ given by $\dim\H^0(X, M)=i$ for $M\in\Bun_{G_m}$. Let $_i\Bunt_{G_m}$ be the restriction of the gerbe $\Bunt_{G_m}\to\Bun_{G_m}$ to this stratum. Define $_i\Bun_{G_n, Y}$ and $_i\Bunt_{G_n, Y}$ similarly. For each $i$ there is a map $_i\rho: {_i\Bun_{G_m}}\to {_i\Bunt_{G_m}}$ sending $M$ to $(M, \cB=\det\H^0(X, M))$ with the induced isomorphism $\cB^2\,\iso\, \det\RG(X, M)$, and similarly for $Y$. The diagram is 2-commutative
$$
\begin{array}{ccc}
_i\Bun_{G_n, Y} & \to & _i\Bun_{G_{2n}}\\
\downarrow && \downarrow\lefteqn{\scriptstyle {_i\rho}}\\
_i\Bunt_{G_n, Y} & \toup{\tilde \pi_n} & _i\Bunt_{G_{2n}}
\end{array}
$$

 The perverse sheaf $\Aut$ is equipped with trivializations for $i=0,1$$$
 _i\rho^*\Aut\,\iso\, \Qlb[\dim(\Bun_{G_{2n}})-i]
$$ 
This gives the desired isomorphism stratum by stratum and also a canonical normalization of the sought-for isomorphism over $_i\Bunt_{G_n, Y}$ for $i=0,1$. It remains to show its existence. 

Recall the group subscheme $P_n\subset G_n$ defined in Section~\ref{Sect_2.1.2}. Consider the commutative diagram
$$
\begin{array}{ccc}
\Bun_{P_n, Y} & \toup{\tilde\nu_{n, Y}} & \Bunt_{G_n, Y}\\
\downarrow\lefteqn{\scriptstyle \pi_{n,P}} && \downarrow\lefteqn{\scriptstyle \tilde\pi_n}\\
\Bun_{P_{2n}} & \toup{\tilde\nu_{2n}} & \Bunt_{G_{2n}},
\end{array}
$$
where we denoted by $\pi_{n,P}$ the map sending $(L\subset M)\in\Bun_{P_n, Y}$ to $(\phi_*L\subset \phi_*M)\in\Bun_{P_{2n}}$. The map $\tilde\nu_{n,Y}$ sends $(L\subset M)$ to $(M,\cB=\det\RG(Y,L))$ equipped with $\cB^2\,\iso\, \det\RG(Y, M)$, and $\tilde\nu_{2n}$ sends $(L_1\subset M_1)$ to $(\cB_1, M_1)$, where $\cB_1=\det\RG(X, L_1)$.

 Write $_c\Bun_{n, Y}\subset \Bun_{n,Y}$ for the open substack given by $\H^0(Y, L)=0$ for $L\in \Bun_{n,Y}$. Let $_c\cV_Y\to {_c\Bun_{n,Y}}$ be the vector bundle with fibre $\Hom(L, \Omega_Y)$ over $L$. Let $_c\Bun_{P_n, Y}$ be the preimage of $_c\Bun_n$ in $\Bun_{P_n, Y}$. Recall the natural map
$f: {_c\Bun_{P_n, Y}}\to \Sym^2 {_c\cV_Y^*}$ over $_c\Bun_{n,Y}$ and the perverse sheaf $S_{P,\psi}$ on $\Sym^2 {_c\cV_Y^*}$ introduced in \cite[Section~5.2]{Ly4}. 

 Let $_c\Bun_{2n}\subset\Bun_{2n}$ be the open substack classifying $L_1\in \Bun_{2n}$ with $\H^0(X, L_1)=0$. One defines $_c\Bun_{P_{2n}}\subset\Bun_{P_{2n}}$ similarly. Let $_c\cV\to {_c\Bun_{2n}}$ be the vector bundle with fibre $\Hom(L_1, \Omega)$ at $L_1$. We have a similar map $f: {_c\Bun_{P_{2n}}}\to \Sym^2 {_c\cV^*}$ over $_c\Bun_{2n}$, and a perverse sheaf $S_{P,\psi}$ on 
$\Sym^2 {_c\cV^*}$. The following diagram commutes
$$
\begin{array}{ccc}
_c\Bun_{P_n, Y} & \toup{f} & \Sym^2 {_c\cV_Y^*}\\
\downarrow\lefteqn{\scriptstyle \pi_{n,P}} && \downarrow\lefteqn{\scriptstyle \pi_{\cV}}\\
_c\Bun_{P_{2n}} & \toup{f} & \Sym^2 {_c\cV^*},
\end{array}
$$
where the map $\pi_{\cV}$ is induced by the morphism $_c\cV^*_Y\to {_c\cV^*}$ sending $(L, v\in\H^1(Y, L))$ to $(\phi_*L, v\in \H^1(X, \phi_*L))$. Since $\pi_{\cV}^*S_{P,\psi}[\dimrel(\pi_{\cV})]\,\iso\, S_{P,\psi}$ naturally, our claim follows.
\end{proof}


\begin{thebibliography}{99}
\bibitem{A+} Dima Arinkin, Dario Beraldo, Justin Campbell, Lin Chen, Yuchen Fu, Dennis Gaitsgory, Quoc Ho, SergeyLysenko,  Sam  Raskin,  Simon  Riche,  Nick  Rozenblyum,  James  Tao,  David  Yang, Yifei  Zhao, Notes from the winter school on local geometric Langlands, 2018. Available at {\tt https://sites.google.com/site/winterlanglands2018/notes-of-talks}
\bibitem{BD} A. Beilinson, V. Drinfeld, Quantization of HitchinÕs integrable system and Hecke eigen-sheaves, preprint available at {\tt http://www.math.uchicago.edu/$\sim$arinkin/langlands}
\bibitem{BG} A. Braverman, D. Gaitsgory, Geometric Eisenstein series, Inv. Math. 150 (2002), 287 - 384
\bibitem{De} P. Deligne, Le d\'eterminant de la cohomologie, Contemp. Math. 67 (1987), 93 - 177 
\bibitem{FG} E. Frenkel, D. Gaitsgory, Local geometric Langlands correspondence and affine Kac-Moody algebras, In \select{Algebraic geometry and number theory}, vol. 253 of Progr. Math., 69 - 260. Birkh\"auser Boston, Boston, MA, 2006.
\bibitem{FGV} E. Frenkel, D. Gaitsgory, K. Vilonen, On the geometric Langlands conjecture, J. Amer. Math. Soc.  15  (2002),  no. 2, 367 - 417
\bibitem{FH91} W. Fulton, J. Harris, Representation theory, a first course, Springer,  Graduate Texts in Mathematics, vol. 129 (2004)
\bibitem{G} D. Gaitsgory, On a vanishing conjecture appearing in the geometric Langlands correspondence, Ann. of Math., 160 (2004), 617 - 682
\bibitem{G1} D. Gaitsgory, Quantum Langlands Correspondence, arXiv:1601.05279
\bibitem{GL1} D. Gaitsgory, S. Lysenko, Metaplectic Whittaker category and quantum groups : the `small' FLE, arXiv:1903.02279
\bibitem{GL2} D. Gaitsgory, S. Lysenko, Parameters and duality for the metaplectic geometric Langlands theory, available at {\tt http://www.math.harvard.edu/~gaitsgde/GL/twistings.pdf}, version April 28, 2019
\bibitem{GN} D. Gaitsgory, D. Nadler, Spherical varieties and Langlands duality, Moscow Math. J.,10, Nu. 1 (2010), 65 - 137
\bibitem{KS} M. Kamgarpour, T. Schedler, Geometrization of principal series representations of reductive groups,  Annales de l'Institut Fourier, t. 65 (2015) no. 5, p. 2273 - 2330 
\bibitem{LO}
Y. Laszlo, M. Olsson, The six operations for sheaves on Artin stacks II: adic coefficients, Publ. Math. IHES, vol. 107, Nu. 1 (2008), 169 - 210
\bibitem{LL} V. Lafforgue, S. Lysenko, Geometric Weil representation: local field case, Compos. Math. 145 (2009), no. 1, 56 - 88
\bibitem{Laum} G. Laumon, Correspondance de Langlands g\'eom\'etrique pour les corps de fonctions, Duke Math. J., vol. 54, Nu. 2 (1987), 309 - 359
\bibitem{Ly3} S. Lysenko, Local geometrized Rankin-Selberg method for $\GL_n$, Duke Math. J. vol. 111, No. 3 (2002)
\bibitem{Ly2} S. Lysenko, Whittaker and Bessel functors for $\GSp_4$, Ann. Inst. Fourier (Grenoble) 56 (2006), 1505 - 1565
\bibitem{Ly4} S. Lysenko, Moduli of metaplectic bundles on curves and theta-sheaves, Annales ENS, 4eme s\'erie, 39 (2006), 415 - 466
\bibitem{Ly1} S. Lysenko, Geometric Waldspurger periods, Compos. Math. 144 (2008), no. 2, 377 - 438
\bibitem{Ly5} S. Lysenko, Geometric theta-lifting for the dual pair $\SO_{2m}, \Sp_{2n}$, Annales ENS, 4eme s\'erie, t. 44 (2011), 427 - 493
\bibitem{Ly6} S. Lysenko, Geometric Whittaker models and Eisenstein series for $\mathrm{Mp}_2$,  arXiv:1211.1596
\bibitem{Ly8} S. Lysenko, Twisted geometric Langlands correspondence for a torus, Int. Math. Res. Notices (2015) (18): 8680 - 8723
\bibitem{Ly7} S. Lysenko, Twisted Whittaker models for metaplectic groups, GAFA 2017, vol. 27, Issue 2, 289 - 372
\bibitem{Mu} D. Mumford, Prym varieties I, Contribution to analysis (Academic Press, New York, 1974), 325 - 350
\bibitem{SV}  Y. Sakellaridis, A. Venkatesh, Periods and harmonic analysis on spherical varieties, Ast\'erisque 396 (2018)
\bibitem{W} J.-L. Waldspurger, Sur les valeurs de certaines fonctions $L$ automorphes en leur centre de sym\'etrie, Comp. Math. 54 (1985), 173 - 242
\bibitem{W1} J.-L. Waldspurger, Correspondance de Shimura, J. Math. Pures Appl. (9) 59 (1980), no. 1, 1 - 132.
\bibitem{W2} J.-L. Waldspurger, Correspondance de Shimura et  et quaternions, Forum Math. 3 (1991), no. 3, 219 - 307
\bibitem{WTG} W. T. Gan, Representation of metaplectic groups, Fifth International Congress of Chinese Mathematicians. Part 1, 2, 155 - 171, AMS/IP Stud. Adv. Math., 51, pt. 1, 2, Amer. Math. Soc., Providence, RI (2012)
\bibitem{Y} Zh. Yun, Motives with exceptional Galois groups and the inverse Galois problem, Inv. Math. (2014), Vol. 196, Issue 2, 267 - 337 
\end{thebibliography}
\end{document}